\documentclass[12pt]{amsart}

\usepackage[english]{babel}

\usepackage{dirtytalk}

\usepackage{amsmath}

\usepackage[mathx,matha]{mathabx}

\usepackage{amstext}

\usepackage{amssymb}

\usepackage{array}

\usepackage{amsthm}

\usepackage{wasysym}

\usepackage{rotating}

\usepackage{dsfont}

\usepackage[latin1]{inputenc}

\usepackage[T1]{fontenc}

\usepackage{mathtools}

\usepackage{appendix}

\usepackage{xcolor}

\usepackage{ulem}

\usepackage[arrow, matrix, curve]{xy}

\usepackage{verbatim}

\usepackage{graphicx}
\graphicspath{{images/}}

\usepackage{subfigure}

\usepackage{hyperref}

\usepackage{cases}

\usepackage[update,prepend]{epstopdf}

\usepackage{bbold}

\usepackage[a]{esvect}

\usepackage{tikz}

\newtheorem{theorem}{Theorem}[section]

\newtheorem*{theorem*}{Theorem}

\theoremstyle{plain}

\newtheorem*{conjecture*}{Conjecture}

\newtheorem{assumption}[theorem]{Assumption}

\newtheorem{corollary}[theorem]{Corollary}

\newtheorem{proposition}[theorem]{Proposition}

\newtheorem{lemma}[theorem]{Lemma}

\theoremstyle{remark}

\newtheorem{condition}{Condition}

\theoremstyle{definition}

\newtheorem{definition}[theorem]{Definition}

\newtheorem{remark}[theorem]{Remark}

\newtheorem{example}[theorem]{Example}

\def\eps{\varepsilon}

\def\bi{\begin{itemize}}

\def\ei{\end{itemize}}

\newcommand{\R}{\mathbb{R}}

\newcommand{\N}{\mathbb{N}}

\newcommand{\Z}{\mathbb{Z}}

\newcommand{\A}{\mathbb{A}}

\renewcommand{\phi}{\varphi}

\def\d{\text{d}}

\DeclareMathOperator{\vol}{vol}

\DeclareMathOperator{\interior}{int}

\DeclareMathOperator{\conv}{conv}

\DeclareMathOperator{\const}{const.}

\DeclareOldFontCommand{\it}{\normalfont\itshape}{\mathit}

\newcommand{\bspm}{\left(\begin{smallmatrix}}\newcommand{\espm}{\end{smallmatrix}\right)}

\newcommand{\bpm}{\begin{pmatrix}}\newcommand{\epm}{\end{pmatrix}}

\def\bs{\begin{satz}}\def\es{\end{satz}}

\def\blem{\begin{lemma}}\def\elem{\end{lemma}}

\def\bthm{\begin{theorem}}\def\ethm{\end{theorem}}

\def\bcor{\begin{corollary}}\def\ecor{\end{corollary}}

\def\beq{\begin{equation}}\def\eeq{\end{equation}}

\def\beqq{\begin{equation*}}\def\eeqq{\end{equation*}}

\def\bal{\begin{align}}\def\eal{\end{align}}

\def\bpf{\begin{proof}}\def\epf{\end{proof}}

\def\bex{\begin{example}}\def\eex{\end{example}}

\def\brem{\begin{remark}}\def\erem{\end{remark}}

\def\bass{\begin{assumption}}\def\eass{\end{assumption}}

\def\bprop{\begin{proposition}}\def\eprop{\end{proposition}}

\def\bdefi{\begin{definition}}\def\edefi{\end{definition}}

\def\bcond{\begin{condition}}\def\econd{\end{condition}}

\def\bconj{\begin{conjecture*}}\def\econj{\end{conjecture*}}

\DeclareFontEncoding{FMS}{}{}

\DeclareFontSubstitution{FMS}{futm}{m}{n}

\DeclareFontEncoding{FMX}{}{}

\DeclareFontSubstitution{FMX}{futm}{m}{n}

\DeclareSymbolFont{fouriersymbols}{FMS}{futm}{m}{n}

\DeclareSymbolFont{fourierlargesymbols}{FMX}{futm}{m}{n}

\DeclareMathDelimiter{\VERT}{\mathord}{fouriersymbols}{152}{fourierlargesymbols}{147}

\def\bi{\begin{itemize}}

\def\ei{\end{itemize}}

\def\ben{\begin{enumerate}}

\def\een{\end{enumerate}}

\usepackage{enumitem}

\usepackage{tcolorbox}

\newtcolorbox{implementation}[2][]{colframe=blue!75!black,colbacktitle=green!10!white,colback=green!10!white,coltitle=green!75!black,title={#2},fonttitle=\bfseries,#1}

\begin{document}

\title[Viterbo's conjecture for Lagrangian products in $\R^4$]{Viterbo's conjecture for Lagrangian products in $\R^4$ and symplectomorphisms to the Euclidean ball}

\author{Daniel Rudolf}


\date{\today}

\maketitle

\begin{abstract}
We use the generalized Minkowski billiard characterization of the EHZ-capacity of Lagrangian products in order to reprove that the $4$-dimensional Viterbo conjecture holds for the Lagrangian products (any triangle/parallelogram in $\mathbb{R}^2$)$\times$(any convex body in $\mathbb{R}^2$) and extend this fact to the Lagrangian products (any trapezoid in $\mathbb{R}^2$)$\times$(any convex body in $\mathbb{R}^2$). Based on this analysis, we classify equality cases of this version of Viterbo's conjecture and prove that most of them can be proven to be symplectomorphic to Euclidean balls. As a by-product, we prove sharp systolic Minkowski billiard / worm problem inequalities. Furthermore, we discuss the Lagrangian products (any convex quadrilateral in $\mathbb{R}^2$)$\times$(any convex body in $\mathbb{R}^2$) for which we show that the truth of Viterbo's conjecture would follow from the positive solution of a challenging Euclidean covering problem. Finally, we show that the flow associated to equality cases of Viterbo's conjecture for Lagrangian products in $\mathbb{R}^4$--which turn out to be convex polytopes--is not Zoll in general, but that a weaker Zoll property, namely, that every characteristic almost everywhere away from lower-dimensional faces is closed and action-minimizing, does apply.
\end{abstract}

\section{Introduction and main results}\label{Sec:equacaIntromainresults}

Viterbo's conjecture is an isoperimetric-type conjecture for symplectic capacities of convex bodies, which was raised by Viterbo in \cite{Viterbo2000}. It says that for any symplectic capacity $c$ and any convex body $C$ in the standard symplectic phase space $\R^{2n}=\R^{2n}(x,y)$ one has
\beq \frac{c(C)}{c\left(B_1^{2n}(0)\right)} \leq \left(\frac{\vol(C)}{\vol\left(B_1^{2n}(0)\right)}\right)^{\frac{1}{n}},\label{eq:Viterboconjecture1}\eeq
where by $B_1^{2n}(0)$ we denote the $2n$-dimensional unit ball of $\R^{2n}$, with equality holding if and only if $C$ is symplectomorphic to a ball. In other words, Viterbo's conjecture states that among the convex bodies in $\R^{2n}$ with given volume, the Euclidean ball has the maximal symplectic capacity. Plugging the known volume for the ball into \eqref{eq:Viterboconjecture1}, we get Viterbo's conjecture in the form of a systolic ratio:
\beq c(C)\leq \left(n! \vol(C)\right)^{\frac{1}{n}} \; \; \Leftrightarrow \; \; \vol(C)\geq \frac{c(C)^n}{n!}.\label{eq:Viterboconjecture2}\eeq

Here, symplectic capacities $c$ are characterized as functions mapping convex bodies in $\R^{2n}(x,y)$ to values in $[0,\infty]$ while they satisfy the following axioms: monotonicity, i.e., for convex bodies $C_1$ and $C_2$ in $\R^{2n}$ with $C_1 \subseteq C_2$ one has $c(C_1) \leq c(C_2)$; $2$-homogeinity, i.e., for $\lambda\neq 0$ and convex body $C$ in $\R^{2n}$ one has $c(\lambda C)=\lambda ^2 c(C)$; and nontriviality, i.e., one has
\beqq c\left(B_1^{2n}(0)\right)=c\left(Z^{2n}_1(0)\right)=\pi,\eeqq
where by $Z^{2n}_1(0)$ we denote the symplectic cylinder
\beqq B_1^2(0)\times \R^{2n-2} = \left\{z=(x_1,...,x_n,y_1,...y_n)\in\R^{2n}: x_1^2+y_1^2\leq 1\right\}.\eeqq

Before we present our results concerning this conjecture, let us first briefly review the current status of research.

Trivially, Viterbo's conjecture holds for $n=1$ in which any symplectic capacity must agree with the area. Viterbo himself proved in \cite{Viterbo2000} his conjecture up to a constant that depends linearly on the dimension using the classical John ellipsoid theorem (see \cite{John1948}). In \cite{ArtsteinAvidanMilmanOstrover2008}, Viterbo's conjecture has been proven up to a dimension-independent constant customizing methods and techniques from asymptotic geometric analysis and adjusting them to the symplectic context.

Besides this, the conjecture is known to hold for certain classes of convex bodies, including ellipsoids and convex Reinhardt domains as has been proven in \cite{Hermann1998}, and, in the case $n=2$ and for the E(keland-)H(ofer-)Z(ehnder)-capacity, for all convex domains which are close enough to a ball as has been proven in \cite{ABHS2018}. Recently, in \cite{AbboBene2020}, this result has been extended to all $n\in\N$ and all symplectic capacities.

In \cite{Balitskiy2018}, it has been proven that the Lagrangian product of a permutohedron and a simplex (properly related to each other) delivers equality in Viterbo's conjecture for the EHZ-capacity. Furthermore, therein, Viterbo's conjecture has been proven for some special cases and interpreted as isoperimetric-like inequalities for billiard trajectories.

We note that Viterbo's conjecture trivially holds for the Gromov width which for a convex body $C\subset\R^{2n}$ can be defined by
\beqq w_G(C):=\sup\left\{\pi r^2 : B^{2n}_r(0)\text{ embeds symplectically into }C\right\}.\eeqq
From this perspective, the long standing open question (see \cite{Hermann1998}, \cite{Hofer1989}, \cite{Viterbo2000}, or \cite{Ostrover2014}) whether all symplectic capacities coincide on the class of convex domains in $\R^{2n}$ becomes important. An affirmative answer to this question would immediately imply Viterbo's conjecture. This immediately implies the relevance of investigating equality cases of Viterbo's conjecture: If we could find an equality case of Viterbo's conjecture which is not symplectomorphic to a Euclidean ball, then this would be a counterexample to this question whether all symplectic capacities coincide on convex domains. This follows directly from the definition of the Gromov width capacity (see also the discussion around Question 5.1 in \cite{Ostrover2014}).

For the special case of the EHZ-capacity $c_{EHZ}$ of Lagrangian products of the form $K\times K^\circ$, where $K\subset\R^n$ is a centrally symmetric convex body and $K^\circ$ its polar body, it has been shown in \cite{ArtKarOst2013} that $c_{EHZ}(K\times K^\circ)=4$ and hence Viterbo's conjecture coincides with the well-known Mahler conjecture (see \cite{Mahler1939}) from convex geometry:
\beqq \nu(K)=\vol(K)\vol(K^\circ)\geq \frac{4^n}{n!},\eeqq
where $\nu(K)$ is called the Mahler volume of $K$.

In what follows, we are mainly concerned with Viterbo's conjecture for the EHZ-capacity of Lagrangian products $K\times T$ in $\R^4$, where $K,T\subset\R^2$ are convex bodies. Viterbo's conjecture for this class of convex sets is motivated, on the one hand, by the connection to Mahler's conjecture, on the other hand, by the fact that it becomes a systolic inequality for Minkowski billiards. Here, the systolic question for Minkowski billiards consists in bounding the minimal length of closed Minkowski billiard trajectories by the volume (here: area) of the billiard table. We note especially that Lagrangian products $K\times T \subset \R^2 \times \R^2$ allow the usage of the following viable connection: if $\phi:\R^2\rightarrow\R^2$ is a linear isomorphism, then $\phi\times \left(\phi^T\right)^{-1}$ is a symplectomorphism from $K\times T$ to $\phi(K)\times\left(\phi^T\right)^{-1}(T)$ (and this is also true in higher dimensions).

Now, let us make our setting more precise. Let $K,T\subset\R^2$ be two convex bodies, i.e., compact convex sets with the origin in their interiors. The EHZ-capacity $c_{EHZ}$ of the convex Lagrangian product $K\times T$ can be defined\footnote{This definition is the outcome of a historically grown study of symplectic capacities. Traced back--recalling that $c_{EHZ}$ in its present form is the generalization (by K\"{u}nzle in \cite{Kuenzle1996}) of a symplectic capacity after applying the dual action functional introduced by Clarke in \cite{Clarke1979}--, the EHZ-capacity denotes the coincidence of the Ekeland-Hofer- and Hofer-Zehnder-capacities, originally constructed in \cite{EkHo1989} and \cite{HoferZehnder1990}, respectively.} by
\beqq c_{EHZ}(K\times T)=\min\{\A(x):x \text{ closed characteristic on }\partial (K\times T)\},\eeqq
where a closed characteristic\footnote{Sometimes also called \textit{generalized} characteristics.} on $\partial (K\times T)$ is an absolutely continuous loop $x$ in $\R^{4}$ satisfying
\beqq \begin{cases}\dot{x}(t)\in \R^+J\partial \mu_{K\times T}(x(t))\quad \text{a.e.}, \\ \mu_{K\times T}(x(t))=1\quad\forall t\in \R/\Z,\end{cases}\eeqq
where $J$ is the symplectic matrix $\begin{pmatrix} 0 & \mathbb{1} \\ -\mathbb{1} & 0 \end{pmatrix}$, $\partial$ the subdifferential-operator, and $\mu$ the Minkowski functional given by
\beqq \mu_{K\times T}(x)=\min\left\{t\geq 0 : x\in t(K\times T)\right\}, \quad x\in\R^{4},\eeqq
and where by $\A(x)$ we denote the loop's action defined by
\beqq \A(x)= -\frac{1}{2}\int_{\R/\Z} \langle J\dot{x}(t),x(t)\rangle \;dt.\eeqq
Within this setting, \eqref{eq:Viterboconjecture2} becomes
\beq \vol(K\times T) \geq \frac{c_{EHZ}(K\times T)^2}{2}.\label{eq:Viterboconjecture3}\eeq

In order to state our first main result concerning \eqref{eq:Viterboconjecture3}, we briefly introduce what we mean by \textit{trapezoids} in $\R^2$: these are the quadrilaterals in $\R^2$ with at least one pair of sides parallel, where we also include triangles as degenerate trapezoids. 

Then, our first main result answers the question, for which Lagrangian products \eqref{eq:Viterboconjecture3} can be proven to be true.

\bthm\label{Thm:answer1}
Let $Q$ be any trapezoid in $\R^2$. Then, Viterbo's conjecture is true for all Lagrangian products
\beq Q \times T, \label{eq:answer1}\eeq
where $T$ is any convex body in $\R^2$.
\ethm

We note that the Lagrangian products $Q\times T$ in \eqref{eq:answer1} include the Lagrangian products
\beq \Delta \times T \quad \text{ and } \quad \square \times T,\label{eq:balitskiyconfigurations}\eeq
where $\Delta$ is any triangle, $\square$ any parallelogram, and $T$ any convex body in $\R^2$, respectively, for which Balitskiy in \cite[Theorems 5.1 and 5.2]{Balitskiy2018} already proved that Viterbo's conjecture is true. In this regard, we remark that all triangles as well as all parallelograms are equivalent modulo affine transformations, respectively. Trapezoids, instead, form a one-dimensional space inside the two-dimensional space of convex quadrilaterals, where convex quadrilaterals are equivalent to \textit{diamonds} by which we mean quadrilaterals whose diagonals intersect perpendicularly and have the same length.

We will prove Theorem \ref{Thm:answer1} by rigorously using the generalized Minkowski billiard characterization of the EHZ-capacity whose generalization to general convex bodies $K\times T\subset \R^n \times \R^n$--without requiring additional conditions on $K$ and $T$--we proved in \cite{Rudolf2022b}\footnote{We remark that the Minkowski billiard characterization of the EHZ-capacity was made explicit for the first time in \cite{ArtOst2012}. They showed this characterization only under the assumption of smoothness and strict convexity of both $K$ and $T$. For further details we refer to the discussion beyond Theorem 1.1 in \cite{Rudolf2022b}.}. This will help us especially in order to classify the equality cases of Viterbo's conjecture for these configurations in a subsequent step. In this context, we remark that Balitskiy's proof in \cite{Balitskiy2018} for the Lagrangian configurations \eqref{eq:balitskiyconfigurations} was partly based on a somewhat different technique involving the representation of the permutohedron in terms of the Voronoi cell of a certain lattice. It is worth mentioning that his technique allowed him to prove this statement for the higher dimensional generalizations, while--based on his technique--it seems harder to handle action-minimizing characteristics through non-smooth boundary points (where our generalized Minkowski billiard characterization of the EHZ-capacity comes into play).

\begin{figure}[h!]
\centering
\def\svgwidth{360pt}
\begingroup%
  \makeatletter%
  \providecommand\color[2][]{%
    \errmessage{(Inkscape) Color is used for the text in Inkscape, but the package 'color.sty' is not loaded}%
    \renewcommand\color[2][]{}%
  }%
  \providecommand\transparent[1]{%
    \errmessage{(Inkscape) Transparency is used (non-zero) for the text in Inkscape, but the package 'transparent.sty' is not loaded}%
    \renewcommand\transparent[1]{}%
  }%
  \providecommand\rotatebox[2]{#2}%
  \newcommand*\fsize{\dimexpr\f@size pt\relax}%
  \newcommand*\lineheight[1]{\fontsize{\fsize}{#1\fsize}\selectfont}%
  \ifx\svgwidth\undefined%
    \setlength{\unitlength}{280.70146338bp}%
    \ifx\svgscale\undefined%
      \relax%
    \else%
      \setlength{\unitlength}{\unitlength * \real{\svgscale}}%
    \fi%
  \else%
    \setlength{\unitlength}{\svgwidth}%
  \fi%
  \global\let\svgwidth\undefined%
  \global\let\svgscale\undefined%
  \makeatother%
  \begin{picture}(1,1.11254328)%
    \lineheight{1}%
    \setlength\tabcolsep{0pt}%
    \put(0,0){\includegraphics[width=\unitlength,page=1]{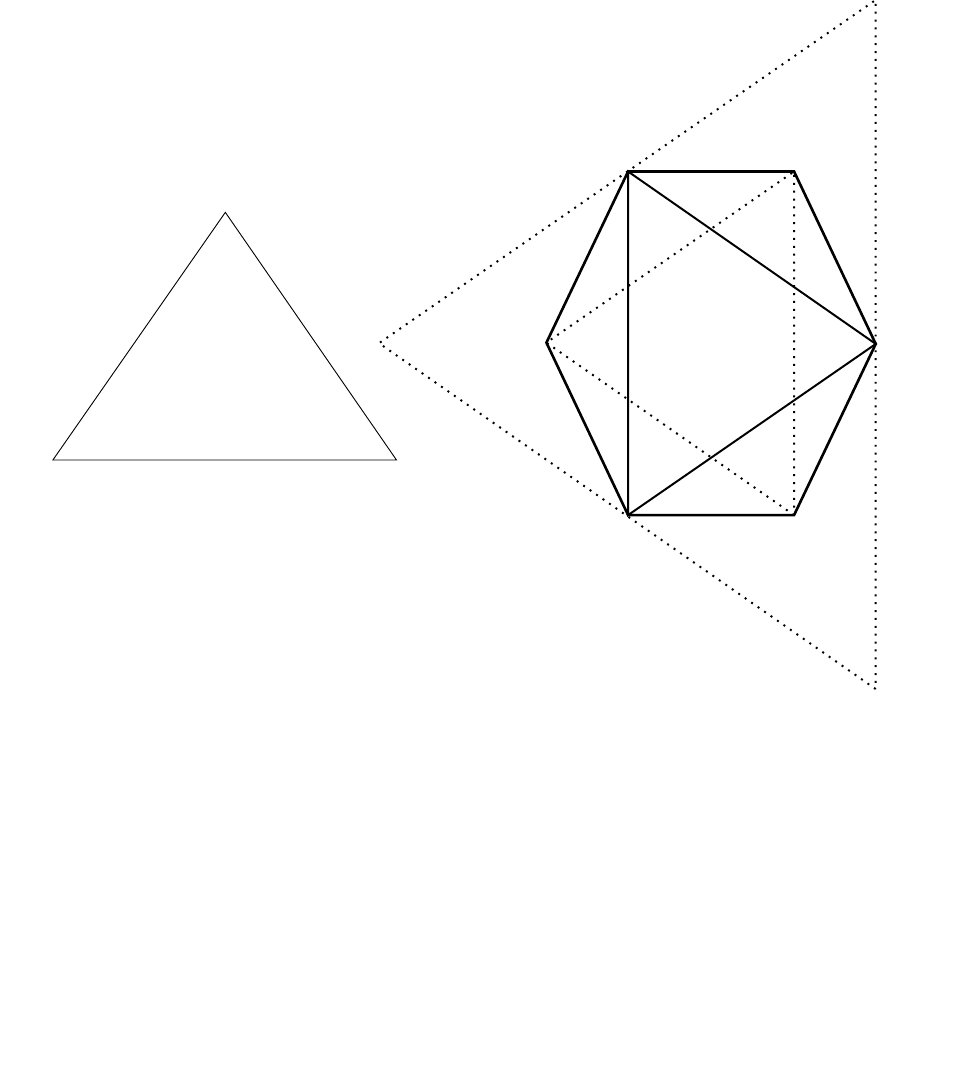}}%
    \put(0.10647673,0.82038746){\color[rgb]{0,0,0}\makebox(0,0)[lt]{\lineheight{1.25}\smash{\begin{tabular}[t]{l}$\Delta$\end{tabular}}}}%
    \put(0.79006621,0.9917556){\color[rgb]{0,0,0}\makebox(0,0)[lt]{\lineheight{1.25}\smash{\begin{tabular}[t]{l}$-2J\Delta$\end{tabular}}}}%
    \put(0.51991788,0.94107564){\color[rgb]{0,0,0}\makebox(0,0)[lt]{\lineheight{1.25}\smash{\begin{tabular}[t]{l}$T$\end{tabular}}}}%
    \put(0.92074631,0.63436181){\color[rgb]{0,0,0}\makebox(0,0)[lt]{\lineheight{1.25}\smash{\begin{tabular}[t]{l}$-J\Delta$\end{tabular}}}}%
    \put(0.92582311,0.87319457){\color[rgb]{0,0,0}\makebox(0,0)[lt]{\lineheight{1.25}\smash{\begin{tabular}[t]{l}$J\Delta$\end{tabular}}}}%
    \put(0,0){\includegraphics[width=\unitlength,page=2]{convexhull.pdf}}%
  \end{picture}%
\endgroup%

\caption[Lagrangian products $\Delta\times T$ which are equality cases of Viterbo's conjecture.]{We assume that $\Delta$'s centroid is the origin. Then, the figure illustrates some of the Lagrangian products $\Delta \times T$ which are equality cases of Viterbo's conjecture. Here, $T$ is one of the volume-minimizing convex hulls of $J\Delta$ and $-J\Delta + t$.}
\label{img:convexhull}
\end{figure}

Our second main result then answers the question, which of the above Lagrangian products for which Viterbo's conjecture holds are equality cases. In order to state this, we need the following preparations: If $\Delta$ is any triangle in $\R^2$, then we define $\mathcal{T}_\Delta$ as the set of volume-minimizing convex hulls
\beq \conv\{J\Delta,-J\Delta +t\}\label{eq:convexhull}\eeq
over all $t\in\R^2$. In other words, $\mathcal{T}_\Delta$ equals
\beqq \left\{\conv\{J\Delta,-J\Delta + t\}:t \text{ minimizes }\vol\left(\conv\{J\Delta,-J\Delta+\widetilde{t}\}\right) \text{ over all }\widetilde{t}\in\R^2\right\}.\eeqq
Under the assumption that $\Delta$'s centroid is the origin (which comes without loss of generality), we remark that $T$ is in $\mathcal{T}_\Delta$ if and only if $t$ in \eqref{eq:convexhull} is in $-J\Delta$ or, equivalently, if $-J\Delta + t$ is subset of $-2J\Delta$ (see Figure \ref{img:convexhull}).
Clearly, $T$ is a hexagon when $-J\Delta+t$ is subset of the interior of $-2J\Delta$ (which means that $t$ in \eqref{eq:convexhull} is in $-J\mathring{\Delta}$) and $T$ is a parallelogram if $-J\Delta + t$ is a subset of $-2J\Delta$ such that it touches its boundary (which means that $t$ in \eqref{eq:convexhull} is in $-J\partial \Delta$). Because of this, we understand $\mathcal{T}_\Delta$ as disjoint union
\beqq \mathcal{T}_\Delta = \mathcal{T}_{\Delta,\varhexagon} \cup \mathcal{T}_{\Delta,\square},\eeqq
where $\mathcal{T}_{\Delta,\varhexagon}$ and $\mathcal{T}_{\Delta,\square}$ can be written as
\beqq \mathcal{T}_{\Delta,\varhexagon} = \left\{\conv\{J\Delta,-J\Delta + t\}:t\in -J\mathring{\Delta}\right\}\eeqq
and
\beqq \mathcal{T}_{\Delta,\square} = \left\{\conv\{J\Delta,-J\Delta + t\}:t\in -J\partial\Delta\right\}.\eeqq
Furthermore, by
\beqq \diamondsuit(a_1,a_2):=\conv\left\{a_1\times [0,1],[0,1]\times a_2\right\},\; a_1,a_2\in[0,1],\eeqq
we denote \textit{diamonds in basic form} (meaning that the diagonals of the diamond are parallel to the horizontal and vertical axes of $\R^2$). For $a_1\in\{0,1\}$ or $a_2\in\{0,1\}$, one gets triangles; for $a_1=a_2=\frac{1}{2}$, parallelograms; and for $a_1\in\{a_2, 1-a_2\}$, trapezoids. A diamond which is not a trapezoid is characterized by diagonals which do not divide proportionally.

Then, our second main result reads:

\bthm\label{Thm:answer2}
\begin{itemize}
\item[(i)] Let $\Delta$ be any triangle in $\R^2$. Then, all Lagrangian products
\beq \Delta \times T\; \text{ with }\; T\in\mathcal{T}_\Delta\label{eq:answer2}\eeq
are equality cases of Viterbo's conjecture. These are the only ones within this configuration (up to scaling and translation in the second component).
\item[(ii)] Let $\square$ be any square in $\R^2$ whose sides are parallel to the horizontal and vertical axes of $\R^2$. Then, all Lagrangian products
\beq \square \times \diamondsuit(a_1,a_2), \quad a_1,a_2\in[0,1],\label{eq:answer2ii}\eeq
are equality cases of Viterbo's conjecture. These are the only ones within this configuration (up to scaling and translation in the second component).
\item[(iii)] Let $Q$ be any trapezoid in $\R^2$ which is neither a triangle nor a parallelogram. Then, there is a parallelogram $P$ in $\R^2$--which we will specify later--such that
\beq Q\times P\label{eq:answer2iii}\eeq
is an equality case of Viterbo's conjecture. This is the only one within this configuration (up to scaling and translation in the second component).
\end{itemize}
\ethm

We note that by using suitable affine transformations, the square-configurations \eqref{eq:answer2ii} can be easily lifted to parallelogram-configurations. Furthermore, we remark that it will turn out that the trapezoid-case (iii) is included within the parallelogram-case (ii) (by finding a suitable symplectomorphism and interchanging $\R^2(x)$ and $\R^2(y)$ whose concatenation maps configurations of the form \eqref{eq:answer2iii} onto configurations of the form \eqref{eq:answer2ii}).

Our third main result answers the question, the interiors of which equality cases presented in Theorem \ref{Thm:answer2} can be proven to be symplectomorphic to a Euclidean ball:

\bthm\label{Thm:answer3}
\begin{itemize}
\item[(i)] Let $\Delta$ be any triangle in $\R^2$. Then, the interiors of the Lagrangian products
\beqq \Delta \times T \; \text{ with } \; T\in \mathcal{T}_{\Delta,\square}\eeqq
are symplectomorphic to a Euclidean ball.
\item[(ii)] Let $\square$ be any square in $\R^2$ whose sides are parallel to the horizontal and vertical axes of $\R^2$. Then, the interiors of the Lagrangian products
\beqq \square \times \diamondsuit(a_1,a_2), \quad a_1,a_2\in[0,1],\eeqq
are symplectomorphic to a Euclidean ball.
\item[(iii)] Let $Q$ be any trapezoid in $\R^2$ which is neither a triangle nor a parallelogram. If $P$ is the parallelogram from Theorem \ref{Thm:answer2}(iii), then the interior of
\beqq Q\times P\eeqq
is symplectomorphic to a Euclidean ball.
\end{itemize}
\ethm

We remark that, so far, it is not clear whether the interiors of the equality cases from Theorem \ref{Thm:answer2}(i)
\beqq \Delta \times T \; \text{ with } \; T\in \mathcal{T}_{\Delta,\varhexagon}\eeqq
are symplectomorphic to a Euclidean ball.\footnote{We note that Ostrover, Ramos, and Sepe recently announced being able to prove that the Lagrangian product of the interior of an equilateral triangle and a regular hexagon is indeed also symplectomorphic to a Euclidean ball. This would immediately imply that all equality cases known so far of Viterbo's conjecture for Lagrangian products in $\R^4$ are symplectomorphic to Euclidean balls.}

We remark that Theorems \ref{Thm:answer2}(iii) and \ref{Thm:answer3}(iii) also hold for $Q$ replaced by any convex quadrilateral whose diagonals do not divide proportionally (together with appropriately chosen parallelogram). But, whether Theorem \ref{Thm:answer1} holds for $Q$ replaced by any convex quadrilateral is not clear so far. However, in Section \ref{Sec:quadrilateral}, we will formulate a conjecture for a covering/container problem which implies Viterbo's conjecture for these sets.

An immediate consequence of Theorems \ref{Thm:answer1} and \ref{Thm:answer2} is the following corollary which provides the implication of the above results concerning Viterbo's conjecture for (sharp) Minkowski billiard (see \eqref{eq:CorMinkBillIneq}) and Minkowski worm (see \eqref{eq:CorMinkWormIneq}) inequalities. Against this background, the following corollary answers the above mentioned systolic question for Minkowski billiard trajectories in the case of Minkowski geometries whose unit balls are trapezoids.

Before we state the corollary, we first have to introduce the corresponding notations: We start by recalling that for convex bodies $K,T\subset\R^n$, we say that a closed polygonal curve\footnote{For the sake of simplicity, whenever we talk of the vertices $q_1,...,q_m$ of a closed polygonal curve, we assume that they satisfy $q_j\neq q_{j+1}$ and $q_j$ is not contained in the line segment connecting $q_{j-1}$ and $q_{j+1}$ for all $j\in\{1,...,m\}$. Furthermore, whenever we settle indices $1,...,m$, then the indices in $\Z$ will be considered as indices modulo $m$.\label{foot:polygonalline}} with vertices $q_1,...,q_m$, $m\geq 2$, on the boundary of $K$ is a closed $(K,T)$-Minkowski billiard trajectory if there are points $p_1,...,p_m$ on $\partial T$ such that
\beq \begin{cases}q_{j+1}-q_j\in N_T(p_j),\\ p_{j+1}-p_j\in -N_K(q_{j+1})\end{cases}\label{eq:Minkreflectionrulepre}\eeq
is satisfied for all $j\in\{1,...,m\}$. Further, we recall that for convex body $K\subset\R^n$, we denote by $F^{cp}_{n+1}(K)$ the set of closed polygonal curves with at most $n+1$ vertices which cannot be translated into $\mathring{K}$. We denote by $M(K,T)$ the set of closed $(K,T)$-Minkowski billiard trajectories and by $M_{n+1}(K,T)$ the set of members of $M(K,T)$ which have at most $n+1$ bouncing points. We further recall that for convex body $T\in\R^2$ a Minkowski worm problem is to find minimizers of
\beqq \min_{K\in A(T,\alpha)} \vol(K),\eeqq
where by $A(T,\alpha)$ we denote the set of convex bodies in $\R^n$ that cover a translate of every closed curve of $\ell_T$-length $\alpha$. Here, the $\ell_T$-length of a closed curve is its arc length measured by the Minkowski functional $\mu_{T^\circ}$.

\bcor\label{Cor:sharpinequalities}
Let $Q$ be any trapezoid in $\R^2$. Then, we have
\beq \min_{q\in M(K,Q)}\ell_Q^2(q)=\min_{q\in M_3(K,Q)}\ell_Q^2(q) \leq 2\vol(K)\label{eq:CorMinkBillIneq}\eeq
for all convex bodies $K\subset\R^2$ and
\beq \min_{K\in A(T,\sqrt{2})}\vol(K)\geq 1.\label{eq:CorMinkWormIneq}\eeq
Both inequalities are sharp.
\ecor

Closed symplectic balls in $\R^4$ can be characterized as the unions of smooth starshaped domains in $\R^4$ which are Zoll, meaning that all the characteristics on the boundary are closed and have the same action. The Lagrangian products appearing as equality cases in Theorem \ref{Thm:answer2} are polytopes and Theorem \ref{Thm:answer3} says that at least in some cases, their interiors are symplectomorphic to an open ball. Since their boundary is not smooth, characteristics have to be understood in the generalized sense mentioned above, but it makes sense to ask whether the Zoll property still holds. In what follows, we will show that this is not the case in general, but the subsequent theorem will show that at least a weaker version of this property holds for all equality cases presented in Theorem \ref{Thm:answer2}: All orbits which do not run through vertices ($0$-faces) or on $1$-or $2$-faces of the Lagrangian products are closed and have the same action. Here, we note that the $3$-faces of a Lagrangian product $K\times T \subset\R^2 \times \R^2$ are given by the products $K\times$($1$-face of $T$) and ($1$-face of $K$)$\times T$, the $2$-faces by the products $K\times$($0$-face of $T$), ($0$-face of $K$)$\times T$, and ($1$-face of $K$)$\times$($1$-face of $T$), the $1$-faces by the products ($1$-face of $K$)$\times$($0$-face of $T$) and ($0$-face of $K$)$\times$($1$-face of $T$), and the $0$-faces by the products ($0$-face of $K$)$\times$($0$-face of $T$). Furthermore, one has the density of the union of paths of the Minkowski billiard and their dual trajectories:

\bthm\label{Thm:Zollproperty}
All equality cases $K\times T$ in $\R^4$ presented in Theorem \ref{Thm:answer2} satisfy the following properties:
\begin{itemize}
\item[(i)] every regular $(K,T)$-Minkowski billiard trajectory, whose dual billiard trajectory in $T$ is regular as well, is closed, has exactly\footnote{We consider only \textit{simple} billiard trajectories, i.e., trajectories which do not pass the same path multiply.} $4$ bouncing points, and is an $\ell_T$-minimizer;
\item[(ii)] the union of paths of the $(K,T)$-Minkowski billiard trajectories from (i) is dense on $K$;
\item[(iii)] the union of paths of the dual billiard trajectories associated to the $(K,T)$-Minkowski billiard trajectories from (ii) is dense on $T$.
\end{itemize}
\ethm

We remark that the statement of Theorem \ref{Thm:Zollproperty} has already been shown for two special Lagrangian products: for $\Delta\times H$, where $\Delta$ is the equilateral triangle and $H$ a properly related regular hexagon; and for $\square \times \diamondsuit$, where $\square$ is the square (centred at the origin) and $\diamondsuit$ its polar. This has been shown in \cite{Balitskiy2016}.

Transferred to the original task of finding the action-minimizing closed characteristic on the boundary of Lagrangian products $K\times T$, Theorem \ref{Thm:Zollproperty} implies the subsequent corollary. For this transfer we use what we have shown within the proofs of Theorems 1.1 and 2.1 in \cite{Rudolf2022a}: Since both $K$ and $T$ are not strictly convex, in general, not every closed characteristic on $\partial (K\times T)$ can be considered as closed $(K,T)$-Minkowski billiard trajectory. The other way around, however, every closed $(K,T)$-Minkowski billiard trajectory is a closed characteristic on $\partial (K\times T)$, and, additionally, the aforementioned Theorem 1.1 in \cite{Rudolf2022a} guarantees that every $\ell_T$-minimizing $(K,T)$-Minkowski billiard trajectory is an action-minimizing closed characteristic on $\partial (K\times T)$.

\bcor\label{Cor:Zollproperty}
All equality cases $K\times T$ in $\R^4$ pesented in Theorem \ref{Thm:answer2} satisfy the following properties:
\begin{itemize}
\item[(i)] every regular characteristic on $\partial (K\times T)$, where by regularity we mean that the characteristic's trajectory runs on the interiors of the facets almost everywhere, is closed, runs over exactly $8$ facets, and minimizes the action;
\item[(ii)] the union of these trajectories of the characteristics from (i) is dense on $\partial (K\times T)$.
\end{itemize}
\ecor

We note that the statement of Theorem \ref{Thm:Zollproperty} is not true for non-regular $(K,T)$-Minkowski billiard trajectories. We will give examples of closed non-regular $(K,T)$-Minkowski billiard trajectories, where $K\times T$ belongs to the equality cases presented in Theorem \ref{Thm:answer2}, whose $\ell_T$-length is not the minimal one. This means that there are non-regular characteristics on $\partial (K\times T)$ whose action is not the minimal one.

We conclude by summarizing our aforementioned results in Table \ref{table:ResultsCollection}. Therein, we have considered the Lagrangian configurations against the background of their basic form representation, namely, against the background that every triangle/parallelogram/trapezoid/qua\-drilateral can be traced back (by applying a suitable affine transformation) to a diamond $\diamondsuit(a_1,a_2)$ for appropriately chosen $a_1,a_2\in [0,1]$.

\begin{table}
\centering
\caption[Overview of the results concerning Viterbo's conjecture against the background of the basic form representation.]{Overview of the results concerning Viterbo's conjecture against the background of the basic form representation. Concerning the meaning of the respective columns: \say{V}: Viterbo's conjecture true?; \say{A}: All equality cases known?; \say{S}: The interiors of all equality cases symplectomorphic to a Euclidean ball?; \say{Z}: All equality cases satisfy the weak Zoll-property? Note that by $R_{\frac{\pi}{4}}$ we denote the rotation by angle $\frac{\pi}{4}$.}\label{table:ResultsCollection}
\begin{tabular}[h]{c||c|c|c|c|c}
Lagrangian product  & V & equality cases for $T=$ & A & S & Z \\
\hline
\hline
$\diamondsuit(a_1,a_2)\times T$  & $\checkmark$ & $T\in\mathcal{T}_{\diamondsuit(a_1,a_2),\square}$ & $\checkmark$ & $\checkmark$ & $\checkmark$ \\
$a_1\in\{0,1\}$ or  & &  &  &  &  \\
$a_2\in\{0,1\}$  & & $T\in\mathcal{T}_{\diamondsuit(a_1,a_2),\varhexagon}$ & $\checkmark$ & ? & $\checkmark$ \\
(triangles) & &  &  &  &  \\
\hline
$\diamondsuit(a_1,a_2)\times T$  & $\checkmark$ & $R_{\frac{\pi}{4}}\diamondsuit(a_1,a_2)$ & $\checkmark$ & $\checkmark$ & $\checkmark$ \\
$a_1=a_2=\frac{1}{2}$  & & for $a_1,a_2\in[0,1]$ & & & \\
(parallelograms) & & & & & \\
\hline
$\diamondsuit(a_1,a_2)\times T$  & $\checkmark$ & $\square = R_{\tfrac{\pi}{4}}\diamondsuit(\tfrac{1}{2},\tfrac{1}{2})$ & $\checkmark$ & $\checkmark$ & $\checkmark$ \\
$a_1\in\{a_2,1-a_2\}$  & & for $a_1\notin\{0,1\}$ (or $a_2\notin\{0,1\}$) & & &  \\
(trapezoids) & & & & & \\
\hline
$\diamondsuit(a_1,a_2)\times T$  & ? & $\square = R_{\tfrac{\pi}{4}}\diamondsuit(\tfrac{1}{2},\tfrac{1}{2})$ & ? & $\checkmark$ & $\checkmark$ \\ 
$a_1,a_2\in [0,1]$  & & for $a_1\notin\{a_2,1-a_2\}$ & & &  \\
(quadrilaterals) & & & & & \\
\end{tabular}
\end{table}

Let us briefly discuss the structure of this paper: In Section \ref{Sec:equacapreliminaries}, we mainly recall some necessary preliminaries. Then, in Section \ref{Sec:Symplectomorphismstoball}, we rigorously show how to construct symplectomorphisms to the Euclidean ball which will be helpful for the further course, where we split the proofs of Theorems \ref{Thm:answer1}, \ref{Thm:answer2}, and \ref{Thm:answer3} into the following three cases: the triangle- (Section \ref{Sec:triangle}), parallelogram- (Section \ref{Sec:parallelogram}), and convex-quadrilateral-case (Section \ref{Sec:quadrilateral}). In Section \ref{Sec:equacaCorollaries}, we prove Corollary \ref{Cor:sharpinequalities}, and finally in Section \ref{Sec:Zollpropertys}, we prove Theorem \ref{Thm:Zollproperty}.

\section{Preliminaries}\label{Sec:equacapreliminaries}


We recall the generalized Minkowski billiard characterization of the EHZ-capacity of convex Lagrangian products:

\bthm[$4$-dimensional version of Theorem 1.1 in \cite{Rudolf2022a}]\label{Thm:Chap7relationship}
Let $K,T\subset\R^2$ be convex bodies. Then, we have
\beqq c_{EHZ}(K\times T) = \min_{q\in F_{3}^{cp}(K)}\ell_T(q) = \min_{p\in F_{3}^{cp}(T)}\ell_K(p) = \min_{q \in M_{3}(K,T)} \ell_T(q).\eeqq
Moreover, we have
\beqq \min_{q\in F_{j}^{cp}(K)}\ell_T(q) = \min_{q \in M_{j}(K,T)} \ell_T(q)\quad \forall j\in\{2,3\}.\eeqq
\ethm

Furthermore, we need the following proposition:

\bprop[Proposition 5.3 in \cite{Rudolf2022b}]\label{Prop:Chap7characteristicgeodesic}
Let $C\subset\R^{2n}$ be a convex body. Let $x$ be any closed characteristic on $\partial C$. Then, the action of $x$ equals its $\ell_{\frac{JC}{2}}$-length:
\beqq \A(x)=\ell_{\frac{JC}{2}}(x).\eeqq
\eprop

We remark that beyond the utility of this proposition relevant to this paper, it implies a noteworthy connection between closed characteristics and closed Finsler geodesics: Every closed characteristic on $\partial C$ can be interpreted as a Finsler geodesic with respect to the Finsler metric determined by $\mu_{2JC^\circ}$ and which is parametrized by arc length. This raises a number of questions; for example, which Finsler geodesics are closed characteristics (we note that, usually, there are more geodesics than those which, by the least action principle and Proposition \ref{Prop:Chap7characteristicgeodesic}, can be associated to closed characteristics) and whether the length-minimizing Finsler geodesics are of this kind. Following this line of thought, would lead to the question whether it is possible to deduce Viterbo's conjecture from systolic inequalities for Finsler geodesics. However, we leave these questions for further research.

\section{Symplectomorphisms to the Euclidean ball in $\R^4$}\label{Sec:Symplectomorphismstoball}

We consider the Lagrangian splitting
\beq \R^4(z) = \R^2(x)\times \R^2(y)\label{eq:Lagrangiansplitting}\eeq
and define for $a\in\R_{>0}$ and $a_1,a_2\in[0,a]$
\beqq \diamondsuit(a,a_1,a_2)=\interior\left(\conv\{[0,a]\times a_2,a_1\times [0,a]\}\right)\subset\R^2(x),\eeqq
the open square
\beqq \square(1)=\left\{(y_1,y_2)\in\R^2:0<y_1,y_2<1\right\}\subset\R^2(y)\eeqq
and the open ball
\beqq B^4_{\sqrt{\frac{a}{\pi}}}=\left\{z=(z_1,z_2)\in\R^4(z) : \pi|z_1|^2+\pi|z_2|^2 <a \right\}\subset\R^4,\eeqq
with $z_i=(x_i,y_i),\; i\in\{1,2\},$ of radius $\sqrt{\frac{a}{\pi}}$ and volume $\frac{a^2}{2}$.

\bthm\label{Thm:Schlenk}
Let $a\in\R_{>0}$ and $a_1,a_2\in[0,a]$. Then
\beq \diamondsuit(a,a_1,a_2)\times\square(1) \stackrel{sympl.}{\cong} B^4_{\sqrt{\frac{a}{\pi}}}.\label{eq:schlenk1}\eeq
\ethm

We assembled the proof of this theorem, on the one hand, from the ideas contained in different parts of Schlenk's work in \cite{Schlenk2005}, more precisely, from Proposition 3.1.2, Lemmata 3.1.5, 3.1.8, 5.3.1, Example 3.1.7, and from the introduction within Section 9.3, on the other hand, from what has been carried out by Latschev, McDuff, and Schlenk in \cite{LatschevMcDuffSchlenk2013}. To the author's knowledge, so far, parts of the proof for the case
\beqq a_1=a_2=0,\eeqq
are written down in \cite{Schlenk2005}, furthermore, the full proof for the diamond-case, i.e.,
\beqq a_1=a_2=\frac{a}{2},\eeqq
by Latschev, McDuff, and Schlenk in \cite[Corollary 4.2]{LatschevMcDuffSchlenk2013}, and, via different methods, by Traynor in \cite{Traynor1995}. Moreover, Latschev, McDuff, and Schlenk proved in \cite[Proposition 4.4]{LatschevMcDuffSchlenk2013} that a similar result holds for so-called distorted diamonds which are distortions of the standard diamond $a_1=a_2=\frac{a}{2}$ and consist of a rectangle, a top and bottom triangle and two flaps (see \cite[Figure 3(II)]{LatschevMcDuffSchlenk2013}). For further symplectic packings in $4$ dimensions, we refer to \cite{Biran1997}, \cite{Biran1999}, and \cite{Karshon1994}.

Therefore, in what follows, it makes sense to give a rigorous proof.

\begin{figure}[h!]
\centering
\def\svgwidth{400pt}
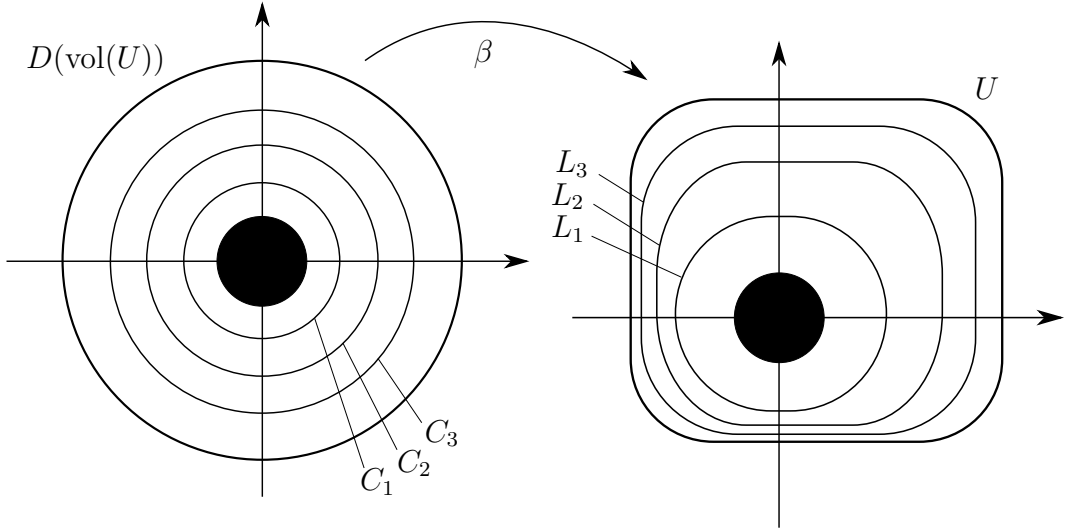
\caption[The admissibility of families of loops]{Here, $U\subset\R^2$ is a simply connected domain with $p=0$. The figure shows a few members of the families of concentric circles and of closed curves in $\mathcal{L}$.
The concentric circles $C_i$ are mapped by $\beta$ to the loops $L_i$ of enclosing the same area while, additionally, a neighbourhood of the origin is preserved.}
\label{img:admissible}
\end{figure}

For the proof of Theorem \ref{Thm:Schlenk}, we need the following preparations: For $a\in\R_{>0}$, we define for $i\in\{1,2\}$ the open discs
\beqq D(a)=\left\{(x_i,y_i):x_i^2+y_i^2<\frac{a}{\pi}\right\}\subset\R^2(z_i)\eeqq
of area $a$ and call a family $\mathcal{L}$ of loops in a simply connected domain $U\subset\R^2$ \textit{admissible} if there is a diffeomorphism
\beqq \beta: D\left(\vol(U)\right)\setminus\{0\} \rightarrow U \setminus \{p\}\eeqq
for some point $p\in U$ such that concentric circles are mapped to elements of $\mathcal{L}$ and in a neighbourhood of the origin $\beta$ is a translation (see Figure \ref{img:admissible}).
Then, we have:
\blem[Lemma 3.1.5 in \cite{Schlenk2005}]\label{Lem:Schlenk}
Let $U$ and $V$ be bounded and simply connected domains in $\R^2$ of equal area and let $\mathcal{L}_U$ and $\mathcal{L}_V$ be admissible families of loops in $U$ and $V$, respectively. Then, there is an area- and orientation-preserving diffeomorphism, i.e., a symplectomorphism, between $U$ and $V$ mapping loops to loops.
\elem 

Then, we come to the proof of Theorem \ref{Thm:Schlenk}:

\bpf[Proof of Theorem \ref{Thm:Schlenk}]
Let $a\in\R_{>0}$ and $a_i\in[0,a]$ for $i\in\{1,2\}$. For every $\eps >0$, we will prove that
\begin{align}
B^4_{\sqrt{\frac{a}{\pi}}} \stackrel{sympl.}{\hookrightarrow} & (1+\eps)\diamondsuit(a,a_1,a_2)\times \square(1) \label{eq:schlenk2}\\
\stackrel{sympl.}{\hookrightarrow} & (1+\eps)\left(\diamondsuit(a,a_1,a_2)\times \square(1)\right). \notag
\end{align}
Together with the below Lemma \ref{Lem:LMSLemma}, this would imply \eqref{eq:schlenk1}.

\begin{figure}[h!]
\centering
\def\svgwidth{400pt}
\begingroup%
  \makeatletter%
  \providecommand\color[2][]{%
    \errmessage{(Inkscape) Color is used for the text in Inkscape, but the package 'color.sty' is not loaded}%
    \renewcommand\color[2][]{}%
  }%
  \providecommand\transparent[1]{%
    \errmessage{(Inkscape) Transparency is used (non-zero) for the text in Inkscape, but the package 'transparent.sty' is not loaded}%
    \renewcommand\transparent[1]{}%
  }%
  \providecommand\rotatebox[2]{#2}%
  \newcommand*\fsize{\dimexpr\f@size pt\relax}%
  \newcommand*\lineheight[1]{\fontsize{\fsize}{#1\fsize}\selectfont}%
  \ifx\svgwidth\undefined%
    \setlength{\unitlength}{380.57992874bp}%
    \ifx\svgscale\undefined%
      \relax%
    \else%
      \setlength{\unitlength}{\unitlength * \real{\svgscale}}%
    \fi%
  \else%
    \setlength{\unitlength}{\svgwidth}%
  \fi%
  \global\let\svgwidth\undefined%
  \global\let\svgscale\undefined%
  \makeatother%
  \begin{picture}(1,0.40186203)%
    \lineheight{1}%
    \setlength\tabcolsep{0pt}%
    \put(0,0){\includegraphics[width=\unitlength,page=1]{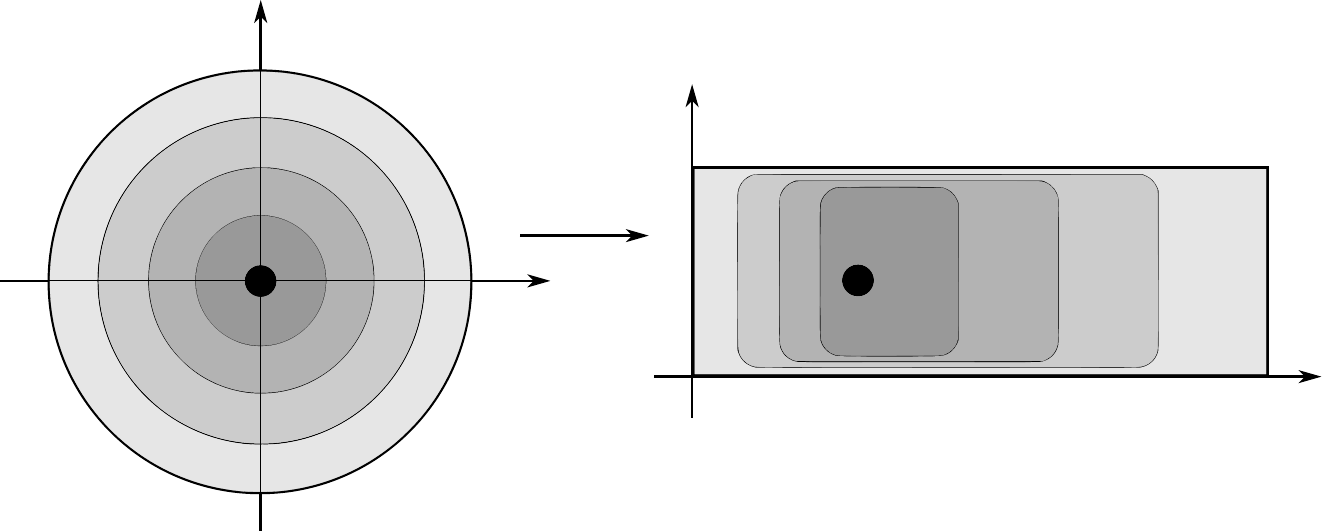}}%
    \put(0.41705963,0.23735365){\color[rgb]{0,0,0}\makebox(0,0)[lt]{\lineheight{1.25}\smash{\begin{tabular}[t]{l}$\alpha_i$\end{tabular}}}}%
    \put(0,0){\includegraphics[width=\unitlength,page=2]{alpha_i.pdf}}%
    \put(0.63858511,0.08547897){\color[rgb]{0,0,0}\makebox(0,0)[lt]{\lineheight{1.25}\smash{\begin{tabular}[t]{l}$\widetilde{a}_i$\end{tabular}}}}%
    \put(0.95237395,0.08888574){\color[rgb]{0,0,0}\makebox(0,0)[lt]{\lineheight{1.25}\smash{\begin{tabular}[t]{l}$a$\end{tabular}}}}%
    \put(0.86284305,0.28562654){\color[rgb]{0,0,0}\makebox(0,0)[lt]{\lineheight{1.25}\smash{\begin{tabular}[t]{l}$R(a)$\end{tabular}}}}%
    \put(0,0){\includegraphics[width=\unitlength,page=3]{alpha_i.pdf}}%
    \put(0.49846114,0.27057629){\color[rgb]{0,0,0}\makebox(0,0)[lt]{\lineheight{1.25}\smash{\begin{tabular}[t]{l}$1$\end{tabular}}}}%
    \put(0.29681082,0.32777606){\color[rgb]{0,0,0}\makebox(0,0)[lt]{\lineheight{1.25}\smash{\begin{tabular}[t]{l}$D(a)$\end{tabular}}}}%
    \put(0,0){\includegraphics[width=\unitlength,page=4]{alpha_i.pdf}}%
  \end{picture}%
\endgroup%

\caption[Visualization of $\alpha_i:D(a)\subset \R^2(z_i)\rightarrow R(a)\subset \R^2(z_i)$]{Visualization of $\alpha_i:D(a)\subset \R^2(z_i)\rightarrow R(a)\subset \R^2(z_i)$. It is $\widetilde{a}_i=\frac{a-\eps'}{a}a_i+\frac{\eps'}{2}$.}
\label{img:alpha_i}
\end{figure}

So, let us prove \eqref{eq:schlenk2}: For $\eps>0$, we define $\eps'$ by
\beqq \frac{2\eps'}{a}=\eps.\eeqq
Further, we define the open rectangle
\beqq R(a)=\{(u,v):u\in(0,a),v\in (0,1)\}\subset\R^2.\eeqq
By using Lemma \ref{Lem:Schlenk}, we construct two area- and orientation-preserving diffeomorphisms (i.e., symplectomorphisms)
\beqq \alpha_i: D(a)\subset \R^2(z_i)\rightarrow R(a)\subset \R^2(z_i),\quad i\in\{1,2\},\eeqq
such that for the first coordinate in the image $R(a)$ we have
\beq a_i-\frac{a_i}{a}\pi|z_i|^2-\eps' \leq \alpha_i(z_i)_1 \leq a_i + \frac{a-a_i}{a}\pi|z_i|^2+\eps'\label{eq:schlenk4}\eeq
for all $z_i\in D(a)$, while $\alpha_i$ is illustrated by Figure \ref{img:alpha_i}.

Let us discuss the construction of the $\alpha_i$s for $i\in\{1,2\}$: Let $i\in\{1,2\}$. For the construction of $\alpha_i$, we extend the construction in \cite[Proof of Lemma 3.1.8.]{Schlenk2005} for $a_i=0$ to general $a_i\in [0,a]$: In an \say{optimal world} we would choose the loops $\widehat{L}_s$, $0<s<1$, in the image $R(a)$ as the boundaries of the rectangles with corners
\beqq \left(a_i-sa_i,0\right),\left(a_i+s(a-a_i),0\right),\left(a_i-sa_i,1\right),\left(a_i+s(a-a_i),1\right).\eeqq
If the family $\widehat{\mathcal{L}}=\{\widehat{L}_s\}$ induces a map $\widehat{\alpha}_i$, we would then have
\beqq a_i-\frac{a_i}{a}\pi|z_i|^2 \leq \widehat{\alpha}_i(z_i)_1 \leq a_i + \frac{a-a_i}{a}\pi|z_i|^2\eeqq
for all $z_i\in D(a)$. The non-admissible family of loops $\widehat{\mathcal{L}}$ can be perturbed to an admissible family of loops $\mathcal{L}$ in such a way that the induced map $\alpha_i$ satisfies \eqref{eq:schlenk4}.

Indeed, choose the translation disc appearing in \cite[Proof of Lemma 3.1.5.]{Schlenk2005} of radius $\frac{\eps'}{8}$ centered at
\beqq (u_0,v_0)=\left(\frac{a-\eps'}{a}a_i+\frac{\eps'}{2},\frac{1}{2}\right).\eeqq
\begin{figure}[h!]
\centering
\def\svgwidth{360pt}
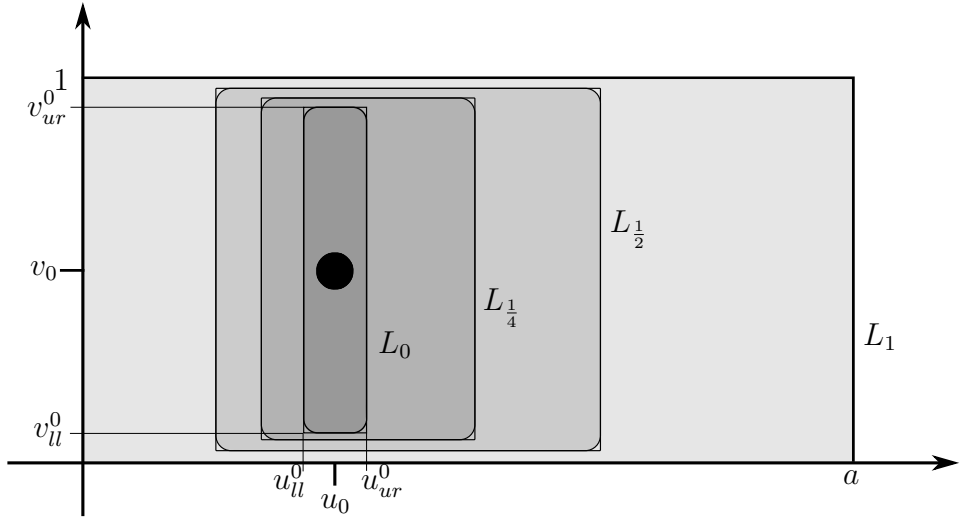
\caption[Visualization of the family of loops in $R(a)$]{Visualization of the family of loops in $R(a)$.}
\label{img:embedding}
\end{figure}
For $r<\frac{\eps'}{8}$, the loops $L(r)$ are therefore the concentric circles centered at $(u_0,v_0)$. in what follows, all rectangles considered have edges parallel to the coordinate axes. We may thus describe a rectangle by specifying its lower left and upper right corner. Let $\widebar{L}_0$ be the boundary of the rectangle with lower left corner
\beqq \left(u_{ll}^0,v_{ll}^0\right)=\left(\frac{a-\eps'}{a}a_i+\frac{\eps'}{4},\frac{\eps'}{4a}\right)\eeqq
and with upper right corner
\beqq \left(u_{ur}^0,v_{ur}^0\right)=\left(\frac{a-\eps'}{a}a_i+\frac{3\eps'}{4},1-\frac{\eps'}{4a}\right),\eeqq
and let $\widebar{L}_1$ be the boundary of $R(a)$. We define a family of loops $\widebar{L}_s$ by linearly interpolating between $\widebar{L}_0$ and $\widebar{L}_1$, i.e., $\widebar{L}_s$ is the boundary of the rectangle with corners
\beqq \left(u_{ll}^s,v_{ll}^s\right)=\left((1-s)\left(\frac{a-\eps'}{a}a_i+\frac{\eps'}{4}\right),(1-s)\frac{\eps'}{4a}\right)\eeqq
and
\beqq \left(u_{ur}^s,v_{ur}^s\right)=\left(\frac{a-\eps'}{a}a_i+\frac{3\eps'}{4}+s\left(a-\left(\frac{a-\eps'}{a}a_i+\frac{3\eps'}{4}\right)\right),1-\frac{\eps'}{4a}+s\frac{\eps'}{4a}\right)\eeqq
for $s\in[0,1]$. Let $\{L_s\}$, $s\in [0,1)$, be the smooth family of smooth loops obtained from $\{\widebar{L}_s\}$ by smoothing the corners as indicated in Figure \ref{img:embedding}. By choosing the smooth corners of $L_s$ more and more rectangular as $s\rightarrow 1$, we can arrange that the set
\beqq \bigcup_{0<s<1}L_s\eeqq
is the domain bounded by $L_0$ and $\widebar{L}_1$. Moreover, by choosing all smooth corners rectangular enough, we can arrange that the area enclosed by $L_s$ and $\widebar{L}_s$ is less than $\frac{\eps'}{4}$. We complete the families of loops
\beqq \{L(r)\} \; \text{ and } \; \{L_s\}\eeqq
to an admissible family of loops $\mathcal{L}$ in $R(a)$ and let
\beqq \alpha_i:D(a)\rightarrow R(a)\eeqq
be the map defined by $\mathcal{L}$. By construction, if $\alpha_i(z_i)_1$ lies on a loop in
\beqq \mathcal{L}\setminus \{L_s\}_{s\in (0,1)}\eeqq
as well as if $\alpha_i(z_i)_1$ lies on a loop in $L_s$ for some $s\in(0,1)$, then \eqref{eq:schlenk4} is satisfied. This completes the construction of a symplectomorphism
\beqq \alpha_i:D(a)\rightarrow R(a)\eeqq
satisfying \eqref{eq:schlenk4}.

Given the symplectomorphisms $\alpha_i$, $i\in\{1,2\}$, as in \cite[Lemma 3.1.8(i)]{Schlenk2005} for $a_1=a_2=0$, for the generalized case, we then conclude that
\beqq \left\{(\alpha_1(z_1)_2,\alpha_2(z_2)_2)\in\R^2:\pi|z_1|^2+\pi|z_2|^2<a\right\}\eeqq
symplectically embeds into $\square(1)$ and
\beqq \left\{(\alpha_1(z_1)_1,\alpha_2(z_2)_1)\in\R^2 : \pi|z_1|^2+\pi|z_2|^2<a \right\}\eeqq
symplectically embeds into
\beqq (1+\eps)\diamondsuit(a,a_1,a_2).\eeqq
Therefore:
\beqq B^4_{\sqrt{\frac{a}{\pi}}} \stackrel{sympl.}{\hookrightarrow} (1+\eps)\left(\diamondsuit(a,a_1,a_2)\times \square(1)\right).\eeqq

\epf

\blem[Lemma 4.3 in \cite{LatschevMcDuffSchlenk2013}]\label{Lem:LMSLemma}
For $a>0$, let $V\subset\R^4$ be a bounded domain such that for each compact subset $K\subset V$, there exists $\widehat{a}<a$ and a symplectic embedding
\beqq \widehat{\phi}:B^4_{\widehat{a}}\rightarrow V \quad \text{such that} \quad K \subset \widehat{\phi}(B^4_{\widehat{a}}).\eeqq
Then, $V$ is symplectomorphic to $B^4_a$.
\elem

We note that the proof of Lemma \ref{Lem:LMSLemma} is based on the results of McDuff in \cite{McDuff1991} which involve the concept of $J$-holomorphic curves.

We remark that, in the proof of Theorem \ref{Thm:Schlenk}, instead of referring to Lemma \ref{Lem:LMSLemma} in order to get rid of the $\eps$ in \eqref{eq:schlenk2}, we also can refer to a later developed idea\footnote{The idea is based on the construction of a one-paramter-family of $\eps$-embeddings.} of Buhovsky (carried out by Pelayo and Ng\d{o}c in \cite{PelayoNgoc2016}) which also holds in higher dimensions. For that, we refer to Lemma 8.2 in \cite{Schlenk2018}.

\section{The triangle-case}\label{Sec:triangle}

We begin with the following proposition:

\bprop\label{Prop:DeltaJDelta}
Let $\Delta\subset\R^2$ be any triangle. Then, we have
\beqq c_{EHZ}(\Delta\times J\Delta)=\vol(\Delta).\eeqq
\eprop

\bpf
Applying Theorem \ref{Thm:Chap7relationship}, we first notice that
\beqq c_{EHZ}\left(\Delta\times J\Delta\right)=\min_{q\in M_3(\Delta,J\Delta)}\ell_{J\Delta}(q).\eeqq

\begin{figure}[h!]
\centering
\def\svgwidth{265pt}
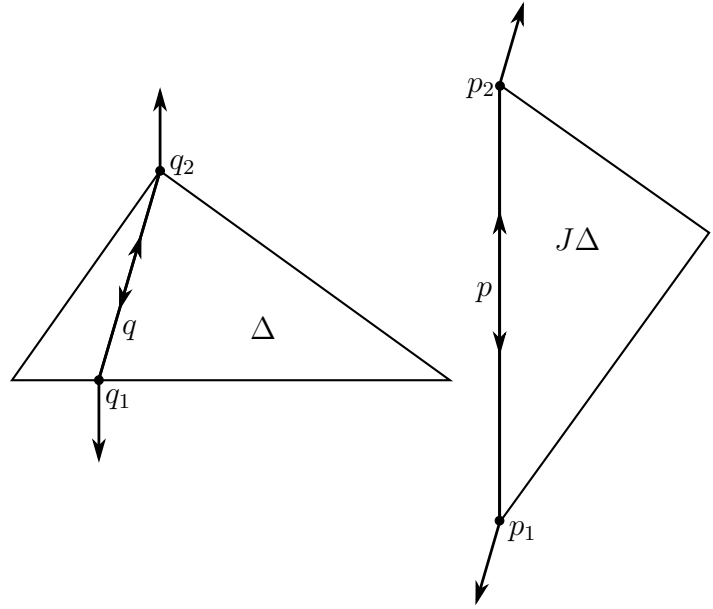
\caption[Visualisation of $\Delta$ and $J\Delta$]{Visualisation of $\Delta$ and $J\Delta$ together with the closed $2$-bouncing $(\Delta,J\Delta)$-Minkowski billiard trajectory $q=(q_1,q_2)$ and its dual billiard trajectory $p=(p_1,p_2)$ in $J\Delta$.}
\label{img:DeltaJDelta}
\end{figure}

We first consider the closed $(\Delta,J\Delta)$-Minkowski billiard trajectories with $2$ bouncing points. Let $q=(q_1,q_2)$ be any of these. Then, since
\beqq N_{\Delta}(q_1)\cap - N_{\Delta}(q_2) \neq \emptyset,\eeqq
we can assume that, without loss of generality, $q_2$ is a vertex of $\Delta$ and $q_1$ lies on the opposite side. In any case, the closed polygonal line $p=(p_1,p_2)$ is one of the associated closed dual billiard trajectories in $J\Delta$ (note that, in general, the closed dual billiard trajectories are not unique), when $p_1$ and $p_2$ are the vertices enclosing that side of $J\Delta$ which is the $J$-rotated side of $\Delta$ that contains $q_1$ (see Figure \ref{img:DeltaJDelta}).
Recalling \cite[Proposition 2.2]{KruppRudolf2022} and $\left(J\Delta\right)^\circ = J \Delta^\circ$, one calculates
\begin{align*}
\ell_{J\Delta}(q)=\mu_{J\Delta^\circ}(q_2-q_1)+\mu_{J\Delta^\circ}(q_1-q_2)&=\langle q_2-q_1,p_1\rangle + \langle q_1-q_2,p_2\rangle \\
&=\langle q_2-q_1,p_1-p_2\rangle \\
&=2\vol(\Delta).
\end{align*}
Therefore, we conclude
\beq \min_{q\in M_2(\Delta,J\Delta)}\ell_{J\Delta}(q)=2\vol(\Delta).\label{eq:DeltaJDelta0}\eeq

Now, we consider the closed $(\Delta,J\Delta)$-Minkowski billiard trajectories with $3$ bouncing points. We distinguish between regular and non-regular ones, where we recall to call a closed polygonal curve with vertices on $\partial \Delta$ \textit{regular}, if all its vertices are smooth boundary points of $\Delta$, otherwise we call it \textit{non-regular}.

Every non-regular closed $(\Delta,J\Delta)$-Minkowski billiard trajectory $q$ with $3$ bouncing points $q_1,q_2,q_3$ has the property that one of its vertices is a vertex of $\Delta$, say $q_2$ in Figure \ref{img:DeltaJDelta}, and another one lies on the opposite side of $\Delta$, say $q_1$ in Figure \ref{img:DeltaJDelta}. This follows from the fact that, otherwise, the normal vectors at $\Delta$ in the bouncing points $q_1,q_2,q_3$ do not surround the origin which is a contradiction to what has been shown within the proof of \cite[Proposition 3.9]{KruppRudolf2022}. Then, the connecting line of these two vertices interpreted as closed polygonal curve in $F_2^{cp}(\Delta)$--in this case, it is the closed $2$-bouncing $(\Delta,J\Delta)$-Minkowski billiard trajectory from above--necessarily has less or equal $\ell_{J\Delta}$-length than $q$ (see \cite[Proposition 2.3(i)]{KruppRudolf2022}). Referring to Theorem \ref{Thm:Chap7relationship} and \eqref{eq:DeltaJDelta0}, we conclude
\beqq \ell_{J\Delta}(q)\geq \min_{q\in F_2^{cp}(\Delta)}\ell_{J\Delta}(q)=\min_{q\in M_2(\Delta,J\Delta)}\ell_{J\Delta}(q)=2\vol(\Delta).\eeqq

\begin{figure}[h!]
\centering
\def\svgwidth{400pt}
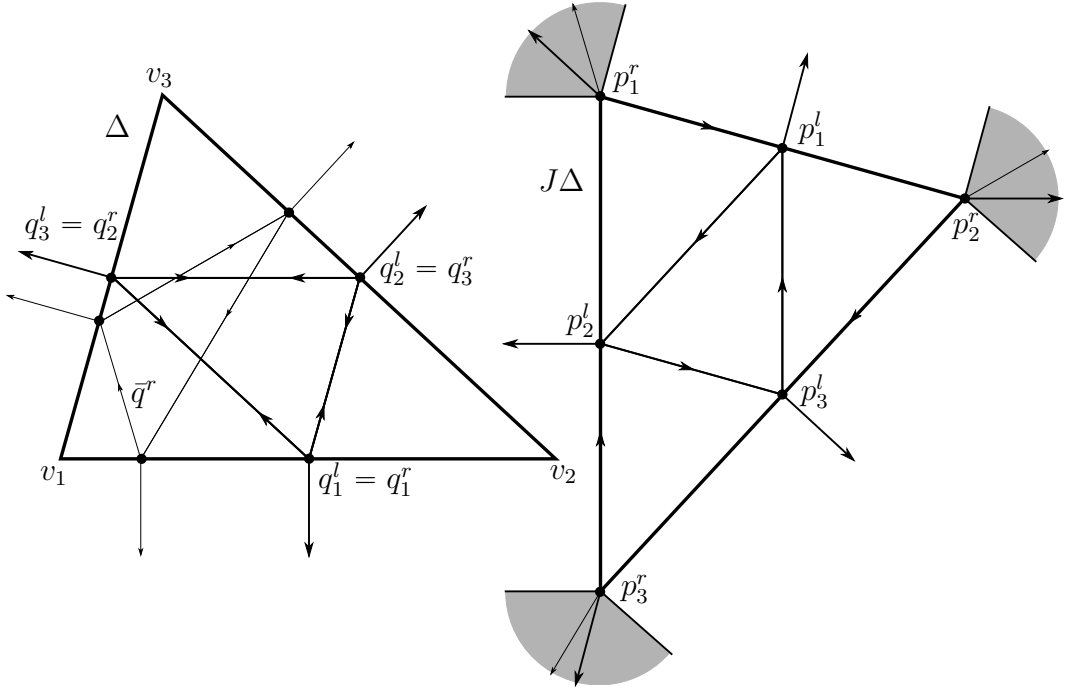
\caption[The figure illustrates the regular closed $3$-bouncing $(\Delta,J\Delta)$-Minkowski billiard trajectories]{The regular closed $3$-bouncing $(\Delta,J\Delta)$-Minkowski billiard trajectories $q^l=(q_1^l,q_2^l,q_3^l)$ and $q^r=(q_1^r,q_2^r,q_3^r)$ together with their closed dual billiard trajectories $p^l=(p_1^l,p_2^l,p_3^l)$ and $p^r=(p_1^r,p_2^r,p_3^r)$, respectively. Besides $q^l$ and $q^r$, $\bar{q}^r$ is another regular closed $(\Delta,J\Delta)$-Minkowski billiard trajectory with closed dual billiard trajectory $p^r$, but which corresponds to different normal vectors in the normal cones at $J\Delta$ in $p_1^r$, $p_2^r$, and $p_3^r$.}
\label{img:DeltaJDelta2}
\end{figure}

\begin{figure}[h!]
\centering
\def\svgwidth{350pt}
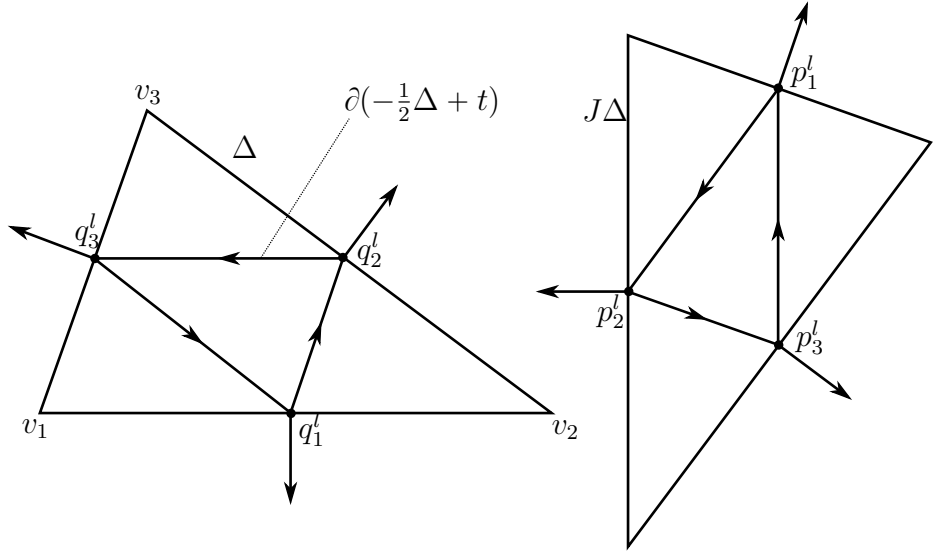
\caption[$q^l$ is the counter-clockwise passed boundary of $\partial(-\frac{1}{2}\Delta +t)$]{Reflecting on the algorithm for the construction of $q^l$, yields that $q^l$ is the counter-clockwise passed boundary of $-\frac{1}{2}\Delta +t$, where $t\in\R^2$ is chosen such that the vertices of $-\frac{1}{2}\Delta +t$ are on $\partial \Delta$. It is clear that $-\frac{1}{2}\Delta$ is the unique (up to translation) minimizing triangle of \eqref{eq:minproblemlambda}. Sidenote: A similar argumentation yields that $p^l$ is the counter-clockwise passed boundary of $-\frac{1}{2}J\Delta+u$, where $u\in\R^2$ is chosen such that the vertices of $-\frac{1}{2}J\Delta + u$ are on $\partial (-J\Delta)$.}
\label{img:MinusTriangleInTriangle}
\end{figure}

It remains to understand the regular closed $3$-bouncing $(\Delta,J\Delta)$-Minkowski billiard trajectories. By referring to the algorithm presented in \cite{KruppRudolf2022}, we notice that these regular closed $3$-bouncing $(\Delta,J\Delta)$-Minkowski billiard trajectories are (partly) determined by the outer normal vectors at the interiors of the sides of $\Delta$. Both orientations, i.e., both the clock- and counter-clockwise order of the three sides of $\Delta$--produce one uniqely determined closed dual billiard trajectory in $J\Delta$, respectively, say $p^l=(p_1^l,p_2^l,p_3^l)$ and $p^r=(p_1^r,p_2^r,p_3^r)$ (see Figure \ref{img:DeltaJDelta2}). These two closed dual billiard trajectories produce--after a choice of normal vectors in the normal cones at $p_1^r,p_2^r,p_3^r$--the following two trajectories in $\Delta$: $q^l=(q_1^l,q_2^l,q_3^l)$ and $q^r=(q_1^r,q_2^r,q_3^r)$. By construction, $q^l$ and $q^r$ are regular closed $(\Delta,J\Delta)$-Minkowski billiard trajectories with $3$ bouncing points.
By varying the normal vectors in the normal cones at $p_1^r$, $p_2^r$, and $p_3^r$, we can find more regular closed $3$-bouncing $(\Delta,J\Delta)$-Minkowski billiard trajectories $\bar{q}^r$ instead of $q^r=(q_1^r,q_2^r,q_3^r)$. Nevertheless, it is sufficient to just concentrate on $q^r$ since all other $\bar{q}^r$ have the same $\ell_{J\Delta}$-length. This follows from the fact that they all have their associated closed dual billiard trajectories in $J\Delta$ in common, i.e., $p^r$, and applying \cite[Proposition 3.4]{KruppRudolf2022}, this means
\beq \ell_{J\Delta}(\bar{q}^r)=\ell_{-\Delta}(p^r)=\ell_{J\Delta}(q^r).\label{eq:starequality}\eeq
Therefore, the two pairs $(q^{l,r},p	^{l,r})$--as indicated in Figure \ref{img:DeltaJDelta2}--are the only representatives of the regular closed $3$-bouncing $(\Delta,J\Delta)$-Minkowski billiard trajectories to which we must refer below. From \cite[Proposition 3.11(iv)]{KruppRudolf2022}, \cite[Lemma 5.4]{Rudolf2022b}, Proposition \ref{Prop:Chap7characteristicgeodesic} (note that the closed characteristic on $\partial \Delta$ is uniquely given by passing through $\partial \Delta$ clockwise), and that in two dimensions the volume coincides with the capacity, it follows
{\allowdisplaybreaks\begin{align*}
\ell_{J\Delta}(q^r)&=\ell_{J\Delta}\left(v_3-v_1\right)+\ell_{J\Delta}\left(v_2-v_3\right)+\ell_{J\Delta}\left(v_1-v_2\right)\\
&=2\left(\ell_{\frac{J\Delta}{2}}\left(v_3-v_1\right)+\ell_{\frac{J\Delta}{2}}\left(v_2-v_3\right)+\ell_{\frac{J\Delta}{2}}\left(v_1-v_2\right)\right)\\
&=2c_{EHZ}(\Delta)\\
&=2\vol(\Delta),
\end{align*}}%
where by $v_1,v_2,v_3$ we denote the vertices of $\Delta$ as indicated in Figure \ref{img:DeltaJDelta2}. In order to calculate $\ell_{J\Delta}(q^l)$, it is useful to note that $q^l$ is the counter-clockwise passed boundary of the minimizing triangle $-\lambda \Delta$ of the minimization problem
\beq \min\left\{\lambda : -\lambda \Delta \in F(\Delta)\right\}.\label{eq:minproblemlambda}\eeq
This is a consequence of applying the algorithm--presented in \cite{KruppRudolf2022}--for manually determining the Minkowski billiard trajectories (see Figure \ref{img:MinusTriangleInTriangle}).
As we will prove below in Lemma \ref{Lem:minproblemlambda}, the minimum in \eqref{eq:minproblemlambda} is $\lambda=\frac{1}{2}$. This implies that
\beqq \ell_{J\Delta}(q^l)=\frac{1}{2}\left(\ell_{J\Delta}\left(v_3-v_1\right)+\ell_{J\Delta}\left(v_2-v_3\right)+\ell_{J\Delta}\left(v_1-v_2\right)\right)=c_{EHZ}(\Delta)=\vol(\Delta) \eeqq 
since $q_2^l-q_1^l$ and $v_3-v_1$, $q_3^l-q_2^l$ and $v_1-v_2$, as well as $q_1^l-q_3^l$ and $v_2-v_3$ are parallel, point into the same direction and differ by a factor $2$, respectively (note that the above characterization of $q^l$ implies that $q_1^l,q_2^l,q_3^l$ each are located on the center of the sides of $\Delta$).

Finally, this implies
\beqq c_{EHZ}(\Delta\times J\Delta)=\min_{q\in M_3(\Delta,J\Delta)}\ell_{J\Delta}(q)=\ell_{J\Delta}(q^l)=\vol(\Delta).\eeqq
\epf

\blem\label{Lem:minproblemlambda}
Let $\Delta$ be any triangle in $\R^2$. Then, one has
\beq \min\left\{\lambda : -\lambda \Delta \in F(\Delta)\right\} = \frac{1}{2}. \label{eq:Lemminproblemlambda}\eeq
\elem

\begin{figure}[h!]
\centering
\def\svgwidth{240pt}
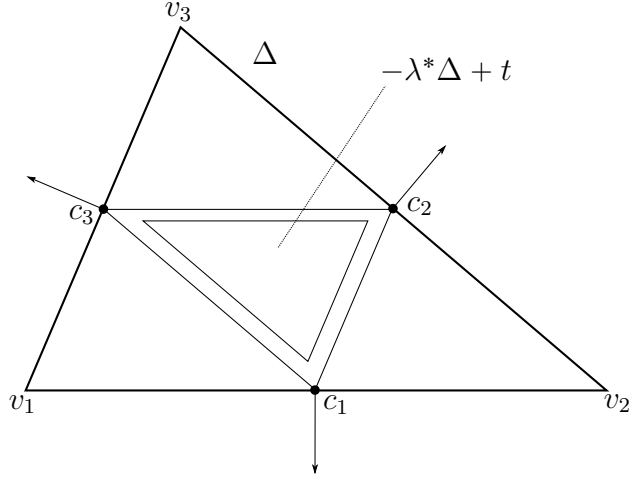
\caption[The triangle with vertices $c_1,c_2,c_3$ is the unique minimizer (up to translation) of \eqref{eq:Lemminproblemlambda}.]{The triangle with vertices $c_1,c_2,c_3$ is the unique minimizer (up to translation) of \eqref{eq:Lemminproblemlambda}. Every other triangle $-\lambda^*\Delta$ with $\lambda^*<\frac{1}{2}$ can be translated such that it lies in the interior of the triangle which has the vertices $c_1,c_2,c_3$.}
\label{img:minproblemlambda}
\end{figure}

\bpf
Let $c_1,c_2,c_3$ be the centers of the sides of $\Delta$ (see Figure \ref{img:minproblemlambda}). By the definition of $c_1,c_2,c_3$, the triangle $c=(c_1,c_2,c_3)$ has the form of the boundary of $-\frac{1}{2}\Delta$ (one easily checks that $c_2-c_1$ and $v_3-v_1$, $c_3-c_2$ and $v_1-v_2$, as well as $c_1-c_3$ and $v_2-v_3$ are parallel, point into the same direction, and differ by a factor $2$, respectively). Since the normal vectors at $\Delta$ in $c_1,c_2,c_3$ surround the origin, $c$ and therefore also $-\frac{1}{2}\Delta$ is in $F(\Delta)$ (see \cite[Lemma 2.1(ii)]{KruppRudolf2020}). This implies
\beq \min\left\{\lambda : -\lambda \Delta \in F(\Delta)\right\} \leq \frac{1}{2}.\label{eq:Lemminproblem}\eeq
If the left side in \eqref{eq:Lemminproblem} is less than $\frac{1}{2}$, then we can find a $t\in\R^2$ such that $-\lambda^*\Delta + t$ is a subset of the interior of the triangle which is given by the vertices $c_1,c_2,c_3$. This implies that $-\lambda^*\Delta$ can be translated into $\mathring{\Delta}$, i.e.,
\beqq -\lambda^*\Delta \notin F(\Delta).\eeqq
This is a contradiction to the assumption that the left side in \eqref{eq:Lemminproblem} is less than $\frac{1}{2}$. Therefore, it follows \eqref{eq:Lemminproblemlambda}.
\epf

Based on Proposition \ref{Prop:DeltaJDelta}, we now can make the following observations. First of all, we note that due to the scale invariance of Viterbo's conjecture in the sense of \cite[Proposition 2.9]{Rudolf2022b}, without loss of generality, we can require the convex body $T\subset\R^2$ in \eqref{eq:Viterboconjecture3} to fulfill
\beqq \vol(T)=2 \vol(\Delta).\eeqq
Then, \eqref{eq:Viterboconjecture3} is equivalent to
\beqq \vol(\Delta)\geq \frac{c_{EHZ}(\Delta\times T)}{2}.\eeqq
Using Proposition \ref{Prop:DeltaJDelta}, this becomes
\beq c_{EHZ}(\Delta\times J\Delta) \geq \frac{c_{EHZ}(\Delta\times T)}{2}.\label{eq:basetrianglecase}\eeq

Based on \eqref{eq:basetrianglecase}, we will proceed in four steps: Lemmata \ref{Lem:TriangleCaseLem1}, \ref{Lem:TriangleCaseLem2}, and \ref{Lem:TriangleCaseLem3}, and Proposition \ref{Prop:TriangleCaseProp4}.

\blem\label{Lem:TriangleCaseLem1}
Let $\Delta$ be any triangle in $\R^2$. Then, we have
\beqq \vol(T^*)=2\vol(\Delta)\eeqq
for all $T^*\in\mathcal{T}_{\Delta}$.
\elem

\blem\label{Lem:TriangleCaseLem2}
Let $\Delta$ be any triangle in $\R^2$. Then, we have
\beqq c_{EHZ}(\Delta\times J\Delta) = \frac{c_{EHZ}(\Delta\times T^*)}{2}\eeqq
for all $T^*\in\mathcal{T}_{\Delta}$.
\elem

\blem\label{Lem:TriangleCaseLem3}
Let $\Delta$ be any triangle in $\R^2$. Then, we have
\beqq c_{EHZ}(\Delta \times T^*)>c_{EHZ}(\Delta \times T),\quad T^*\in\mathcal{T}_{\Delta},\eeqq
for all convex bodies $T\notin\mathcal{T}_{\Delta}$ with $\vol(T)=2\vol(\Delta)$.
\elem

\bprop\label{Prop:TriangleCaseProp4}
Let $\Delta$ be any triangle in $\R^2$. Then, we have that 
\beqq \mathring{\Delta}\times \mathring{T}^* \stackrel{sympl.}{\cong} B^4_{\sqrt{\frac{a}{\pi}}} \eeqq
for all $T^*\in\mathcal{T}_{\Delta,\square}$.
\eprop

Combining Lemmata \ref{Lem:TriangleCaseLem1}, \ref{Lem:TriangleCaseLem2}, and \ref{Lem:TriangleCaseLem3}, and Proposition \ref{Prop:TriangleCaseProp4}, implies Theorem \ref{Thm:answer1} for the case when $Q$ in \eqref{eq:answer1} is any triangle in $\R^2$, Theorem \ref{Thm:answer2}(i), and Theorem \ref{Thm:answer3}(i).

So, let us prove the above statements one after the other:

\bpf[Proof of Lemma \ref{Lem:TriangleCaseLem1}]
Let $\Delta$ be any triangle in $\R^2$ that, without loss of generality, is centred at the origin. Then, having Figure \ref{img:convexhull} in mind, the set of volume-minimizing convex hulls
\beq \conv\{J\Delta,-J\Delta +t\}, \quad t\in\R^2,\label{eq:Lem1convexhulls}\eeq
is the set of these convex hulls \eqref{eq:Lem1convexhulls} for which $t$ is in $-J\Delta$, i.e., for which $-J\Delta +t$ is a subset of $-2J\Delta$. If $t$ is in $-J\partial\Delta$, i.e., if $-J\Delta +t$ touches the boundary of $-2J\Delta$, then the convex hulls are parallelograms; if $t$ is in $-J\mathring{\Delta}$, i.e., if $-J\Delta +t$ lies in the interior of $-2J\Delta$, then the convex hulls are hexagons.

If the convex hulls are parallelograms, then it clearly follows that their volume--as product of one side-length (which is the length of one side of the $J$-rotated $\Delta$) and the associated height (which corresponds to the height of the aforementiond side of the $J$-rotated $\Delta$)--is two times the volume of $\Delta$.

\begin{figure}[h!]
\centering
\def\svgwidth{400pt}
\begingroup%
  \makeatletter%
  \providecommand\color[2][]{%
    \errmessage{(Inkscape) Color is used for the text in Inkscape, but the package 'color.sty' is not loaded}%
    \renewcommand\color[2][]{}%
  }%
  \providecommand\transparent[1]{%
    \errmessage{(Inkscape) Transparency is used (non-zero) for the text in Inkscape, but the package 'transparent.sty' is not loaded}%
    \renewcommand\transparent[1]{}%
  }%
  \providecommand\rotatebox[2]{#2}%
  \newcommand*\fsize{\dimexpr\f@size pt\relax}%
  \newcommand*\lineheight[1]{\fontsize{\fsize}{#1\fsize}\selectfont}%
  \ifx\svgwidth\undefined%
    \setlength{\unitlength}{258.90932108bp}%
    \ifx\svgscale\undefined%
      \relax%
    \else%
      \setlength{\unitlength}{\unitlength * \real{\svgscale}}%
    \fi%
  \else%
    \setlength{\unitlength}{\svgwidth}%
  \fi%
  \global\let\svgwidth\undefined%
  \global\let\svgscale\undefined%
  \makeatother%
  \begin{picture}(1,0.36993758)%
    \lineheight{1}%
    \setlength\tabcolsep{0pt}%
    \put(0,0){\includegraphics[width=\unitlength,page=1]{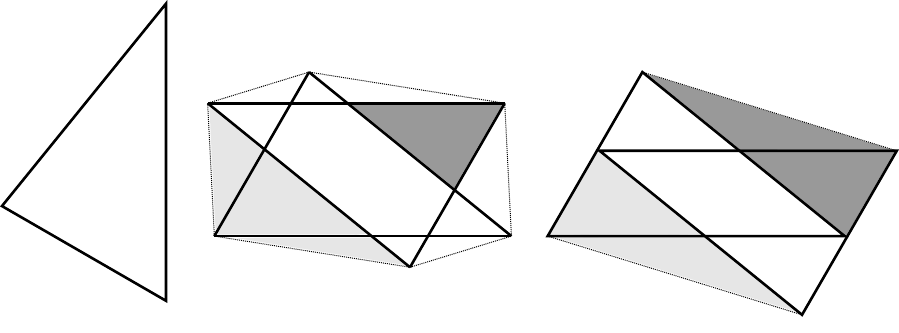}}%
    \put(0.07403141,0.16254947){\color[rgb]{0,0,0}\makebox(0,0)[lt]{\lineheight{1.25}\smash{\begin{tabular}[t]{l}$\Delta$\end{tabular}}}}%
    \put(0.27515458,0.33091002){\color[rgb]{0,0,0}\makebox(0,0)[lt]{\lineheight{1.25}\smash{\begin{tabular}[t]{l}$J\Delta$\end{tabular}}}}%
    \put(0.4308704,0.00860769){\color[rgb]{0,0,0}\makebox(0,0)[lt]{\lineheight{1.25}\smash{\begin{tabular}[t]{l}$-J\Delta+t'$\end{tabular}}}}%
    \put(0,0){\includegraphics[width=\unitlength,page=2]{volumecalculationtriangle.pdf}}%
    \put(0.63748107,0.31446834){\color[rgb]{0,0,0}\makebox(0,0)[lt]{\lineheight{1.25}\smash{\begin{tabular}[t]{l}$J\Delta$\end{tabular}}}}%
    \put(0,0){\includegraphics[width=\unitlength,page=3]{volumecalculationtriangle.pdf}}%
    \put(0.70758257,0.00386613){\color[rgb]{0,0,0}\makebox(0,0)[lt]{\lineheight{1.25}\smash{\begin{tabular}[t]{l}$-J\Delta+t''$\end{tabular}}}}%
    \put(0,0){\includegraphics[width=\unitlength,page=4]{volumecalculationtriangle.pdf}}%
  \end{picture}%
\endgroup%

\caption[The volumes of the parallelograms as well as the volumes of the hexagons equal two times the volume of $\Delta$.]{The volumes of the parallelograms as well as the volumes of the hexagons equal two times the volume of $\Delta$. The triangle $-J\Delta+t'$ can be translated along one of the sides of $J\Delta$ to $-J\Delta +t''$ such that $\conv\{J\Delta,-J\Delta +t''\}$ is a parallelogram.}
\label{img:volumecalculationtriangle}
\end{figure}

If the convex hulls are hexagons, i.e., say, if we consider $-J\Delta +t'$ in Figure \ref{img:volumecalculationtriangle}, where $t'$ is in $-J\mathring{\Delta}$, or, equivalently, $-J\Delta+t'$ is subset of $-2J\mathring{\Delta}$, then one can translate $-J\Delta+t'$ along one of the sides of $J\Delta$ to $-J\Delta+t''$ such that $\conv\{J\Delta,-J\Delta +t''\}$ is a parallelogram. Clearly, one has
\beqq \vol\left(\conv\{J\Delta,-J\Delta +t'\}\right) = \vol\left(\conv\{J\Delta,-J\Delta +t''\}\right)\eeqq
since the volumes of the two light grey areas, the volumes of the two dark grey areas, as well as the volumes of the two white enclosed areas coincide (the length of the base sides as well es the heights do not differ, respectively). Therefore, one has that the volumes of the hexagons equal two times the volume of $\Delta$.
\epf

\bpf[Proof of Lemma \ref{Lem:TriangleCaseLem2}]
Let $\Delta$ be any triangle in $\R^2$ and $T^*$ any convex body in $\mathcal{T}_\Delta$. We write
\beqq T^*=\conv\{J\Delta,-J\Delta+t^*\}.\eeqq
In order to calculate $c_{EHZ}(\Delta\times T^*)$, we recall that Theorem \ref{Thm:Chap7relationship} implies that
\beqq c_{EHZ}\left(\Delta\times T^*\right)=\min_{q\in M_3(\Delta,T^*)}\ell_{T^*}(q).\eeqq

\begin{figure}[h!]
\centering
\def\svgwidth{400pt}
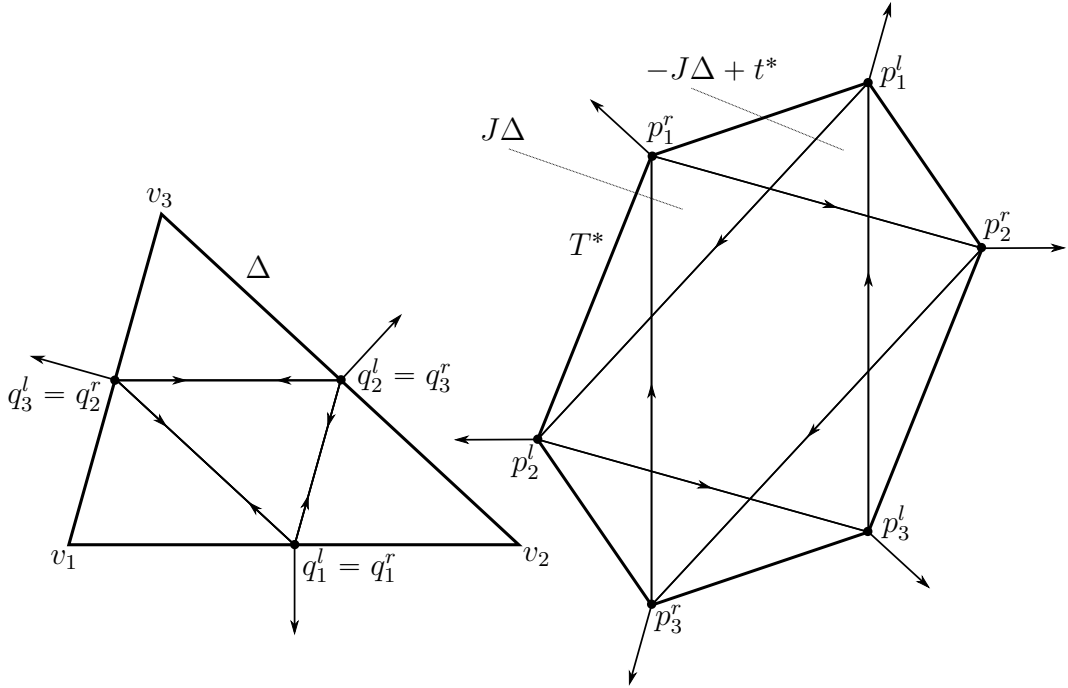
\caption[The figure illustrates the regular closed $3$-bouncing $(\Delta,T^*)$-Minkowski billiard trajectories]{The regular closed $(\Delta,T^*)$-Minkowski billiard trajectories $q^l=(q_1^l,q_2^l,q_3^l)$ and $q^r=(q_1^r,q_2^r,q_3^r)$ together with their closed dual billiard trajectories $p^l=(p_1^l,p_2^l,p_3^l)$ and $p^r=(p_1^r,p_2^r,p_3^r)$, respectively. The orbits of $p^l$ and $p_r$ can be represented by $-J\partial \Delta + t^*$ and $J\partial \Delta$, respectively.}
\label{img:DeltaTstern}
\end{figure}

We begin by considering the regular closed $(\Delta,T^*)$-Minkowski billiard trajectories with $3$ bouncing points. Referring to the algorithm presented in \cite{KruppRudolf2022}, as in the proof of Proposition \ref{Prop:DeltaJDelta}, every regular closed $(\Delta,T^*)$-Minkowski billiard trajectory with $3$ bouncing points has its bouncing points on the interiors of the three sides of $\Delta$. Consequently, the orbits of the corresponding closed dual billiard trajectories in $T^*$ are given by $J\partial \Delta$ and $-J\partial \Delta + t^*$ (which, in Figure \ref{img:DeltaTstern}, are representd by $p^l=(p_1^l,p_2^l,p_3^l)$ and $p^r=(p_1^r,p_2^r,p_3^r)$).
As in Proposition \ref{Prop:DeltaJDelta}, when searching for $\ell_{T^*}$-minimizing regular closed $(\Delta,T^*)$-Minkowski billiard trajectories with $3$ bouncing points, then it is enough to just concentrate on the two regular closed $(\Delta,T^*)$-Minkowski billiard trajectories
\beqq q^l=(q_1^l,q_2^l,q_3^l)\; \text{ and } \; q^r=(q_1^r,q_2^r,q_3^r)\eeqq
as indicated in Figure \ref{img:DeltaTstern} (they are characterized by the fact that $q_1^l=q_1^r$, $q_2^l=q_3^r$, and $q_3^l=q_2^r$ are the centres of the sides $[v_1,v_2]$, $[v_2,v_3]$, and $[v_3,v_1]$--as consequence of the choice of the normal vectors in $p_1^{r,l}$, $p_2^{r,l}$, $p_3^{r,l}$ normal to $p_3^{r,l}-p_2^{r,l}$, $p_1^{r,l}-p_3^{r,l}$, and $p_2^{r,l}-p_1^{r,l}$). Any other regular closed $(\Delta,T^*)$-Minkowski billiard trajectory--arising from a different choice of normal vectors in the normal cones at $p_1^l$, $p_2^l$, and $p_3^l$, as well as at $p_1^r$, $p_2^r$, and $p_3^r$--has the same length as $q^l$ and $q^r$, respectively (since $p^l$ and $p^r$ are the only closed dual billiard trajectories in $T^*$ that differ from each other; see \eqref{eq:starequality} and its explanation in the proof of Proposition \ref{Prop:DeltaJDelta}). So, it remains to calculate $\ell_{T^*}(q^l)$ and $\ell_{T^*}(p^l)$. Referring to \cite[Proposition 3.4]{KruppRudolf2022}, the proof of Propositon \ref{Prop:DeltaJDelta} (due to Lemma \ref{Lem:minproblemlambda}, the $\ell_{-\Delta}$-length of $p^l$ in Figure \ref{img:DeltaTstern} equals two times the $\ell_{-\Delta}$-length of $p^l$ in Figure \ref{img:DeltaJDelta2}; the $\ell_{-\Delta}$-length of $p^r$ in Figure \ref{img:DeltaTstern} equals the $\ell_{-\Delta}$-length of $p^r$ in Figure \ref{img:DeltaJDelta2}), and the fact that determining the length can be clearly traced back to a calculation solely depending on the positions of the bouncing points of $q$ and its dual billiard trajectory $p$, we calculate
\beqq \ell_{T^*}(q^l)=\ell_{-\Delta}(p^l)=2\ell_{J\Delta}(q^l)=2 \vol(\Delta)\eeqq
and
\beqq \ell_{T^*}(q^r)=\ell_{-\Delta}(p^r)=\ell_{J\Delta}(q^r)=2\vol(\Delta).\eeqq

in what follows, we show that there are neither non-regular closed $3$-bouncing $(\Delta,T^*)$-Minkowski billiard trajectories nor closed $2$-bouncing $(\Delta,T^*)$-Minkowski billiard trajectories which have a smaller $\ell_{T^*}$-length than $2\vol(\Delta)$.

Considering the closed $2$-bouncing $(\Delta,T^*)$-Minkowski billiard trajectories, from
\beqq \min_{q\in M_2(\Delta,T^*)}\ell_{T^*}(q)=\min_{q\in F_2^{cp}(\Delta)}\ell_{T^*}(q)\eeqq
(see Theorem \ref{Thm:Chap7relationship}), $J\Delta \subseteq T^*$, \cite[Proposition 3.11(iii)]{KruppRudolf2022}, and \eqref{eq:DeltaJDelta0}, we conclude
\beqq \min_{q\in M_2(\Delta,T^*)}\ell_{T^*}(q)=\min_{q\in F_2^{cp}(\Delta)}\ell_{T^*}(q) \geq \min_{q\in F_2^{cp}(\Delta)}\ell_{J\Delta}(q)=\min_{q\in M_2(\Delta,J\Delta)}\ell_{J\Delta}(q)=2\vol(\Delta).\eeqq

Considering the non-regular closed $3$-bouncing $(\Delta,T^*)$-Minkowski billiard trajectories, we show--similarly to the proof of Proposition \ref{Prop:DeltaJDelta}--that for every non-regular closed $3$-bouncing $(\Delta,T^*)$-Minkowski billiard trajectory $q$, there is a closed $2$-bouncing $(\Delta,T^*)$-Minkowski billiard trajectory with less or equal $\ell_{T^*}$-length. In fact, $q$ has the property that one of its bouncing points is a vertex of $\Delta$ and another one lies on the opposide side of $\Delta$. This follows from the fact that, otherwise, the normal vectors at $\Delta$ in the bouncing points of $q$ do not surround the origin which is a contradiction to what has been shown within the proof of \cite[Proposition 3.9]{KruppRudolf2022}. Then, the connecting line of the abovementioned two vertices interpreted as closed polygonal curve in $F_2^{cp}(\Delta)\subseteq F_3^{cp}(\Delta)$ necessarily has less or equal $\ell_{T^*}$-length than $q$ (see \cite[Proposition 2.3(i)]{KruppRudolf2022}). Referring to Theorem \ref{Thm:Chap7relationship}, we conclude
\beqq \ell_{T^*}(q)\geq \min_{q\in F_2^{cp}(\Delta)}\ell_{T^*}(q)=\min_{q\in M_2(\Delta,T^*)}\ell_{T^*}(q).\eeqq

Finally, together with Proposition \ref{Prop:DeltaJDelta}, we conclude
\beqq c_{EHZ}(\Delta\times T^*)=\min_{q\in M_3(\Delta,T^*)}\ell_{T^*}(q)=2\vol(\Delta)=2c_{EHZ}(\Delta\times J\Delta).\eeqq
\epf

\bpf[Proof of Lemma \ref{Lem:TriangleCaseLem3}]
Let $\Delta$ be any triangle in $\R^2$ and $T\subset\R^2$ any convex body fulfilling $T\notin \mathcal{T}_\Delta$ and $\vol(T)=2\vol(T)$. Since $\mathcal{T}_\Delta$ is the set of volume-minimizing convex hulls
\beqq \conv\{J\Delta,-J\Delta+t\},\; t\in\R^2,\eeqq
there are $\lambda,\mu\in (0,1]$ such that
\beqq \lambda J\Delta \in F(T)\; \text{ and } \; -\mu J \Delta \in F(T)\eeqq
while at least one of the factors $\lambda,\mu$ can be chosen less than $1$. Referring to Theorem \ref{Thm:Chap7relationship}, the proof of Lemma \ref{Lem:TriangleCaseLem2}, and \cite[Proposition 2.3(i)$\&$(iii)$\&$(iv)]{KruppRudolf2022}, this implies
\begin{align*}
c_{EHZ}(\Delta\times T^*)=\min_{p\in F_3^{cp}(T^*)}\ell_{\Delta}(p) &=\ell_{\Delta}\left(\partial^r(J\Delta)\right)\left[=\ell_{\Delta}\left(\partial^l \left(-J\Delta\right)\right)\right]\\
&>\min\left\{\ell_{\Delta}\left(\partial^r(\lambda J\Delta)\right),\ell_{\Delta}\left(\partial^l \left(-\mu J\Delta\right)\right)\right\}\\
&\geq \min_{p\in F_3^{cp}(T)}\ell_\Delta(p)\\
& = c_{EHZ}(\Delta\times T)
\end{align*}
for every $T^*\in \mathcal{T}_\Delta$, where we denote by
\beqq \partial^r(J\Delta) \; \text{ and } \; \partial^l(-J\Delta)\eeqq
the closed polygonal curves which have $J\Delta$ and $-J\Delta$ as clockwise and counter-clockwise passing orbits, respectively.
\epf

\bpf[Proof of Proposition \ref{Prop:TriangleCaseProp4}]
We proceed in two steps: First, we prove the statement of this proposition for one specially chosen triangle $\Delta$ whose centroid is the origin, secondly, we justify why this is sufficient.

So, let us start with the first step for the triangle $\Delta$ given by the vertices
\beqq \left(-\frac{1}{2},-\frac{1}{2\sqrt{3}}\right),\;\left(\frac{1}{2},-\frac{1}{2\sqrt{3}}\right),\; \left(0,\frac{1}{\sqrt{3}}\right).\eeqq
By construction, $\Delta$'s centroid is the origin. This implies that $\mathcal{T}_\Delta$, which is the set of the volume-minimizing convex hulls
\beqq \conv\{J\Delta,-J\Delta +t\},\quad t\in \R^2,\eeqq
can be given as the set
\beqq \left\{\conv\{J\Delta,-J\Delta+t\} : t\in -J\Delta\right\}\eeqq
(see the description of Figure \ref{img:convexhull}), where the set of convex bodies $T^*$ for which it is claimed that $\mathring{\Delta}\times \mathring{T}^*$ is symplectomorphic to a ball, i.e., the set of convex bodies $T^*\in\mathcal{T}_\Delta$ which are not hexagons, is
\beq \mathcal{T}_{\Delta,\square}=\left\{\conv\{J\Delta,-J\Delta+t\} : t\in -J\partial\Delta\right\}.\label{eq:hhhhhh}\eeq
Now, the goal is to prove that for given $T^*\in\mathcal{T}_{\Delta,\square}$, i.e., for given $t\in -J\partial\Delta$, there  are $a\in\R >0$ and $a_1=a_1(t),a_2=a_2(t)\in [0,a]$ such that
\beq \mathring{\Delta}\times \mathring{T}^*(t) \stackrel{sympl.}{\cong} \diamondsuit(a,a_1(t),a_2(t))\times \square(1)\label{eq:hhhhhtoprove}\eeq
which then by Theorem \ref{Thm:Schlenk} implies
\beqq \mathring{\Delta}\times \mathring{T}^*(t) \stackrel{sympl.}{\cong} B^4_{\sqrt{\frac{a}{\pi}}}.\eeqq

For that, we divide $-J\partial\Delta$ into the three edges $[v_1,v_2]$, $[v_2,v_3]$, and $[v_3,v_1]$, where
\beqq v_1=\left(\frac{1}{2\sqrt{3}},-\frac{1}{2}\right),\;v_2=\left(\frac{1}{2\sqrt{3}},\frac{1}{2}\right),\;v_3=\left(-\frac{1}{\sqrt{3}},0\right),\eeqq
and prove \eqref{eq:hhhhhtoprove} for $t$ out of every single of them.

If $t\in [v_1,v_2]$, then
\beqq t=(t_1,t_2) \text{ with } t_1=\frac{1}{2\sqrt{3}},\; t_2\in \left[-\frac{1}{2},\frac{1}{2}\right].\eeqq
We define $\phi_t$ as linear map given by the matrix
\beqq A_t=\bpm \frac{\sqrt{3}}{2} & t_2 \\ 0 & 1\epm\eeqq
and the translation
\begin{align*}
s=s_1&\times s_2:\R^2(\widetilde{x})\times \R^2(\widetilde{y})\rightarrow\R^2(\widetilde{x})\times \R^2(\widetilde{y}),\\
&(\widetilde{x},\widetilde{y})\mapsto \left(\widetilde{x}+\left(\frac{1}{2},\frac{1}{2\sqrt{3}}\right),\widetilde{y}+\left( \frac{1}{2\sqrt{3}},\frac{1}{2}\right)\right),
\end{align*}
and prove that
\begin{align*}
\left(\left(\phi_t\times \left(\phi_t^T\right)^{-1}\right)\circ s\right)\left(\mathring{\Delta}\times \mathring{T}^*(t)\right)&=\phi_t\left(s_1\left(\mathring{\Delta}\right)\right)\times \left(\phi_t^T\right)^{-1}\left(s_2\left(\mathring{T}^*(t)\right)\right)\\
&=\diamondsuit(a,a_1(t),a_2(t))\times \square(1)
\end{align*}
for
\beqq a=\frac{\sqrt{3}}{2},\; a_1(t)=\frac{\sqrt{3}}{4}+\frac{\sqrt{3}}{2}t_2\in\left[0,\frac{\sqrt{3}}{2}\right],\; a_2=0,\eeqq
where we notice that
\beqq \left(\phi_t\times \left(\phi_t^T\right)^{-1}\right)\circ s\eeqq
is a symplectomorphism as composition of the two symplectomorphism
\beqq \phi_t\times \left(\phi_t^T\right)^{-1} \; \text{ and } \; s,\eeqq
see \cite[Proposition 2.9]{Rudolf2022b}.
\begin{figure}[h!]
\centering
\def\svgwidth{365pt}
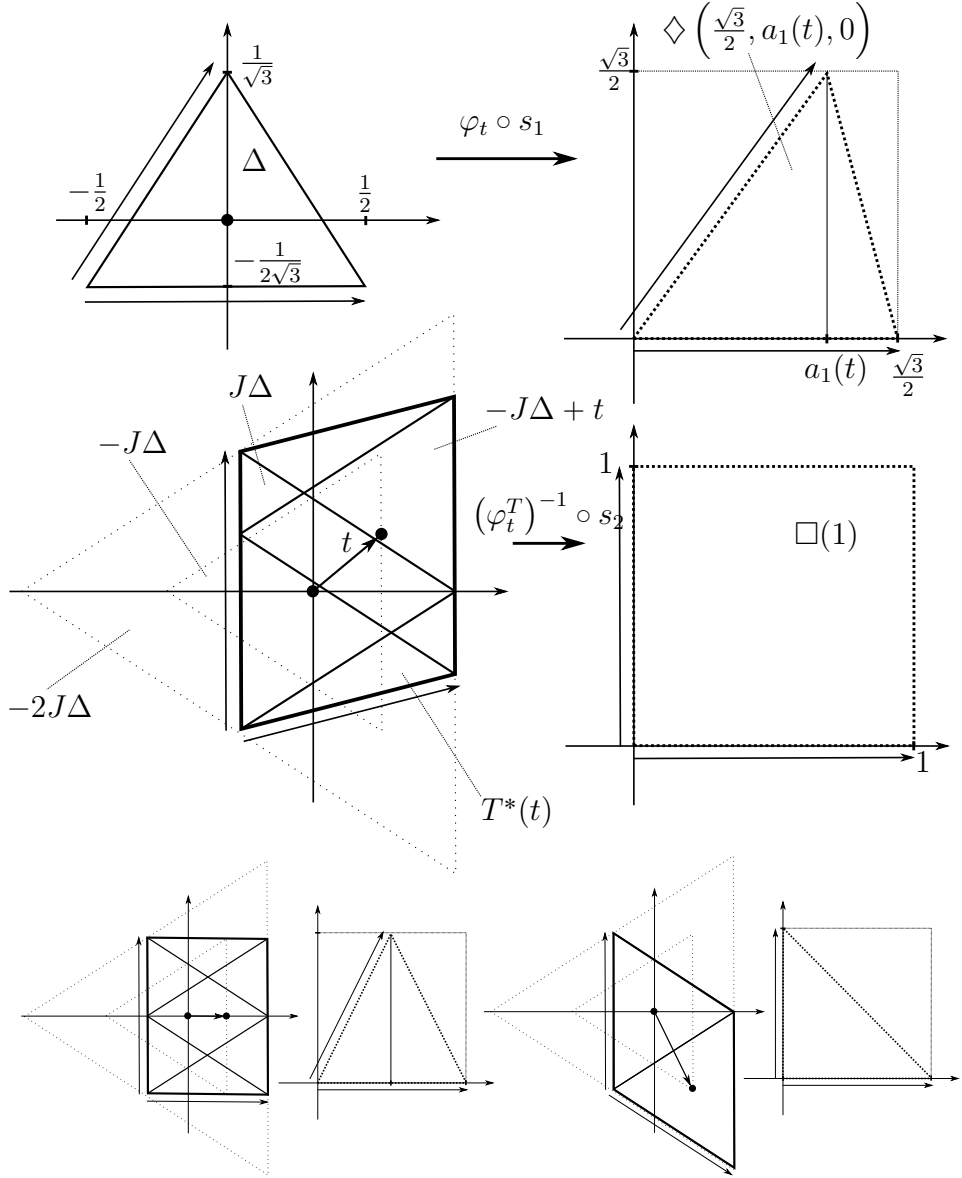
\caption[Illustration of the diffeomorphisms $\phi_t\circ s_1$ and $\left(\phi_t^T\right)^{-1}\circ s_2$ for the first edge]{Illustration of the diffeomorphisms $\phi_t\circ s_1$ and $\left(\phi_t^T\right)^{-1}\circ s_2$ for the case $t\in \left[v_1,v_2\right]$. It is $a_1(t)=\frac{\sqrt{3}}{4}+\frac{\sqrt{3}}{2}t_2$ with $t_2\in\left[-\frac{1}{2},\frac{1}{2}\right]$. In the upper picture $t_2$ is chosen to be $\frac{1}{4}$, in the lower left $0$, and in the lower right $-\frac{1}{2}$. For the upper picture this implies that $a_1(t)$ equals $\frac{3\sqrt{3}}{8}$, for the lower left to $\frac{\sqrt{3}}{4}$, and for the lower right to $0$.}
\label{img:firstedge}
\end{figure}

Indeed (see Figure \ref{img:firstedge}), on the one hand, we have
\beqq \phi_t\left(s_1\left(\mathring{\Delta}\right)\right)=\diamondsuit\left(\frac{\sqrt{3}}{2},\frac{\sqrt{3}}{4}+\frac{\sqrt{3}}{2}t_2,0\right)\eeqq
since $s_1$ translates the lower left corner of $\Delta$ into the origin and
\beqq A_t\bpm 1 \\ 0\epm = \bpm \frac{\sqrt{3}}{2}\\0\epm,\;\; A_t\bpm \frac{1}{2}\\\frac{\sqrt{3}}{2}\epm= \bpm \frac{\sqrt{3}}{4}+\frac{\sqrt{3}}{2}t_2\\ \frac{\sqrt{3}}{2}\epm,\eeqq
and, on the other hand, we have
\beqq \left(\phi_t^T\right)^{-1}\left(s_2\left(\mathring{T}^*(t)\right)\right)=\square(1)\eeqq
since $s_2$ translates the lower left corner of $T^*(t)$ into the origin and with
\beqq \left(A_t^T\right)^{-1}=\bpm \frac{2}{\sqrt{3}} & 0\\ -\frac{2t_2}{\sqrt{3}} & 1\epm\eeqq
we have
\beqq \left(A_t^T\right)^{-1}\bpm \frac{\sqrt{3}}{2}\\t_2\epm = \bpm 1\\0\epm\; \text{ and }\; \left(A_t^T\right)^{-1} \bpm 0\\1\epm=\bpm 0\\1\epm.\eeqq

If $t\in[v_2,v_3]$, then
\beqq t=(t_1,t_2) \text{ with } t_1\in\left[-\frac{1}{\sqrt{3}},\frac{1}{2\sqrt{3}}\right],\; t_2= \frac{t_1}{\sqrt{3}}+\frac{1}{3}.\eeqq
We define $\phi_t$ as linear map given by the matrix
\beqq A_t=\bpm \frac{\sqrt{3}}{2} & \frac{1}{2} \\ t_1- \frac{1}{2\sqrt{3}} & \frac{5}{2}-\frac{t_1}{\sqrt{3}}\epm\eeqq
and the translation
\begin{align*}
s_t&=s_1(t_1)\times s_2:\R^2(\widetilde{x})\times \R^2(\widetilde{y})\rightarrow \R^2(\widetilde{x})\times \R^2(\widetilde{y}),\\
& (\widetilde{x},\widetilde{y})\mapsto \left(\widetilde{x}+\left(\frac{1}{2},\frac{1}{\sqrt{3}}-t_1\right),\widetilde{y}+\left(\frac{1}{2\sqrt{3}},\frac{1}{2}\right)\right),
\end{align*} 
and prove that
\begin{align*}
\left(\left(\phi_t\times \left(\phi_t^T\right)^{-1}\right)\circ s_t\right)\left(\mathring{\Delta}\times \mathring{T}^*(t)\right)&=\phi_t\left(s_1(t_1)\left(\mathring{\Delta}\right)\right)\times \left(\phi_t^T\right)^{-1}\left(s_2\left(\mathring{T}^*(t)\right)\right)\\
&=\diamondsuit(a,a_1(t),a_2(t))\times \square(1)
\end{align*}
for
\beqq a=\frac{\sqrt{3}}{2},\; a_1=\frac{\sqrt{3}}{2},\; a_2(t)=\frac{1}{2\sqrt{3}}-t_1.\eeqq
\begin{figure}[h!]
\centering
\def\svgwidth{365pt}
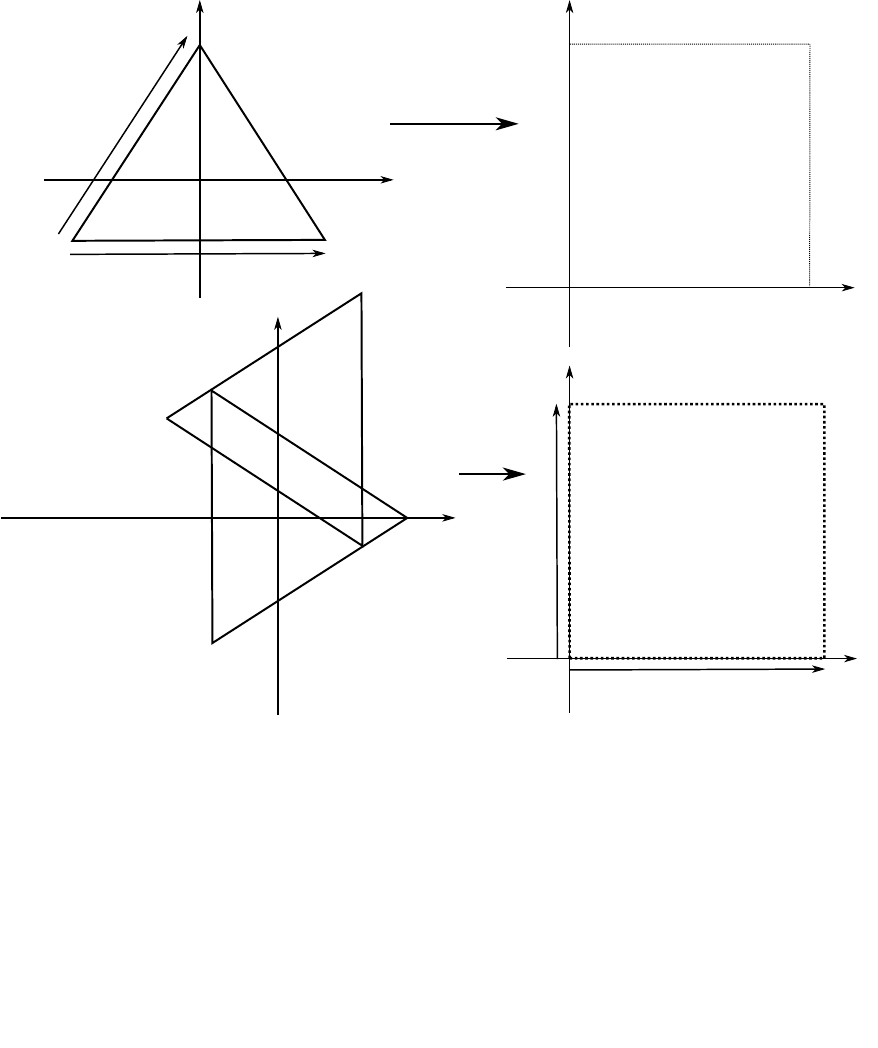
\caption[Illustration of the diffeomorphisms $\phi_t\circ s_1(t_1)$ and $\left(\phi_t^T\right)^{-1}\circ s_2$ for the second edge]{Illustration of the diffeomorphisms $\phi_t\circ s_1(t_1)$ and $\left(\phi_t^T\right)^{-1}\circ s_2$ for the case $t\in \left[v_2,v_3\right]$. It is $a_2(t)=\frac{1}{2\sqrt{3}}-t_1$ with $t_1\in[-\frac{1}{\sqrt{3}},\frac{1}{2\sqrt{3}}]$. In the upper picture $t_1$ is chosen to be $\frac{1}{8\sqrt{3}}$ (which is at $\frac{3}{4}$ of the intervall $[-\frac{1}{\sqrt{3}},\frac{1}{2\sqrt{3}}]$), in the lower left $-\frac{1}{2\sqrt{3}}$, and in the lower right $-\frac{1}{4\sqrt{3}}$ (which is at $\frac{1}{2}$ of the intervall $[-\frac{1}{\sqrt{3}},\frac{1}{2\sqrt{3}}]$).}
\label{img:secondedge}
\end{figure}
Indeed (see Figure \ref{img:secondedge}), on the one hand, we have
\beqq \phi_t\left(s_1(t_1)\left(\mathring{\Delta}\right)\right)=\diamondsuit\left(\frac{\sqrt{3}}{2},\frac{\sqrt{3}}{2},\frac{1}{2\sqrt{3}}-t_1\right)\eeqq
since $s_1(t_1)$ translates the lower left corner of $\Delta$ to $\left(0,\frac{1}{2\sqrt{3}}-t_1\right)$ and
\beqq A_t\bpm 1 \\ 0\epm = \bpm \frac{\sqrt{3}}{2}\\ t_1\epm,\;\; A_t\bpm \frac{1}{2}\\\frac{\sqrt{3}}{2}\epm = \bpm \frac{\sqrt{3}}{2}\\t_1+\frac{\sqrt{3}}{2}\epm,\eeqq
and, on the other hand, we have
\beqq \left(\phi_t^T\right)^{-1}\left(s_2\left(\mathring{T}^*(t)\right)\right)=\square(1)\eeqq
since $s_2$ translates the lower left corner of $T^*(t)$ into the origin and with
\beqq \left(A_t^T\right)^{-1}=\bpm \frac{2t_1}{3}+\frac{2}{\sqrt{3}} & -\frac{2t_1}{\sqrt{3}} \\ t_1 & \frac{t_1}{\sqrt{3}}+1\epm\eeqq
we have
\beqq \left(A_t^T\right)^{-1}\bpm \frac{\sqrt{3}}{2}\\\frac{1}{2}\epm = \bpm 1\\0\epm\; \text{ and }\; \left(A_t^T\right)^{-1}\bpm t_1-\frac{\sqrt{3}}{2} \\ \frac{t_1}{\sqrt{3}}+\frac{1}{2}\epm = \bpm 0 \\ 1\epm.\eeqq

If $t\in[v_3,v_1]$, then
\beqq t=(t_1,t_2) \text{ with } t_1\in\left[-\frac{1}{\sqrt{3}},\frac{1}{2\sqrt{3}}\right],\; t_2=-\frac{t_1}{\sqrt{3}}-\frac{1}{3}.\eeqq
We define $\phi_t$ as linear map given by the matrix
\beqq A_t=\bpm -t_1+\frac{1}{2\sqrt{3}} & \frac{5}{6}+\frac{t_1}{\sqrt{3}} \\ -\frac{\sqrt{3}}{2} & \frac{1}{2}\epm\eeqq
and the translation
\begin{align*}
s=&s_1\times s_2:\R^2(\widetilde{x})\times \R^2(\widetilde{y})\rightarrow \R^2(\widetilde{x})\times \R^2(\widetilde{y}),\\
&(\widetilde{x},\widetilde{y})\mapsto \left(\widetilde{x}+\left( \frac{1}{2}, \frac{2}{\sqrt{3}}\right),\widetilde{y}+\left( -t_1-\frac{1}{2\sqrt{3}},\frac{t_1}{\sqrt{3}}+\frac{5}{6}\right)\right)
\end{align*}
and prove that
\begin{align*}
\left(\left(\phi_t\times \left(\phi_t^T\right)^{-1}\right)\circ s\right)\left(\mathring{\Delta}\times \mathring{T}^*(t)\right)&=\phi_t\left(s_1(\mathring{\Delta})\right)\times \left(\phi_t^T\right)^{-1}\left(s_2\left(\mathring{T}^*(t)\right)\right)\\
&=\diamondsuit(a,a_1(t),a_2)\times \square(1)
\end{align*}
for
\beqq a=\frac{\sqrt{3}}{2},\; a_1(t)=-t_1+\frac{1}{2\sqrt{3}}\in\left[0,\frac{\sqrt{3}}{2}\right],\; a_2=\frac{\sqrt{3}}{2}.\eeqq
\begin{figure}[h!]
\centering
\def\svgwidth{365pt}
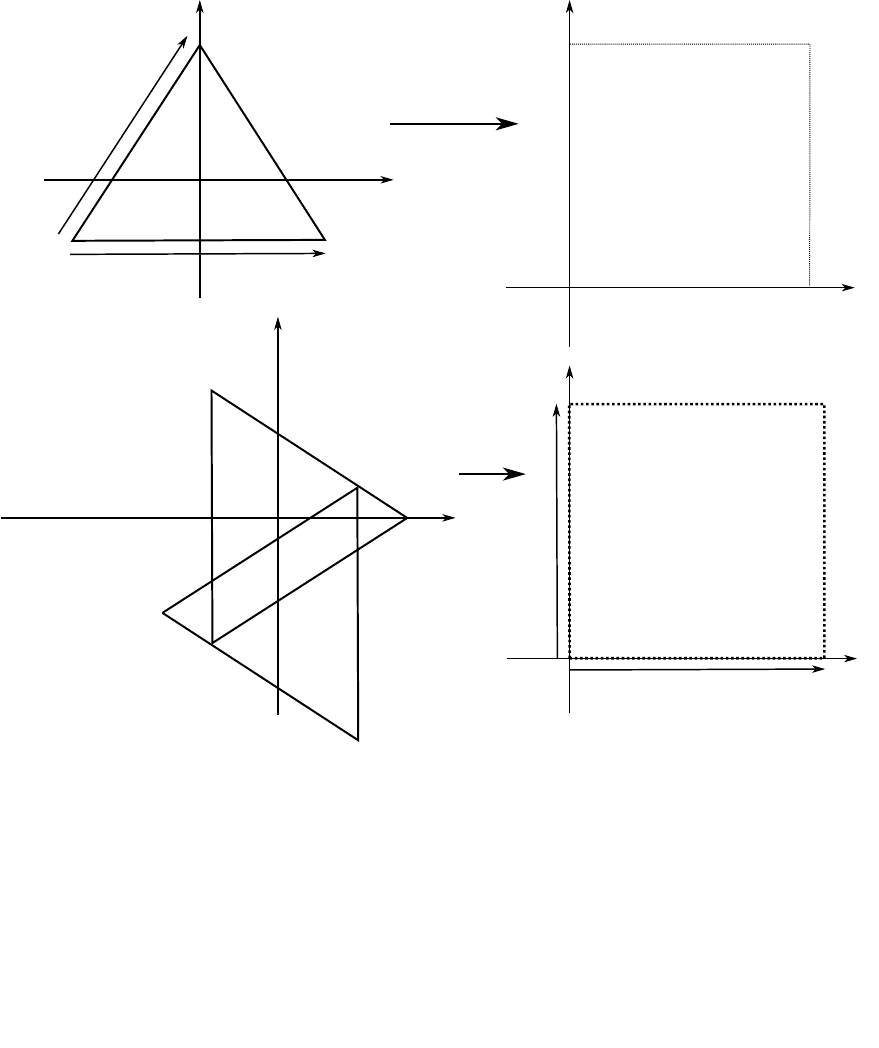
\caption[Illustration of the diffeomorphisms $\phi_t\circ s_1$ and $\left(\phi_t^T\right)^{-1}\circ s_2$ for the third edge]{Illustration of the diffeomorphisms $\phi_t\circ s_1$ and $\left(\phi_t^T\right)^{-1}\circ s_2$ for the case $t\in \left[v_3,v_1\right]$. It is $a_1(t)=-t_1+\frac{1}{2\sqrt{3}}$ with $t_1\in[-\frac{1}{\sqrt{3}},\frac{1}{2\sqrt{3}}]$. In the upper picture $t_1$ is chosen to be $\frac{1}{8\sqrt{3}}$ (which is at $\frac{3}{4}$ of the intervall $[-\frac{1}{\sqrt{3}},\frac{1}{2\sqrt{3}}]$), in the lower left $-\frac{1}{2\sqrt{3}}$, and in the lower right $-\frac{1}{4\sqrt{3}}$ (which is at $\frac{1}{2}$ of the intervall $[-\frac{1}{\sqrt{3}},\frac{1}{2\sqrt{3}}]$).}
\label{img:thirdedge}
\end{figure}
Indeed (see Figure \ref{img:thirdedge}), on the one hand, we have
\beqq \phi_t\left(s_1\left(\mathring{\Delta}\right)\right)=\diamondsuit\left(\frac{\sqrt{3}}{2},-t_1+\frac{1}{2\sqrt{3}},\frac{\sqrt{3}}{2}\right)\eeqq
since $t_1$ translates the lower left vertex of $\Delta$ to $\left(0,\frac{\sqrt{3}}{2}\right)$ and
\beqq A_t\bpm 1\\0\epm = \bpm -t_1 + \frac{1}{2\sqrt{3}} \\ -\frac{\sqrt{3}}{2}\epm,\;\; A_t\bpm \frac{1}{2}\\\frac{\sqrt{3}}{2}\epm = \bpm \frac{\sqrt{3}}{2} \\ 0\epm,\eeqq
and, on the other hand, we have
\beqq \left(\phi_t^T\right)^{-1}\left(s_2\left(\mathring{T}^*(t)\right)\right)=\square(1)\eeqq
since $t_2$ translates the lower vertex of $T^*(t)$ into the origin and with
\beqq \left(A_t^T\right)^{-1}=\bpm \frac{1}{\sqrt{3}} & 1 \\ -\frac{5}{3\sqrt{3}}-\frac{2t_1}{3} & -\frac{2t_1}{\sqrt{3}}+\frac{1}{3}\epm\eeqq
we have
\beqq \left(A_t^T\right)^{-1} \bpm -t_1 + \frac{1}{2\sqrt{3}} \\ \frac{t_1}{\sqrt{3}}+\frac{5}{6}\epm = \bpm 1 \\ 0\epm \; \text{ and }\; \left(A_t^T\right)^{-1} \bpm -\frac{\sqrt{3}}{2}\\\frac{1}{2}\epm = \bpm 0\\ 1\epm.\eeqq

\underline{Note in proof:} When comparing Figures \ref{img:firstedge} and \ref{img:thirdedge} for the choice $t=(\frac{1}{2\sqrt{3}},-\frac{1}{2})$, respectively, one can notice that there are two different Lagrangian products $\diamondsuit \times \square$ (differing by the different $\diamondsuit$s, i.e., $\diamondsuit(\frac{\sqrt{3}}{2},0,\frac{\sqrt{3}}{2})$ and $\diamondsuit(\frac{\sqrt{3}}{2}, 0, 0)$). However, this is not a problem since one can easily check that the symplectomorphism $J\times \left(J^T\right)^{-1}=J\times J$ maps these two Lagrangian products on each other (up to a translation of $\square$).

This completes the first step. Now, as second step, let us justify why this is sufficient: Let $\widetilde{\Delta}$ be any other triangle in $\R^2$ which, without loss of generality, satisfies  $\vol(\widetilde{\Delta})=\vol(\Delta)$ and whose centroid is the origin (note that translations are symplectomorphisms--see the proof of \cite[Proposition 2.8]{Rudolf2022b}) such that
\beqq \widetilde{T}=\conv\left\{J\widetilde{\Delta},-J\widetilde{\Delta} + \widetilde{t}\right\}\eeqq
is in $\mathcal{T}_{\widetilde{\Delta},\square}$. It is enough to show that there is a symplectomorphism $\psi$ and a $T\in\mathcal{T}_{\Delta}$ with
\beqq T=\conv\{J\Delta,-J\Delta+t\}\;\text{ and } \; t\in -J\partial\Delta\eeqq
such that
\beq \psi \left(\mathring{\Delta} \times \mathring{T}\right) = \mathring{\widetilde{\Delta}}\times \mathring{\widetilde{T}}.\label{eq:symplgoal}\eeq
Let $\phi$ be the linear transformation for which
\beqq \phi(\Delta)=\widetilde{\Delta}\eeqq
and let
\beqq T=\conv\left\{J\Delta,-J\Delta + \phi^T\left(\widetilde{t}\right)\right\}.\eeqq
Then,
\beqq \psi=\phi \times \left(\phi^T\right)^{-1}\eeqq
is a symplectomorphism fulfilling \eqref{eq:symplgoal}.

Indeed, using the linearity of $\phi$, we calculate
\begin{align*}
\psi(\Delta\times T)&=\phi(\Delta)\times \left(\phi^T\right)^{-1}\left(\conv\left\{J\Delta,-J\Delta + \phi^T\left(\widetilde{t}\right)\right\}\right)\\
&=\widetilde{\Delta} \times \conv\left\{\left(\phi^T\right)^{-1}(J\Delta),\left(\phi^T\right)^{-1}(-J\Delta) +\widetilde{t}\right\}\\
& \stackrel{(\star)}{=} \widetilde{\Delta}\times \conv\left\{J\phi(\Delta),-J\phi(\Delta)+\widetilde{t}\right\}\\
&= \widetilde{\Delta} \times \conv\left\{J\widetilde{\Delta},-J\widetilde{\Delta}+\widetilde{t}\right\}\\
&= \widetilde{\Delta}\times \widetilde{T},
\end{align*}
where for equality $(\star)$, we used
\beq \phi^T J \phi = \det(\phi) J \; \text{ and } \; \phi^T (-J) \phi = \det(\phi)(-J)\label{eq:detsympl}\eeq
and $\det(\phi)=1$. Here, \eqref{eq:detsympl} holds for general linear transformation $\phi$ because one can easily check that both $\phi^T (\pm J) \phi$ as well as $\det(\phi)(\pm J)$ map $(a,b)\in\R^2$ to
\beqq \left(\pm\det(\phi)b,\mp\det(\phi)a\right)=\det(\phi)(\pm J)(a,b)\in\R^2\eeqq
for all linear transformations $\phi$. It remains to show that $T$ is in $\mathcal{T}_\Delta$ such that
\beq \phi^T\left(\widetilde{t}\right)\in -J\partial\Delta.\label{eq:sympltriangl1}\eeq
Applying $\left(\phi^T\right)^{-1}$ on both sides of \eqref{eq:sympltriangl1} and using \eqref{eq:detsympl}, implies that \eqref{eq:sympltriangl1} is equivalent to
\beqq \widetilde{t}\in \left(\phi^T\right)^{-1}(-J\partial\Delta)=-J\phi(\partial\Delta)=-J\partial\widetilde{\Delta}\eeqq
which holds by assumption. Therefore, \eqref{eq:sympltriangl1} is in fact satisfied.
\epf

\section{The parallelogram-case}\label{Sec:parallelogram}

As we have already noted in the introduction, by using suitable affine transformations $\phi$ and the fact that $\phi\times\left(\phi^T\right)^{-1}$ are symplectomorphisms, the square-configurations $\square\times T$ can be easily lifted to parallelogram-configurations. So, in what follows, we only have to treat the square-configurations.

First, we obviously note that
\beqq \vol\left(\diamondsuit(a_1,a_2)\right)=\frac{1}{2}\quad \forall a_1,a_2\in[0,1].\eeqq

Then, we proceed in four steps:

\blem\label{Lem:squarestep0}
Let $\square$ be any square in $\R^2$. Then, we have
\beq c_{EHZ}(\square\times\diamondsuit(a_1,a_2))=\min_{q\in M_2(\square,\diamondsuit(a_1,a_2))}\ell_{\diamondsuit(a_1,a_2)}(q)\label{eq:squarestep0}\eeq
for all $a_1,a_2\in[0,1]$.
\elem

\blem\label{Lem:squarestep1}
Let $\square$ be any square in $\R^2$. Then, we have
\beqq \vol\left(\square\times \diamondsuit(a_1,a_2)\right)=\frac{c_{EHZ}(\square\times \diamondsuit(a_1,a_2))^2}{2}\eeqq
for all $a_1,a_2\in[0,1]$.
\elem

\blem\label{Lem:squarestep2}
Let $\square$ be any square in $\R^2$. Then, we have
\beqq c_{EHZ}(\square\times\diamondsuit(a_1,a_2))>c_{EHZ}(\square\times T),\quad a_1,a_2\in[0,1],\eeqq
for all convex bodies $T\subset\R^2$ with 
\beqq \vol(T)=\frac{1}{2}\; \text{ and } \; T\neq \diamondsuit(a_1',a_2')+t \;\; \forall a_1',a_2'\in[0,1],t\in\R^2.\eeqq
\elem

\blem\label{Lem:squarestep3}
Let $\square$ be any square in $\R^2$ with side length $a>0$. Then, we have
\beq \mathring{\square}\times \mathring{\diamondsuit}(a_1,a_2) \stackrel{sympl.}{\cong} B^4_{\sqrt{\frac{a^3}{\pi}}}\label{eq:squarestep3}\eeq
for all $a_1,a_2\in[0,1]$.
\elem

Combining Lemmata \ref{Lem:squarestep0}, \ref{Lem:squarestep1}, \ref{Lem:squarestep2}, and \ref{Lem:squarestep3}, implies Theorem \ref{Thm:answer1} for the case when $Q$ in \eqref{eq:answer1} is any parallelogram in $\R^2$ and Theorem \ref{Thm:answer2}(ii) and \ref{Thm:answer3}(ii).

We remark that Lemma \ref{Lem:squarestep1} is a direct consequence of Lemma \ref{Lem:squarestep3} (since any ball in $\R^4$ is an equality case of the $4$-dimensional Viterbo conjecture--which is invariant under symplectomorphisms). Nevertheless, we will give an independet proof of Lemma \ref{Lem:squarestep1} since it prepares a conceptual understanding of the calculation of $c_{EHZ}(\square\times\diamondsuit(a_1,a_2))$ wich is important for the proof of Lemma \ref{Lem:squarestep2}.

So, let us prove the above lemmata one after the other.

\bpf[Proof of Lemma \ref{Lem:squarestep0}]
Recalling Theorem \ref{Thm:Chap7relationship}, we have
\beqq c_{EHZ}(\square\times\diamondsuit(a_1,a_2))=\min_{q\in M_3(\square,\diamondsuit(a_1,a_2))}\ell_{\diamondsuit(a_1,a_2)}(q).\eeqq
We show that for every $(\square,\diamondsuit(a_1,a_2))$-Minkowski billiard trajectory $q$ with $3$ bouncing points, there is a $\widetilde{q}\in M_2(\square,\diamondsuit(a_1,a_2))$ with
\beqq \ell_{\diamondsuit(a_1,a_2)}(\widetilde{q})\leq  \ell_{\diamondsuit(a_1,a_2)}(q).\eeqq
This would imply \eqref{eq:squarestep0}.

So, let $q$ be a $(\square,\diamondsuit(a_1,a_2))$-Minkowski billiard trajectory with $3$ bouncing points $q_1,q_2,q_3$. Then, we conclude that $\pi_1(q)$ and $\pi_2(q)$, where $\pi_i$, $i\in\{1,2\}$, projects $\R^2$ onto its $i$-th coordinate axis, reflect the full horizontal and vertical diameter of $\square$, respectively. Otherwise, the convex cone generated by every choice of normal vectors $n_\square(q_1)$, $n_\square(q_2)$, and $n_\square(q_3)$ would not contain the origin in its interior, i.e., the convex hull of $n_\square(q_1)$, $n_\square(q_2)$, and $n_\square(q_3)$ would not contain the origin. This would be a contradiction to what has been shown within the proof of \cite[Proposition 3.9]{KruppRudolf2022}. Without loss of generality, we assume
\beqq \pi_1(q)=[\pi_1(q_1),\pi_1(q_2)].\eeqq
Then, the closed polygonal curve $q^*$ defined by the two vertices $q_1$ and $q_2$ is in $F_2(\square)$ and has less or equal $\ell_{\diamondsuit(a_1,a_2)}$-length than $q$ (see \cite[Proposition 2.3(i)]{KruppRudolf2022}). Applying the addition of Theorem \ref{Thm:Chap7relationship}, i.e.,
\beqq \min_{q\in F_2^{cp}(\square)} \ell_{\diamondsuit(a_1,a_2)}(q) = \min_{q\in M_2(\square,\diamondsuit(a_1,a_2))} \ell_{\diamondsuit(a_1,a_2)}(q),\eeqq
yields the existence of a $\widetilde{q}\in M_2(\square,\diamondsuit(a_1,a_2))$ with
\beqq \ell_{\diamondsuit(a_1,a_2)}(\widetilde{q})\leq  \ell_{\diamondsuit(a_1,a_2)}(q^*) \leq \ell_{\diamondsuit(a_1,a_2)}(q).\eeqq
\epf

We remark that the proof of Lemma \ref{Lem:squarestep0} shows that the truth of \eqref{eq:squarestep0} is not restricted to Lagrangian configurations $\square\times \diamondsuit(a_1,a_2)$. One could replace $\diamondsuit(a_1,a_2)$ by any body $T \subset\R^2$ when guaranteeing its convexity. This will be of importance for the proof of Lemma \ref{Lem:squarestep2}.

\bpf[Proof of Lemma \ref{Lem:squarestep1}]
Let $a$ be the side length of $\square$. Then, we first note that
\beqq \vol(\square\times\diamondsuit(a_1,a_2))=\vol(\square)\vol(\diamondsuit(a_1,a_2))=\frac{a^2}{2}\eeqq
for all $a_1,a_2\in[0,1]$. Referring to Lemma \ref{Lem:squarestep0}, it therefore suffices to prove
\beq \min_{q\in M_2(\square,\diamondsuit(a_1,a_2))} \ell_{\diamondsuit(a_1,a_2)}(q)=a\label{eq:squarestep1goal}\eeq
for all $a_1,a_2\in[0,1]$.

In fact, let $q=(q_1,q_2)$ be a closed $(\square,\diamondsuit(a_1,a_2))$-Minkowski billiard trajectory.

\begin{figure}[h!]
\centering
\def\svgwidth{295pt}
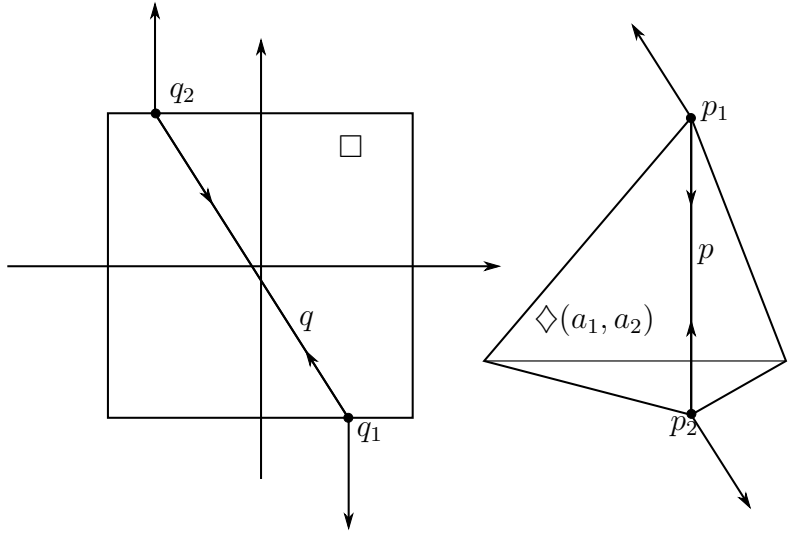
\caption[Illustration of $\square\times \diamondsuit(a_1,a_2)$ for the case $q_1$ and $q_2$ lie on the interiors of two vertically/horizontally opposite edges of $\square$.]{Illustration of $\square\times \diamondsuit(a_1,a_2)$ for the case $q_1$ and $q_2$ lie on the interiors of two vertically opposite edges of $\square$. The dual billiard trajectory $p=(p_1,p_2)$ in $\diamondsuit(a_1,a_2)$ reflects the vertical diameter of $\diamondsuit(a_1,a_2)$.}
\label{img:sequence_1}
\end{figure}

If $q_1$ and $q_2$ lie on the interiors of two vertically/horizontally opposite edges of $\square$, then the Minkowski billiard reflection rule (see \eqref{eq:Minkreflectionrulepre}) implies that the closed dual billiard trajectory $p=(p_1,p_2)$ in $\diamondsuit(a_1,a_2)$ reflects the vertical/horizontal diameter of $\diamondsuit(a_1,a_2)$ (see Figure \ref{img:sequence_1}).

Using \cite[Proposition 2.2]{KruppRudolf2022} and the properties of the inner product, it therefore follows
{\allowdisplaybreaks[0] \begin{align*}
\ell_{\diamondsuit(a_1,a_2)}(q) & = \mu_{\diamondsuit(a_1,a_2)^\circ}(q_2-q_1) + \mu_{\diamondsuit(a_1,a_2)^\circ}(q_1-q_2)\\
&=\langle q_2-q_1,p_1\rangle + \langle q_1-q_2,p_2\rangle\\
&=\langle q_2-q_1,p_1-p_2\rangle\\
&=a.
\end{align*}}%

If $q_1$ and $q_2$ do not lie on the interiors of two vertically/horizontally opposite edges of $\square$, then $q_1$ and $q_2$ lie at least somewhere on vertically/horizontally opposite edges of $\square$--with $q_1$, $q_2$, or both as vertices of $\square$--all other configurations are excluded by what has been shown in the proof of \cite[Proposition 3.9]{KruppRudolf2022}.

\begin{figure}[h!]
\centering
\def\svgwidth{325pt}
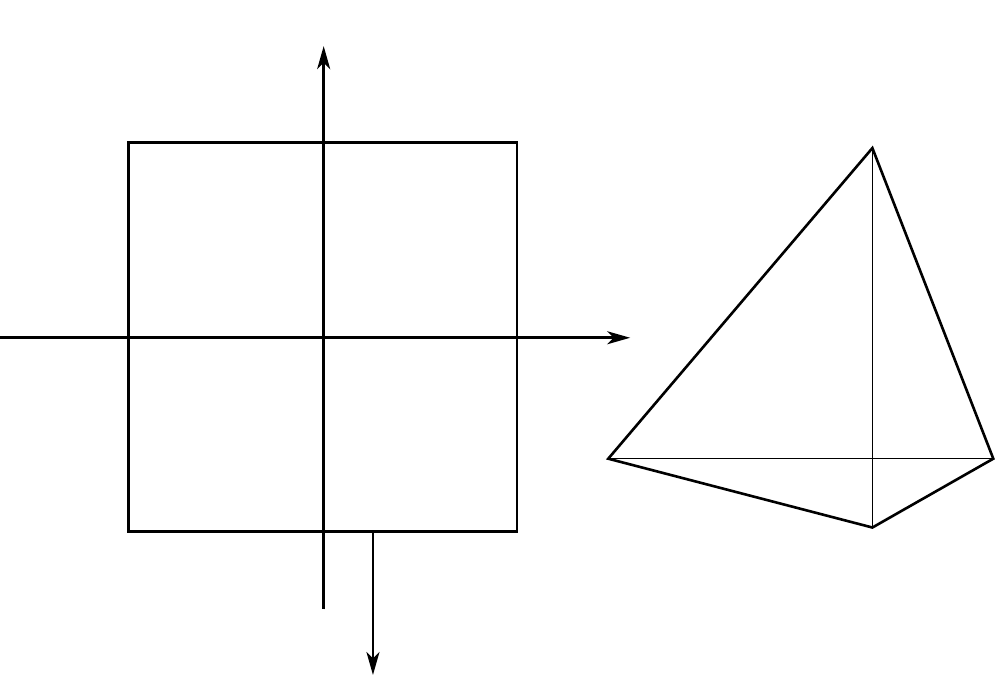
\caption[Illustration of $\square\times \diamondsuit(a_1,a_2)$ for the case $q_1$ and $q_2$ lie on two vertically/horizontally opposite edges of $\square$, where one of them is a vertex.]{Illustration of $\square\times \diamondsuit(a_1,a_2)$ for the case $q_1$ and $q_2$ lie on two vertically opposite edges of $\square$, where $q_2$ is a vertex. The closed dual billiard trajectory $p=(p_1,p_2)$ in $\diamondsuit(a_1,a_2)$ equals the closed dual billiard trajectory which corresponds to a Minkowski billiard trajectory having both--$q_1$ and $q_2$--on the interiors of the edges of $\square$.}
\label{img:problemsolving1}
\end{figure}

If not both--$q_1$ and $q_2$--are vertices of $\square$ (see Figure \ref{img:problemsolving1}), then it follows by the Minkowski billiard reflection rule that $q$'s corresponding closed dual billiard trajectory in $\diamondsuit(a_1,a_2)$ equals the one which corresponds to Minkowski billiard trajectories with both bouncing points on the interiors of the respective edges of $\square$. In this case, one also gets $\ell_{\diamondsuit(a_1,a_2)}(q)=a$.

\begin{figure}[h!]
\centering
\def\svgwidth{360pt}
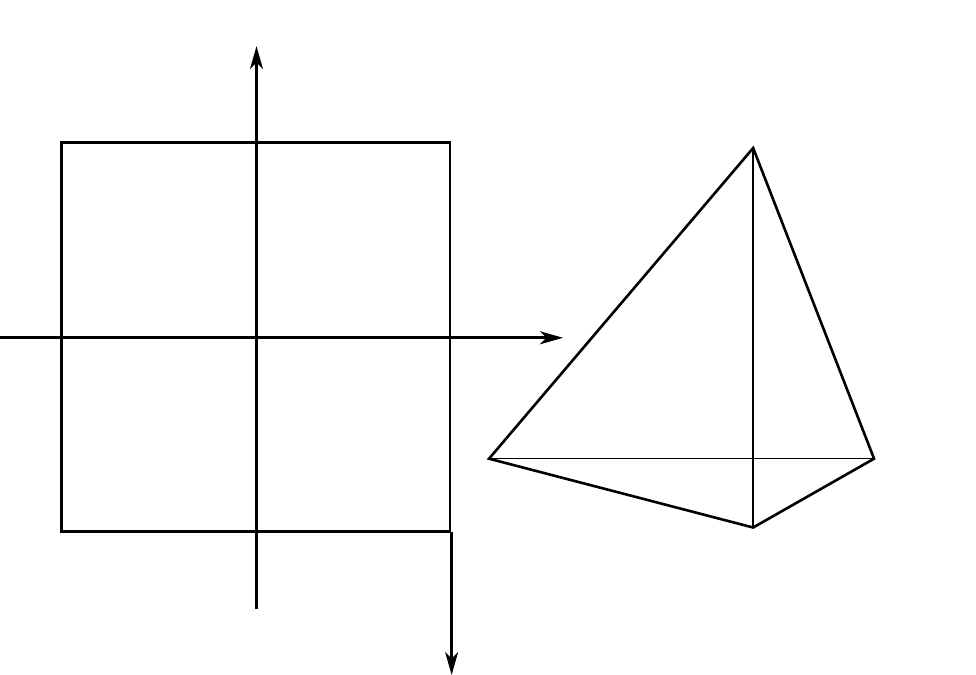
\caption[Illustration of $\square\times \diamondsuit(a_1,a_2)$ for the case $q_1$ and $q_2$ lie on two vertically/horizontally opposite edges of $\square$, where both are vertices.]{Illustration of $\square\times \diamondsuit(a_1,a_2)$ for the case $q_1$ and $q_2$ lie on two vertically opposite edges of $\square$, where both are vertices of $\square$. The closed dual billiard trajectory $p=(p_1,p_2)$ in $\diamondsuit(a_1,a_2)$ does not coincide with a closed dual billiard trajectory $p'=(p_1',p_2')$ which corresponds to a closed Minkowski billiard trajectory $q'=(q_1',q_2')$ with both bouncing points on the interiors of two opposite edges of $\square$.}
\label{img:problemsolving2}
\end{figure}

If both--$q_1$ as well as $q_2$--are opposite vertices of $\square$, then it can happen that $q$'s corresponding closed dual billiard trajectory $p$ in $\diamondsuit(a_1,a_2)$ does not coincide with a closed dual billiard trajectory $p'$ which corresponds to a closed Minkowski billiard trajectory $q'$ that has both bouncing points on the interiors of two opposite edges of $\square$. In Figure \ref{img:problemsolving2}--which, without loss of generality, can serve as a pattern of these cases--one can see that $q$ and $q'$ enclose a sequence of closed polygonal curves $q^n=(q_1^n,q_2^n)$ in $F(\square)$ with $q_1^n=q_1$ and $q_2^n\rightarrow q_2$ ($n\rightarrow\infty$). These closed polygonal curves cannot satisfy the Minkowsi billiard reflection rule what implies that $\ell_{\diamondsuit(a_1,a_2)}(q^n)\geq a$ (due to Theorem \ref{Thm:Chap7relationship}). By the continuity of the length-functional, this implies that $q$ (as limit of $q^n$ with respect to the Hausdorff topology) cannot have $\ell_{\diamondsuit(a_1,a_2)}$-length less than $a$.
\epf

\bpf[Proof of Lemma \ref{Lem:squarestep2}]
Let $T\subset\R^2$ be a convex body with
\beqq \vol(T)=\frac{1}{2}\; \text{ and } \; T\neq \diamondsuit(a_1',a_2')+t \;\; \forall a_1',a_2'\in[0,1],t\in\R^2.\eeqq
The set
\beq \left\{\diamondsuit(a_1,a_2) : a_1,a_2\in[0,1]\right\}\label{eq:squarestep2}\eeq
characterizes the volume-minimizing convex bodies in $\R^2$ (their volume is $\frac{1}{2}$) that contain translates of
\beqq \left[\left(-\frac{1}{2},0\right),\left(\frac{1}{2},0\right)\right] \; \text{ and } \; \left[\left(0,-\frac{1}{2}\right),\left(0,\frac{1}{2}\right)\right].\eeqq
By definition, $T$ is not a member of \eqref{eq:squarestep2}, which implies that there is a $\lambda <1$ such that
\beq \lambda  \left[\left(-\frac{1}{2},0\right),\left(\frac{1}{2},0\right)\right] \in F(T)\; \text{ or } \; \lambda \left[\left(0,-\frac{1}{2}\right),\left(0,\frac{1}{2}\right)\right]\in F(T).\label{eq:suqarestep2twointervals}\eeq
From this, we conclude
\beqq c_{EHZ}(\square \times T) = \min_{q\in M_2(\square,T)}\ell_T(q)=\min_{q\in F_2(\square)}\ell_T(q)=\min_{p\in F_2(T)}\ell_{\square}(p) \leq \lambda a < a, \eeqq
where in the first equality, we used the remark beyond the proof of Lemma \ref{Lem:squarestep0}, in the second and third, Theorem \ref{Thm:Chap7relationship}, and in the second to last inequality, the fact that the intervals in \eqref{eq:suqarestep2twointervals} can be understood as closed polygonal curves with two vertices and $\ell_\square(p)$ can be calculated by again using \cite[Proposition 2.2]{KruppRudolf2022} and the properties of the inner product. Referring
to Lemma \ref{Lem:squarestep0} and \eqref{eq:squarestep1goal}, we therefore conclude
\beqq c_{EHZ}(\square \times T) < c_{EHZ}(\square\times\diamondsuit(a_1,a_2)).\eeqq
\epf

\bpf[Proof of Lemma \ref{Lem:squarestep3}]
Let $c=(c_1,c_2)\in\R^2$ be the center of $\square$. For $t\in\R^2$, we define the translation
\beq s_t: \R^2 \rightarrow \R^2, \quad x\mapsto x-t.\label{eq:translationdefinitionst}\eeq
Then, for $u=(\frac{1}{2},\frac{1}{2})$ and considering the Lagrangian splitting \eqref{eq:Lagrangiansplitting}, we define the map
\beqq \psi:\R^2(x)\times\R^2(y)\rightarrow\R^2(x)\times\R^2(y)\eeqq
\beqq \psi = \left(\mathbb{1}_{\R^2}\times s_{-u}\right)\circ J \circ \left(s_c\times \mathbb{1}_{\R^2}\right)\eeqq
and claim that $\psi$ is a symplectomorphism with
\beq \psi\left(\mathring{\square}\times \mathring{\diamondsuit}(a_1,a_2)\right)=\diamondsuit\left(a,\frac{a_1}{a},\frac{a_2}{a}\right)\times \square(1)\stackrel{\text{sympl.}}{\cong}B^4_{\sqrt{\frac{a^3}{\pi}}}.\label{eq:squarestep3proof0}\eeq
This would prove \eqref{eq:squarestep3}.

Obviously, $\psi$--as composition of symplectomorphisms--is a symplectomorpism. In order to show \eqref{eq:squarestep3proof0}, we notice (considering the notation in Section \ref{Sec:Symplectomorphismstoball}: $\mathring{\diamondsuit}(a_1,a_2)=\diamondsuit(1,a_1,a_2)$ and $\square(1)$ is the open square centred at $(\frac{1}{2},\frac{1}{2})$ with side length $1$) that
\beqq \frac{1}{a}\psi\left(\mathring{\square}\times \mathring{\diamondsuit}(a_1,a_2)\right)=\left(\frac{1}{a}\mathring{\diamondsuit}(a_1,a_2)\right)\times \square(1)=\diamondsuit\left(a,\frac{a_1}{a},\frac{a_2}{a}\right)\times \square(1)\eeqq
(we used $-s_c(\square)=s_c(\square)$) which due to Theorem \ref{Thm:Schlenk} is symplectomorphic to $B^4_{\sqrt{\frac{a}{\pi}}}$. Therefore, it follows
\beqq \psi\left(\mathring{\square}\times \mathring{\diamondsuit}(a_1,a_2)\right) \stackrel{\text{sympl.}}{\cong} a B^4_{\sqrt{\frac{a}{\pi}}}=B^4_{\sqrt{\frac{a^3}{\pi}}}\eeqq
and consequently
\beqq \mathring{\square}\times \mathring{\diamondsuit}(a_1,a_2) \stackrel{\text{sympl.}}{\cong} B^4_{\sqrt{\frac{a^3}{\pi}}}.\eeqq
\epf

\section{The convex-quadrilateral-case}\label{Sec:quadrilateral}

In what follows, we will regard a general convex quadrilateral $Q$ in $\R^2$ as image of an accordingly chosen $\diamondsuit(a_1,a_2)=\diamondsuit(a_1,a_2)(Q)$, $a_1,a_2\in[0,1]$, under an affine transformation $\phi=\phi(Q)$. $Q$ will be uniquely described by the choice of $\diamondsuit(a_1,a_2)(Q)$ and $\phi(Q)$ (and the other way around: $Q$ uniquely determines the corresponding $\diamondsuit(a_1,a_2)$ and $\phi$).

\begin{figure}[h!]
\centering
\def\svgwidth{300pt}
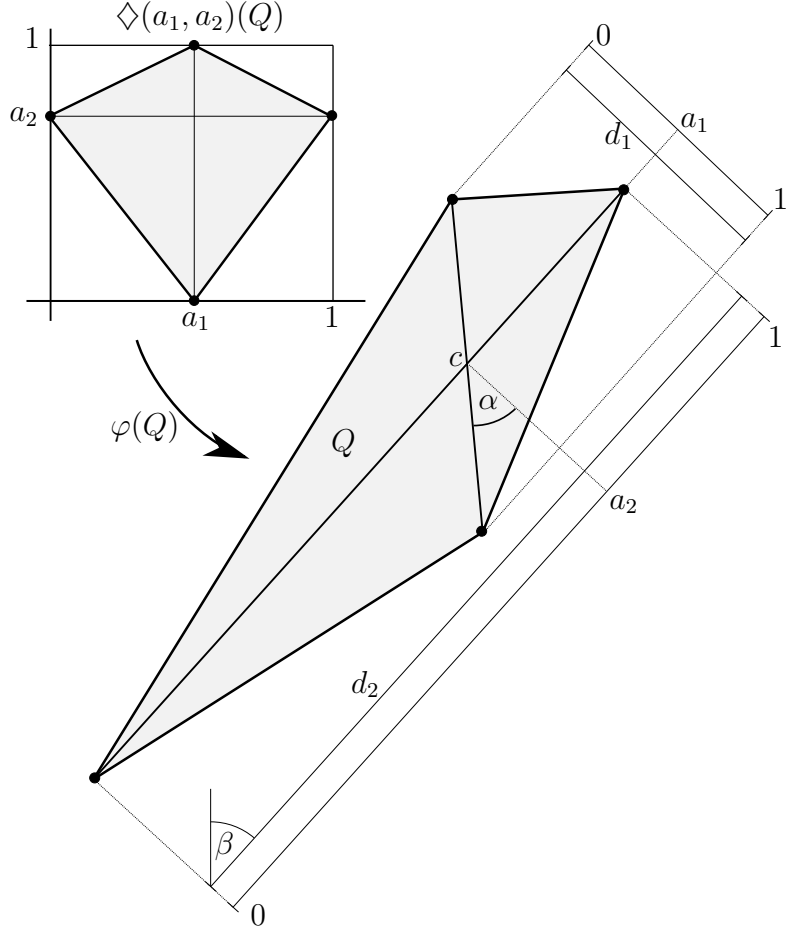
\caption[Illustration of $Q_{(\alpha,\beta,a_1,a_2,d_1,d_2,c_1,c_2)}$]{Illustration of $Q_{(\alpha,\beta,a_1,a_2,d_1,d_2,c_1c_2)}$. Here, $\alpha$ (representing the shear) and $\beta$ (representing the rotation) are chosen to have a negative sign. $d_1$ and $d_2$ represent the horizontal and vertical scaling (before shear and rotation). $c=(c_1,c_2)$ represents the translation.}
\label{img:quadrilateral}
\end{figure} 

Let us make this more precise: Let $Q$ be any convex quadrilateral in $\R^2$. Then, $Q$ can be described by the variables $a_1,a_2,c_1,c_2\in\R$, $d_1,d_2 > 0$, and $\alpha,\beta\in [0,2\phi]$ as indicated in Figure \ref{img:quadrilateral}. We write $Q=Q_{(\alpha,\beta,a_1,a_2,d_1,d_2,c_1,c_2)}$. Here, $c=(c_1,c_2)$ represents the translation. $a_1$ and $a_2$ describe the proportional relationship of the intersection of the two diagonals of $Q$. $d_1$ and $d_2$ represent the horizontal and vertical scaling, respectively, and $\alpha$ and $\beta$ the shear and the rotation, respectively. The order of these operations can be easily taken from Figure \ref{img:quadrilateral}.

Then we start by proving the following lemma:

\blem\label{Lem:QCaseLem1}
Let $Q_{(\alpha,\beta,a_1,a_2,d_1,d_2,c_1,c_2)}$ be any convex quadrilateral in $\R^2$. Then, one has
\beqq \mathring{Q}_{(\alpha,\beta,a_1,a_2,d_1,d_2,c_1,c_2)}\times A_{(\alpha,\beta,d_1,d_2)}\square(1)\stackrel{\text{sympl.}}{\cong} \diamondsuit(1,a_1,a_2)\times \square(1) \stackrel{\text{sympl.}}{\cong} B^4_{\sqrt{\frac{1}{\pi}}},\eeqq
where $A_{(\alpha,\beta,d_1,d_2)}$ is a linear transformation given by
\beq A_{(\alpha,\beta,d_1,d_2)}=R_{-\beta}\left(\widetilde{A}^{-1}_{(\alpha,d_1,d_2)}\right)^T\label{eq:QCaseLem1lineartransformation}\eeq
with
\beqq R_{-\beta}= \bpm \cos(-\beta) & -\sin (-\beta) \\ \sin(-\beta) & \cos(-\beta)\epm\; \text{ and }\; \widetilde{A}_{(\alpha,d_1,d_2)}=\bpm d_1 & 0 \\ d_1\tan(\alpha) & d_2 \epm.\eeqq
\elem


Clearly, this lemma proves Theorem \ref{Thm:answer2}(iii) (except for the remark concerning the case when $Q$ is a trapezoid--what we will deal with later) and what remained to be proven in Theorem \ref{Thm:answer3}. We also remark that $A_{(\alpha,\beta,d_1,d_2)}$ in \eqref{eq:QCaseLem1lineartransformation} specifies $A$ in \eqref{eq:answer2iii}.



\bpf[Proof of Lemma \ref{Lem:QCaseLem1}]
We notice that $R_{-\beta}$ and $\widetilde{A}_{(\alpha,\beta,d_1,d_2)}$ are defined in such a way that the above mentioned $\phi(Q)$ can be written as their composition--up to translations (which are defined in accordance with \eqref{eq:translationdefinitionst}):
\beqq \phi(Q)=s_{(-c_1,-c_2)}\circ L_{R_{-\beta}}\circ L_{\widetilde{A}_{(\alpha,\beta,d_1,d_2)}}\circ s_{(a_1,a_2)}.\eeqq
Here, the linear map $L_{R_{-\beta}}$--represented by the matrix $R_{-\beta}$--is responsible for the corresponding rotation by angle $-\beta$ and the linear map $L_{\widetilde{A}_{(\alpha,\beta,d_1,d_2)}}$--represented by the matrix $\widetilde{A}_{(\alpha,\beta,d_1,d_2)}$--is responsible for the corresponding scaling and shearing. So, we can write
\beqq \left(s_{(-c_1,-c_2)} \circ L_{R_{-\beta}}\circ L_{\widetilde{A}_{(\alpha,\beta,d_1,d_2)}}\circ s_{(a_1,a_2)}\right) (\diamondsuit(a_1,a_2))= Q_{(\alpha,\beta,a_1,a_2,d_1,d_2,c_1,c_2)}.\eeqq
This implies that
\begin{align*}
&\left(s_{(-c_1,-c_2)}\times \mathbb{1}_{\R^2}\right) \circ\left(\left(L_{R_{-\beta}}\circ L_{\widetilde{A}_{(\alpha,\beta,d_1,d_2)}}\right) \times \left(\left(L_{R_{-\beta}}\circ L_{\widetilde{A}_{(\alpha,\beta,d_1,d_2)}}\right)^T\right)^{-1}\right)\\
& \quad \quad \quad \quad \quad \quad \quad \quad \quad \quad \quad \quad\circ \left(s_{(a_1,a_2)}\times \mathbb{1}_{\R^2}\right)\\
= & \left(s_{(-c_1,-c_2)}\times \mathbb{1}_{\R^2}\right) \circ\left( L_{R_{-\beta} \widetilde{A}_{(\alpha,\beta,d_1,d_2)}} \times \left(L_{R_{-\beta} \widetilde{A}_{(\alpha,\beta,d_1,d_2)}}^T\right)^{-1}\right)\circ \left(s_{(a_1,a_2)}\times \mathbb{1}_{\R^2}\right)\\
=& \left(s_{(-c_1,-c_2)}\times \mathbb{1}_{\R^2}\right) \circ \left(L_{R_{-\beta} \widetilde{A}_{(\alpha,\beta,d_1,d_2)}} \times L_{\left(\left(R_{-\beta} \widetilde{A}_{(\alpha,\beta,d_1,d_2)}\right)^T\right)^{-1}}\right) \circ \left(s_{(a_1,a_2)}\times \mathbb{1}_{\R^2}\right)
\end{align*}
is a symplectomorphism producing
\beqq \diamondsuit(1,a_1,a_2)\times \square(1) \stackrel{\text{sympl.}}{\cong } \mathring{Q}_{(\alpha,\beta,a_1,a_2,d_1,d_2,c_1,c_2)} \times A_{(\alpha,\beta,d_1,d_2)}\square(1).\eeqq
Here, we used
\beqq \left(\left(R_{-\beta} \widetilde{A}_{(\alpha,\beta,d_1,d_2)}\right)^T\right)^{-1}=\left(R_{-\beta}^T\right)^{-1}\left(\widetilde{A}^T_{(\alpha,\beta,d_1,d_2)}\right)^{-1}=R_{-\beta}\left(\widetilde{A}^T_{(\alpha,\beta,d_1,d_2)}\right)^{-1}\eeqq
and \eqref{eq:QCaseLem1lineartransformation}. The rest of what is to be proven follows from Theorem \ref{Thm:Schlenk}.
\epf

\begin{figure}[h!]
\centering
\def\svgwidth{370pt}
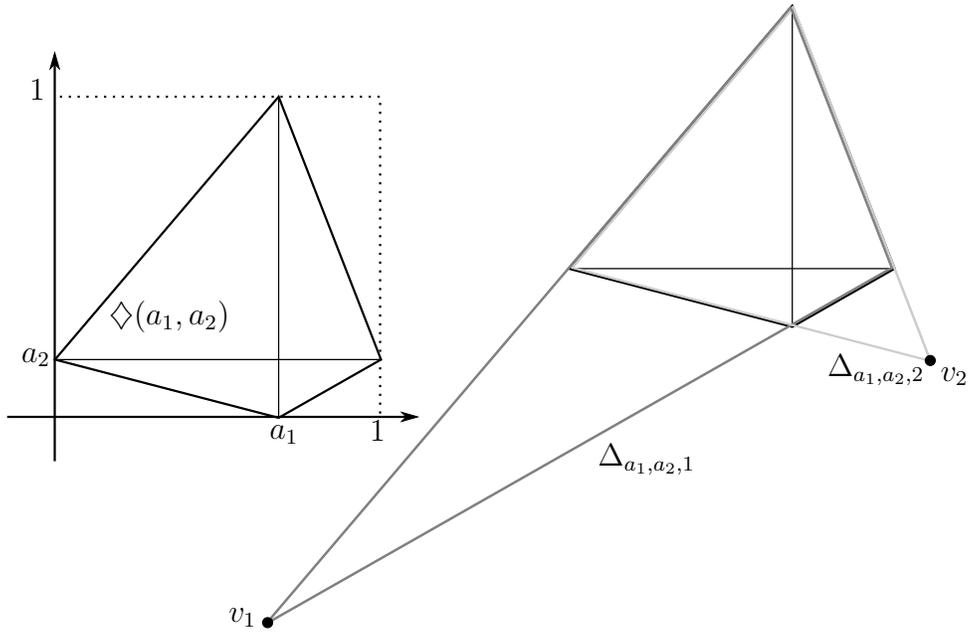
\caption[$\diamondsuit(a_1,a_2)$ is the intersection of the two uniquely determined triangles $\Delta_{a_1,a_2,1}$ and $\Delta_{a_1,a_2,2}$.]{$\diamondsuit(a_1,a_2)$ is the intersection of the two uniquely determined triangles $\Delta_{a_1,a_2,1}$ and $\Delta_{a_1,a_2,2}$.}
\label{img:diamondtriangleintersection}
\end{figure}

In order to show what remained to be proven in Theorem \ref{Thm:answer1} and for the convex-quadrilateral-case in general, we have to make some preparations: 

For any $\diamondsuit(a_1,a_2)$ with $a_1,a_2\in (0,1)$ we make the following observations:

If $a_1\notin\{a_2,1-a_2\}$, then there are two uniquely determined triangles $\Delta_{a_1,a_2,1}$ and $\Delta_{a_1,a_2,2}$ in $\R^2$ such that
\beqq \diamondsuit(a_1,a_2)=\Delta_{a_1,a_2,1}\cap \Delta_{a_1,a_2,2}\eeqq
(see Figure \ref{img:diamondtriangleintersection}). We let $v_1$ and $v_2$ be the vertices of $\Delta_{a_1,a_2,1}$ and $\Delta_{a_1,a_2,2}$ which are not contained in $\diamondsuit(a_1,a_2)$, respectively.
Then, we will prove that there are $\lambda_1,\lambda_2 >0$ such that $\pm\lambda_1 J\Delta_{a_1,a_2,1}$ and $\pm\lambda_2 J\Delta_{a_1,a_2,2}$ can be uniquely translated into $\square$ (the square with side length $1$ and centred at $(\frac{1}{2},\frac{1}{2})$) such that all their vertices are on $\partial \square$ and that $\pm\lambda_1 Jv_1$ and $\pm\lambda_2 J v_2$ are translated into vertices of $\square$, respectively.

\begin{figure}[h!]
\centering
\def\svgwidth{360pt}
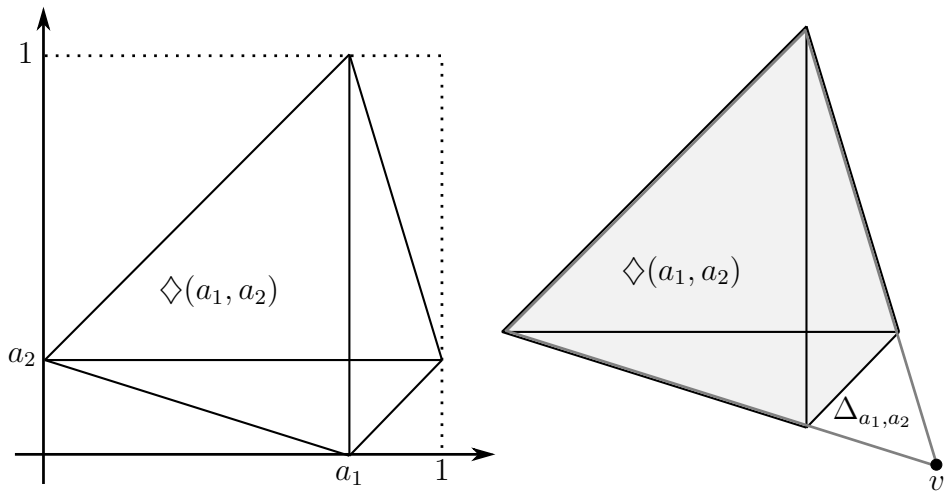
\caption[$\Delta_{a_1,a_2}$ is the unique volume-minimizing triangle-cover of $\diamondsuit(a_1,a_2)$.]{$\Delta_{a_1,a_2}$ is the unique volume-minimizing triangle-cover of $\diamondsuit(a_1,a_2)$.}
\label{img:diamondtrianglecover}
\end{figure}

If $a_1\in\{a_2,1-a_2\}$ with $(a_1,a_2)\neq (\frac{1}{2},\frac{1}{2})$, then there is one uniquely determined triangle--we call it $\Delta_{a_1,a_2}$--which is a volume-minimizing triangle-cover of $\diamondsuit(a_1,a_2)$ (see Figure \ref{img:diamondtrianglecover}). We let $v$ be the vertex of $\Delta_{a_1,a_2}$ which is not contained in $\diamondsuit(a_1,a_2)$.
Then, we will prove that there is a $\lambda >0$ such that $\pm \lambda J\Delta_{a_1,a_2}$ can be uniquely translated into $\square$ such that all their vertices are on $\partial\square$ and that $\pm\lambda J v$ are translated into a vertex of $\square$, repectively. Finally, we denote by $d$ the diagonal segment of $\square$ which is orthogonal to the two parallel sides of $\diamondsuit(a_1,a_2)$.

Then, we can prove the following proposition:

\bprop\label{Prop:quadrilateral12}
Let $Q_{(\alpha,\beta,a_1,a_2,d_1,d_2,c_1,c_2)}$ be any convex quadrilateral in $\R^2$.
\begin{itemize}
\item[(i)] If $a_1\notin\{a_2,1-a_2\}$, then Viterbo's conjecture is true for all Lagrangian products
\beq Q_{(\alpha,\beta,a_1,a_2,d_1,d_2,c_1,c_2)}\times T, \label{eq:quadrilateral1}\eeq
where $T$ can be any convex body in $\R^2$, if $\square$ is a volume-minimizing convex hull of translates of
\beq \pm\lambda_1 J\Delta_{a_1,a_2,1}\; \text{ and }\;\pm\lambda_2 J\Delta_{a_1,a_2,2}, \label{eq:translatesof1}\eeq
or, in other words, if $\square$ is in \beqq \big\{\conv\{\pm \lambda_1 J\Delta_{a_1,a_2,1}+\widetilde{t}_{\pm,1},\pm \lambda_2 J\Delta_{a_1,a_2,2}+\widetilde{t}_{\pm,2}\}:\widetilde{t}_{\pm,1},\widetilde{t}_{\pm,2} \text{ minimize... }\quad \quad\quad\quad\quad \quad\quad\quad\eeqq\beqq\quad \quad \quad \quad \quad ...\vol\left(\conv\{\pm \lambda_1 J\Delta_{a_1,a_2,1}+t_{\pm,1},\pm \lambda_2 J\Delta_{a_1,a_2,2}+t_{\pm,2}\}\right)\text{ over all }t_{\pm,1},t_{\pm,2}\in\R^2\big\}.\eeqq
\item[(ii)] If $a_1\in\{a_2,1-a_2\}$ with $(a_1,a_2)\neq (\frac{1}{2},\frac{1}{2})$, then Viterbo's conjecture is true for all Lagrangian products
\beqq Q_{(\alpha,\beta,a_1,a_2,d_1,d_2,c_1,c_2)}\times T, \eeqq
where $T$ can be any convex body in $\R^2$, if $\square$ is a volume-minimizing convex hull of translates of
\beq \pm\lambda J\Delta_{a_1,a_2}\;\text{ and }\; d.\label{eq:translatesof2}\eeq
\end{itemize}
\eprop

We remark that the set of equality cases of Viterbo's conjecture for the configuration
\beqq Q_{(\alpha,\beta,a_1,a_2,d_1,d_2,c_1,c_2)}\times T\eeqq
is determined by the set of convex bodies which are volume-minimizing convex hulls of translates of the convex sets in \eqref{eq:translatesof1} and \eqref{eq:translatesof2}, respectively.

We start proving Proposition \ref{Prop:quadrilateral12} by showing the following proposition:

\bprop\label{Prop:quadrilateral122}
Let $Q_{(\alpha,\beta,a_1,a_2,d_1,d_2,c_1,c_2)}$ be any covex quadrilateral in $\R^2$.
\begin{itemize}
\item[(i)] If $a_1\notin\{a_2,1-a_2\}$, then we find $\lambda_1,\lambda_2>0$ such that $\pm\lambda_1 J\Delta_{a_1,a_2,1}$ and $\pm\lambda_2 J\Delta_{a_1,a_2,2}$ can be uniquely translated into $\square$ such that all their vertices are boundary points of $\square$ and that $\pm\lambda_1 Jv_1$ and $\pm\lambda_2 J v_2$ are translated into vertices of $\square$, respectively.
\item[(ii)] If $a_1\in\{a_2,1-a_2\}$ with $(a_1,a_2)\neq(\frac{1}{2},\frac{1}{2})$, then we find $\lambda >0$ such that $\pm \lambda J\Delta_{a_1,a_2}$ can be uniquely translated into $\square$ such that all their vertices are boundary points of $\square$ and that $\pm\lambda J v$ are translated into vertices of $\square$, repectively.
\end{itemize}
\eprop

\bpf
\underline{Ad(i)}: As will be seen in the further course of the proof, it is sufficient to prove that one can find a $\lambda_{+,2}$ such that $\lambda_{+,2}J\Delta_{a_1,a_2,2}$ can be uniquely translated into $\square$ such that all its vertices are boundary points of $\square$ with the property that $\lambda_{+,2}Jv_2$ is translated into a vertex of $\square$.

\begin{figure}[h!]
\centering
\def\svgwidth{370pt}
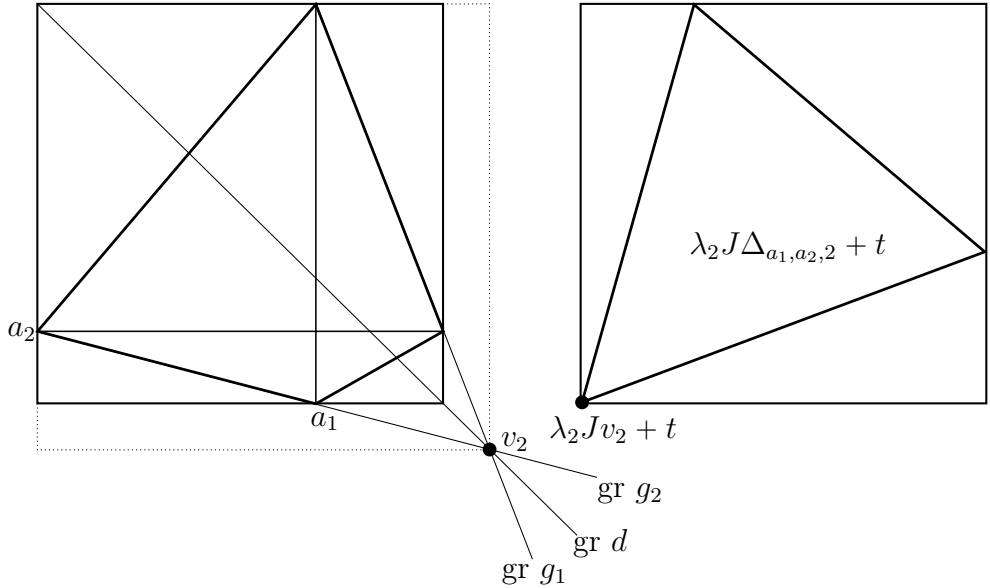
\caption[$\Delta_{a_1,a_2,2}$'s vertex $v_2$ lies on one of the extended diagonals of $\square$.]{$\Delta_{a_1,a_2,2}$'s vertex $v_2$ lies on one of the extended diagonals of $\square$. By rotating $\Delta_{a_1,a_2,2}$ by $J$ and scaling it by the factor $\lambda_2$, the image can be translated (by $t\in\R^2$) such that all its vertices lie on $\partial\square$ while satisfying that $\lambda_2Jv_2 + t$ is a vertex of $\square$.}
\label{img:vertexextendeddiagonal}
\end{figure}

For this, we make the following observation: By construction of $\Delta_{a_1,a_2,2}$, its vertex $v_2$ lies on one of the extended diagonals of $\square$ (see Figure \ref{img:vertexextendeddiagonal}).

Let us prove this rigorously: We extend the edges of $\Delta_{a_1,a_2,2}$ through $v_2$ and understand these lines as graphs of affine functions $g_1$ and $g_2$ which depend on the horizontal coordinate $x$, respectively. We now show that $v_2$--as intersection point of $g_1$ and $g_2$--lies on one of the extended diagonals of the square which we will also understand as an affine function depending on $x$, denoted by $d$. A short calculation yields
\beqq g_1(x)=-\frac{1-a_2}{1-a_1}x+\frac{1-a_1a_2}{1-a_1}\quad \text{ and }\quad g_2(x)=-\frac{a_2}{a_1}x+a_2.\eeqq
As consequence, one derives
\beqq v_2=\left(\frac{a_1(1-a_2)}{a_1-a_2},\frac{a_2(a_1-1)}{a_1-a_2}\right).\eeqq
A trivial calculation yields that $v_2$ lies in fact on the graph of
\beqq d(x)=1-x.\eeqq

Therefore, the vertices of $\Delta_{a_1,a_2,2}$ lie on an extended square. By the above calculations, the side length of the extended square is $\frac{a_1(1-a_2)}{a_1-a_2}$. By rotating $\Delta_{a_1,a_2,2}$ by $J$ and scaling it by $\lambda_{+,2}=\frac{a_1-a_2}{a_1(1-a_2)}$, $\lambda_{+,2}J\Delta_{a_1,a_2,2}$ clearly can be uniquely translated into $\square$ such that all its vertices are on $\partial\square$ while satisfying that $\lambda_{+,2}Jv_2$ is a vertex of $\square$.

That the same argument applies for the other instances is clear from the fact that $\square$ and $\diamondsuit(a_1,a_2)$ can be rotated by $\frac{\pi}{2}$-angle-steps in order to produce situations as rigorously discussed above. 

\underline{Ad(ii)}: The proof runs analogously to the one for (i).
\epf

We remark that from the proof of Proposition \ref{Prop:quadrilateral122} it follows that
\beq \lambda_1=\frac{1-a_2-a_1}{(1-a_2)(1-a_1)} \quad \text{ and } \quad \lambda_2= \frac{a_1-a_2}{a_1(1-a_2)}\label{eq:lambdas}\eeq
(The assignment of these factors is based on the definitions of $\Delta_{a_1,a_2,1}$ and $\Delta_{a_1,a_2,2}$ in Figure \ref{img:diamondtriangleintersection}. Additionally, one has to note that these factors won't be negative due to the condition on the relation between $a_1$ and $a_2$ which stems from the fact that the choice of $\Delta_{a_1,a_2,1}$ and $\Delta_{a_1,a_2,2}$ is not arbitrary.). The value of $\lambda$ results accordingly.

Then, we come to the proof of Propositon \ref{Prop:quadrilateral12}:

\bpf[Proof of Proposition \ref{Prop:quadrilateral12}]
\underline{Ad (i)}: We first note that it is sufficient to prove the theorem's statement for \eqref{eq:quadrilateral1} replaced by
\beqq \diamondsuit(a_1,a_2)\times T,\eeqq
where $T$ can be any convex body in $\R^2$. This is due to the facts that one can find an affine transformation $\phi=\phi(Q)$ mapping $Q_{(\alpha,\beta,a_1,a_2,d_1,d_2,c_1,c_2)}$ onto $\diamondsuit(a_1,a_2)$ and that $\phi\times \left(\phi^T\right)^{-1}$ is a symplectomorphism--under which Viterbo's conjecture is invariant.

We already know from Theorem \ref{Thm:answer2}(ii) that $\square \times \diamondsuit(a_1,a_2)$ (and therefore also $\diamondsuit(a_1,a_2)\times\square$) is an equality case of Viterbo's conjecture. The idea now is the following: We reexamine the equality case $\square\times \diamondsuit(a_1,a_2)$ by using the Minkowski billiard characterization of its EHZ-capacity in order to show that $\square$ is the convex hull of some of the $\ell_{\diamondsuit(a_1,a_2)}$-minimizing closed $(\square,\diamondsuit(a_1,a_2))$-Minkowski billiard trajectories (whose orbits will turn out to be translates of the boundaries of $\pm\lambda_1 J\Delta_{a_1,a_2,1}$ and $\pm\lambda_2 J\Delta_{a_1,a_2,2}$, respectively). Under the condition that $\square$ moreover is a volume-minimizing convex cover of the translates of the aforementioned $\ell_{\diamondsuit(a_1,a_2)}$-minimizing closed $(\square,\diamondsuit(a_1,a_2))$-Minkowski billiard trajectories in general, this would imply that for any other convex body $T$ in $\R^2$ with $\vol(T)=\vol(\square)$ the aforementioned $\ell_{\diamondsuit(a_1,a_2)}$-minimizing closed $(\square,\diamondsuit(a_1,a_2))$-Minkowski billiard trajectories cannot all be translated into the interior of $T$ (otherwise, a scaled copy of $T$ would be a counterexample to $\square$ being a volume-minimizing convex cover of the translates of the minimizing Minkowski billiard trajectories). Referring to Theorem \ref{Thm:Chap7relationship}, this would imply that the EHZ-capacity of $T\times \diamondsuit(a_1,a_2)$ is less or equal than the EHZ-capacity of $\square \times \diamondsuit(a_1,a_2)$ (since the EHZ-capacity of $T\times\diamondsuit(a_1,a_2)$ equals the minimal $\ell_{\diamondsuit(a_1,a_2)}$-length of the closed polygonal curves which cannot be translated into the interior of $T$). This finally would imply that the EHZ-capacity of $\square\times \diamondsuit(a_1,a_2)$ is greater or equal than the EHZ-capacity of any configuration $\diamondsuit(a_1,a_2)\times T$, where $T$ is any convex body in $\R^2$ with $\vol(T)=\vol(\square)$. This finally would imply that Viterbo's conjecture is true for any configuration $\diamondsuit(a_1,a_2)\times T$, where $T$ is any convex body in $\R^2$.

\begin{figure}[h!]
\centering
\def\svgwidth{400pt}
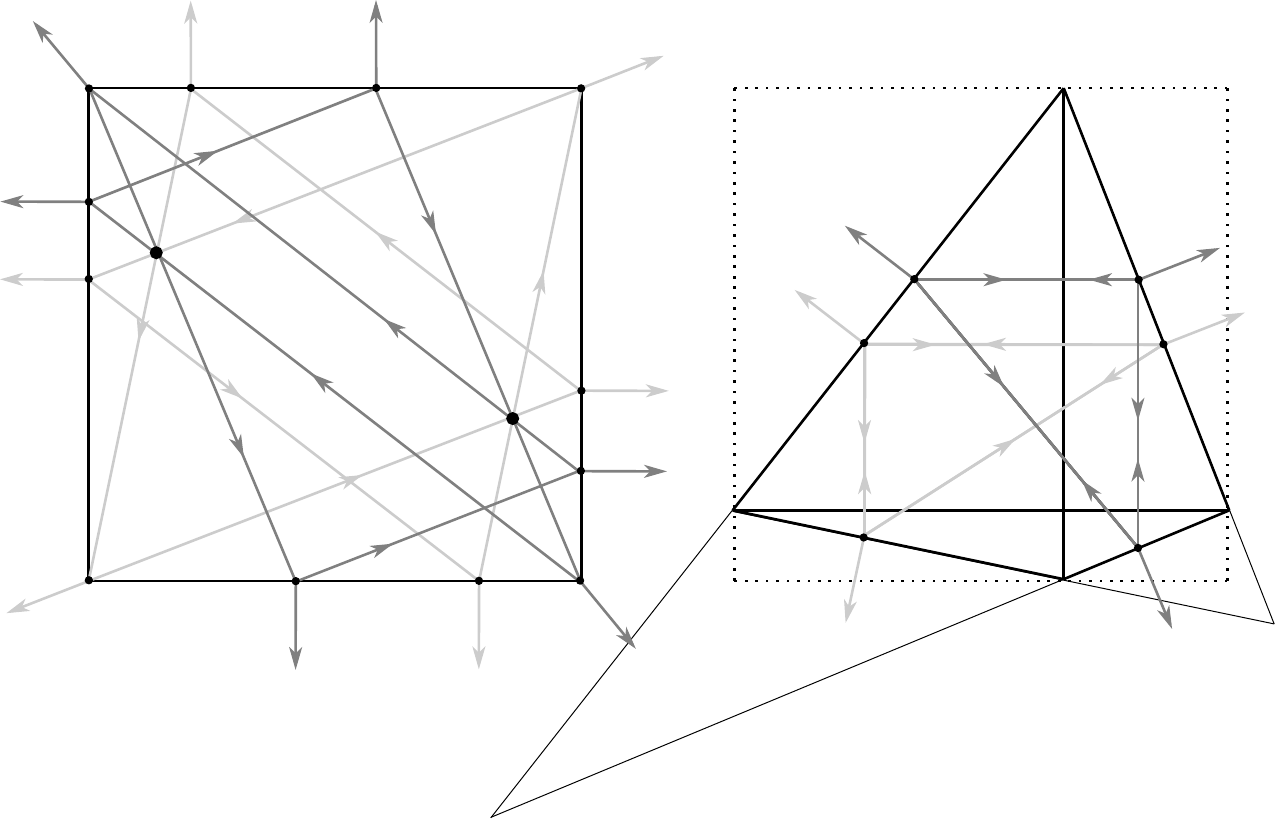
\caption[The four $\ell_{\diamondsuit(a_1,a_2)}$-minimizing closed $(\square,\diamondsuit(a_1,a_2))$-Minkowski billiard trajectories whose convex hull is $\square$.]{Illustration of the four $\ell_{\diamondsuit(a_1,a_2)}$-minimizing closed $(\square,\diamondsuit(a_1,a_2))$-Minkowski billiard trajectories whose convex hull is $\square$. An interesting sidenote: there are two symmetry-points in which all four trajectories intersect.}
\label{img:fourMbt}
\end{figure}

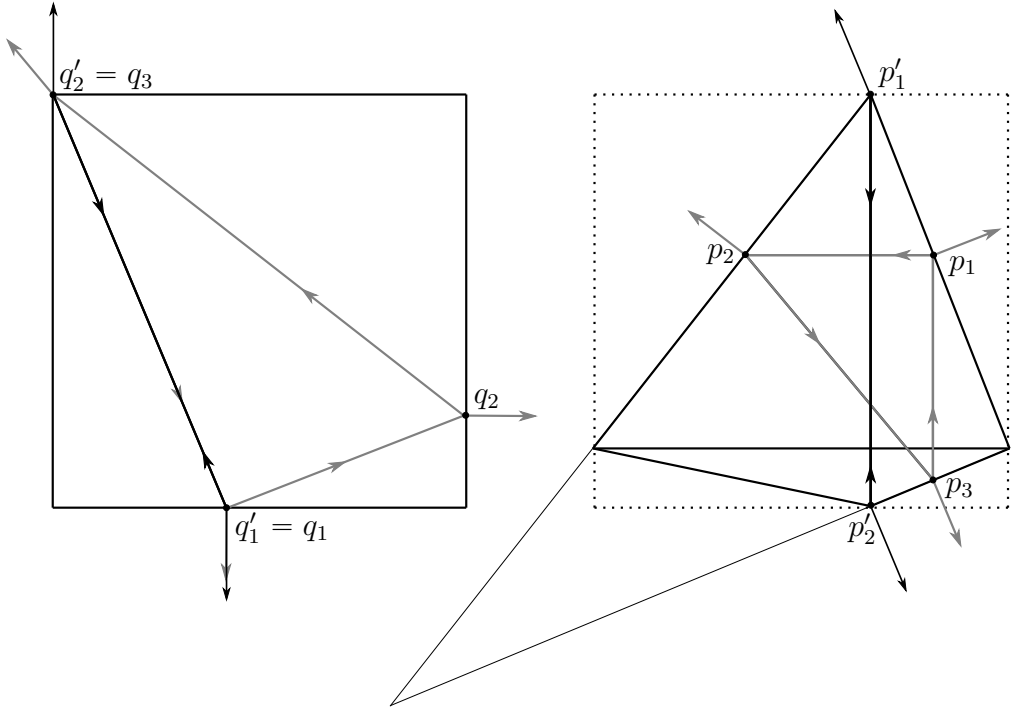
\begin{figure}[h!]
\centering
\def\svgwidth{380pt}
\begingroup%
  \makeatletter%
  \providecommand\color[2][]{%
    \errmessage{(Inkscape) Color is used for the text in Inkscape, but the package 'color.sty' is not loaded}%
    \renewcommand\color[2][]{}%
  }%
  \providecommand\transparent[1]{%
    \errmessage{(Inkscape) Transparency is used (non-zero) for the text in Inkscape, but the package 'transparent.sty' is not loaded}%
    \renewcommand\transparent[1]{}%
  }%
  \providecommand\rotatebox[2]{#2}%
  \newcommand*\fsize{\dimexpr\f@size pt\relax}%
  \newcommand*\lineheight[1]{\fontsize{\fsize}{#1\fsize}\selectfont}%
  \ifx\svgwidth\undefined%
    \setlength{\unitlength}{344.90784717bp}%
    \ifx\svgscale\undefined%
      \relax%
    \else%
      \setlength{\unitlength}{\unitlength * \real{\svgscale}}%
    \fi%
  \else%
    \setlength{\unitlength}{\svgwidth}%
  \fi%
  \global\let\svgwidth\undefined%
  \global\let\svgscale\undefined%
  \makeatother%
  \begin{picture}(1,0.70209445)%
    \lineheight{1}%
    \setlength\tabcolsep{0pt}%
    \put(0,0){\includegraphics[width=\unitlength,page=1]{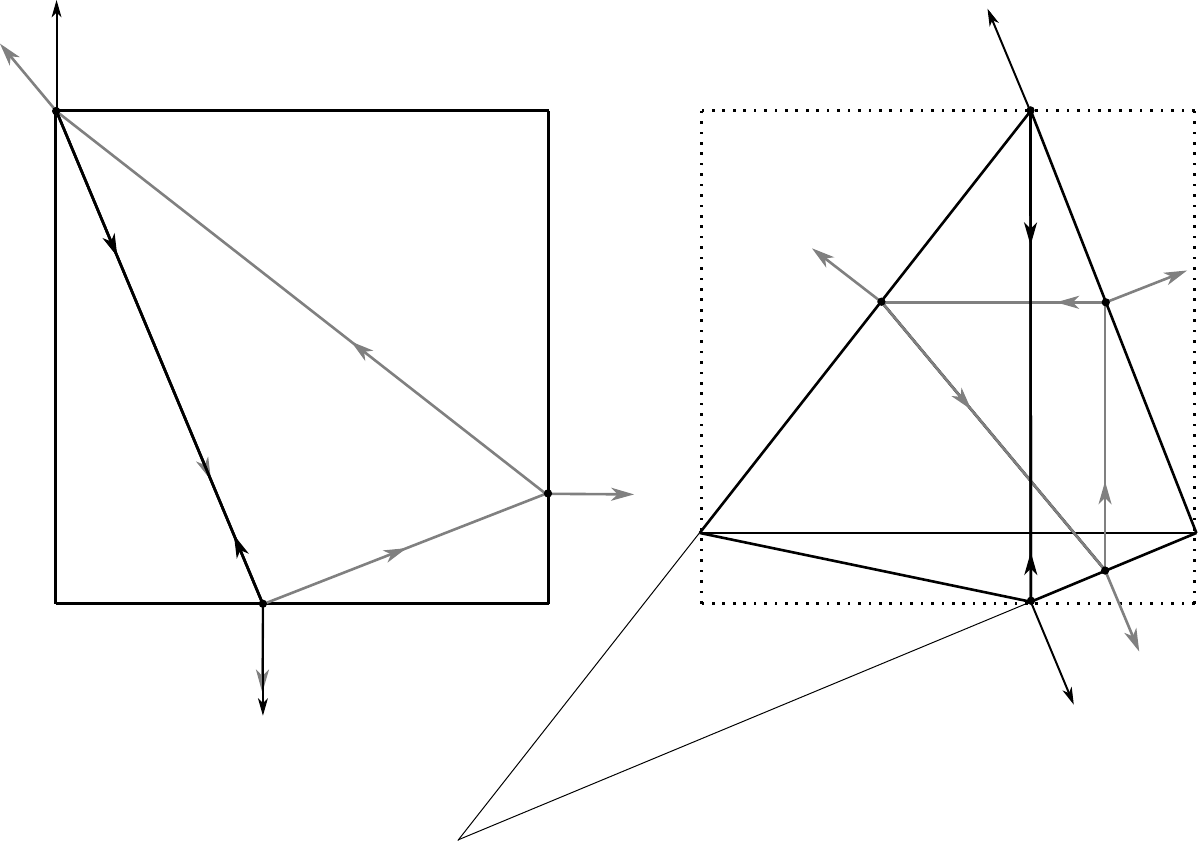}}%
    \put(0.22788207,0.17077201){\color[rgb]{0,0,0}\makebox(0,0)[lt]{\lineheight{1.25}\smash{\begin{tabular}[t]{l}$q_1'=q_1$\end{tabular}}}}%
    \put(0.05293899,0.62121255){\color[rgb]{0,0,0}\makebox(0,0)[lt]{\lineheight{1.25}\smash{\begin{tabular}[t]{l}$q_2'=q_3$\end{tabular}}}}%
    \put(0.46574225,0.30200696){\color[rgb]{0,0,0}\makebox(0,0)[lt]{\lineheight{1.25}\smash{\begin{tabular}[t]{l}$q_2$\end{tabular}}}}%
    \put(0.86933558,0.62283338){\color[rgb]{0,0,0}\makebox(0,0)[lt]{\lineheight{1.25}\smash{\begin{tabular}[t]{l}$p_1'$\end{tabular}}}}%
    \put(0.83865945,0.16968443){\color[rgb]{0,0,0}\makebox(0,0)[lt]{\lineheight{1.25}\smash{\begin{tabular}[t]{l}$p_2'$\end{tabular}}}}%
    \put(0.93573529,0.21317429){\color[rgb]{0,0,0}\makebox(0,0)[lt]{\lineheight{1.25}\smash{\begin{tabular}[t]{l}$p_3$\end{tabular}}}}%
    \put(0.9388417,0.43528329){\color[rgb]{0,0,0}\makebox(0,0)[lt]{\lineheight{1.25}\smash{\begin{tabular}[t]{l}$p_1$\end{tabular}}}}%
    \put(0.69809419,0.44654402){\color[rgb]{0,0,0}\makebox(0,0)[lt]{\lineheight{1.25}\smash{\begin{tabular}[t]{l}$p_2$\end{tabular}}}}%
  \end{picture}%
\endgroup%

\caption[The closed $(\square,\diamondsuit(a_1,a_2))$-Minkowski billiard trajectory $q=(q_1,q_2,q_3)$ is $\ell_{\diamondsuit(a_1,a_2)}$-minimizing with $\ell_{\diamondsuit(a_1,a_2)}(q)=\ell_{\diamondsuit(a_1,a_2)}(q')=1$.]{The closed $(\square,\diamondsuit(a_1,a_2))$-Minkowski billiard trajectory $q=(q_1,q_2,q_3)$ is $\ell_{\diamondsuit(a_1,a_2)}$-minimizing with $\ell_{\diamondsuit(a_1,a_2)}(q)=\ell_{\diamondsuit(a_1,a_2)}(q')=1$.}
\label{img:oneMbt}
\end{figure}

So, it remains to show that there are $\ell_{\diamondsuit(a_1,a_2)}$-minimizing closed $(\square,\diamondsuit(a_1,a_2))$-Minkowski billiard trajectories whose convex hull is $\square$ and that the orbits of these trajectories can be represented by translates of the boundaries of $\pm\lambda_1 J\Delta_{a_1,a_2,1}$ and $\pm\lambda_2 J\Delta_{a_1,a_2,2}$: From the proof of Lemma \ref{Lem:squarestep1} we know that
\beq c_{EHZ}(\square\times \diamondsuit(a_1,a_2))=1,\label{eq:kp}\eeq
where the horizontal and vertical diameter of $\square$ represents the orbits of $\ell_{\diamondsuit(a_1,a_2)}$-minimizing closed $(\square,\diamondsuit(a_1,a_2))$-Minkowski billiard trajectories which have two bouncing points. However, these are not the only ones.  Let us consider the closed $(\square,\diamondsuit(a_1,a_2))$-Minkowski billiard trajectories which have three bouncing points. By considering Proposition \ref{Prop:quadrilateral122}(i) and the algorithm of how to construct closed Minkowski billiard trajectories manually (see \cite{KruppRudolf2022}), we can affirm that one gets the four $(\square,\diamondsuit(a_1,a_2))$-Minkowski billiard trajectories which are indicated in Figure \ref{img:fourMbt}.
Due to the Minkowski billiard reflection rule, their orbits can be represented by translates of the boundaries of $\pm\lambda_1 J\Delta_{a_1,a_2,1}$ and $\pm\lambda_2 J\Delta_{a_1,a_2,2}$. One clearly sees that their convex hull is $\square$. It remains to prove their $\ell_{\diamondsuit(a_1,a_2)}$-minimality: We will do this exemplary with the trajectory $q=(q_1,q_2,q_3)$--indicated in Figure \ref{img:oneMbt}; for the remaining three it is the same.
The $\ell_{\diamondsuit(a_1,a_2)}$-length of $q$ is
\beqq \ell_{\diamondsuit(a_1,a_2)}(q)=\mu_{\diamondsuit(a_1,a_2)^\circ}(q_2-q_1)+\mu_{\diamondsuit(a_1,a_2)^\circ}(q_3-q_2)+\mu_{\diamondsuit(a_1,a_2)^\circ}(q_1-q_3).\eeqq
Because of
\beqq \mu_{\diamondsuit(a_1,a_2)^\circ}(q_1-q_3)=\mu_{\diamondsuit(a_1,a_2)^\circ}(q_1'-q_2')\eeqq
and
\beq \mu_{\diamondsuit(a_1,a_2)^\circ}(q_2-q_1)+\mu_{\diamondsuit(a_1,a_2)^\circ}(q_3-q_2)=\mu_{\diamondsuit(a_1,a_2)^\circ}(q_3-q_1)=\mu_{\diamondsuit(a_1,a_2)^\circ}(q_2'-q_1'),\label{eq:fe}\eeq
where the first equality in \eqref{eq:fe} follows from a property of the Minkowski functional shown in \cite[Lemma 5.4]{Rudolf2022b} together with the fact that the normal vectors in $p_1$ and $p_2$ represent the rays between which the normal cone in $p_1'$ is enclosed. Therefore, one has
\beqq \ell_{\diamondsuit(a_1,a_2)}(q)=\mu_{\diamondsuit(a_1,a_2)^\circ}(q_2'-q_1')+\mu_{\diamondsuit(a_1,a_2)^\circ}(q_1'-q_2')=\ell_{\diamondsuit(a_1,a_2)}(q')=\langle q_2'-q_1',p_1'-p_2'\rangle =1\eeqq
which together with \eqref{eq:kp} implies that $q$ is an $\ell_{\diamondsuit(a_1,a_2)}$-minimizing closed $(\square,\diamondsuit(a_1,a_2))$-Minkowski billiard trajectory.

\begin{figure}[h!]
\centering
\def\svgwidth{390pt}
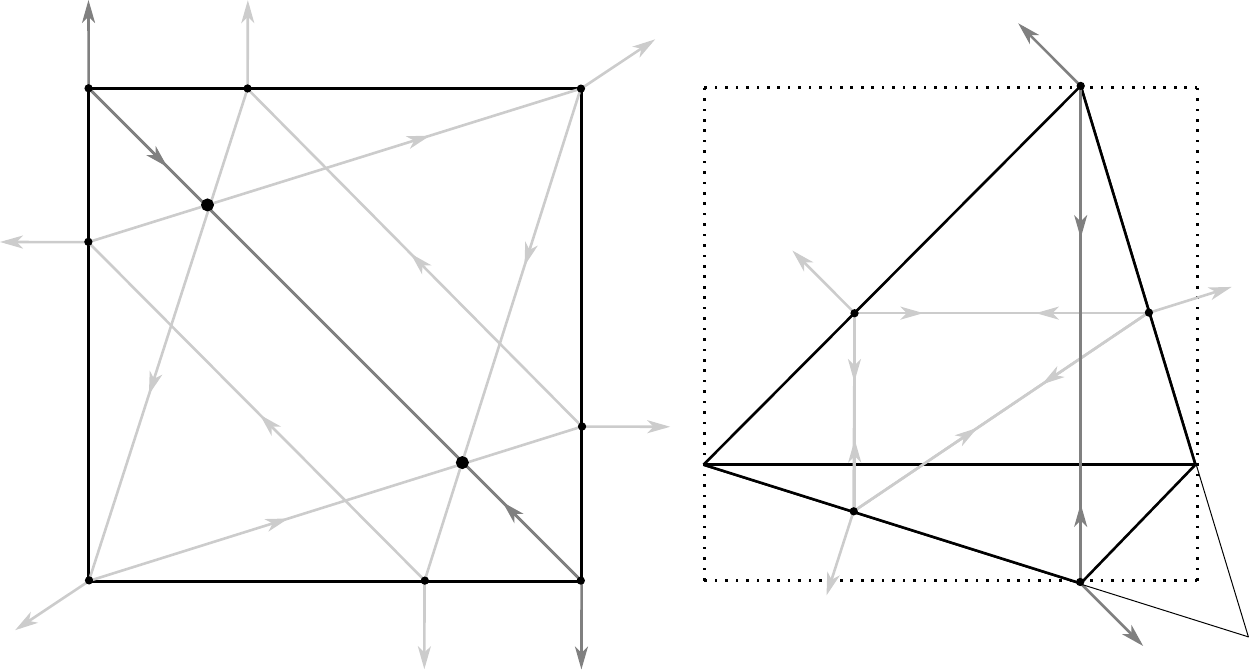
\caption[The three $\ell_{\diamondsuit(a_1,a_2)}$-minimizing closed $(\square,\diamondsuit(a_1,a_2))$-Minkowski billiard trajectories whose convex hull is $\square$.]{Illustration of the three $\ell_{\diamondsuit(a_1,a_2)}$-minimizing closed $(\square,\diamondsuit(a_1,a_2))$-Minkowski billiard trajectories whose convex hull is $\square$. An interesting sidenote: there are two symmetry-points in which all three trajectories intersect.}
\label{img:threeMbt}
\end{figure}

\underline{Ad (ii)}: The argumentation runs exactly like in (i)--with Proposition \ref{Prop:quadrilateral122}(i) replaced by Proposition \ref{Prop:quadrilateral122}(ii). The only difference in this case is, when reexamine the equality case $\square\times\diamondsuit(a_1,a_2)$ via the Minkowski billiard characterization of the EHZ-capacity, one notice that $\square$ is the convex hull of translates of $\pm\lambda J\Delta_{a_1,a_2}$ and $d$ (see Figure \ref{img:threeMbt}). Their $\ell_{\diamondsuit(a_1,a_2)}$-minimality is guaranteed by the arguments in the proofs of Lemma \ref{Lem:squarestep2} and Proposition \ref{Prop:quadrilateral122}(i).
\epf

\begin{figure}[h!]
\centering
\def\svgwidth{280pt}
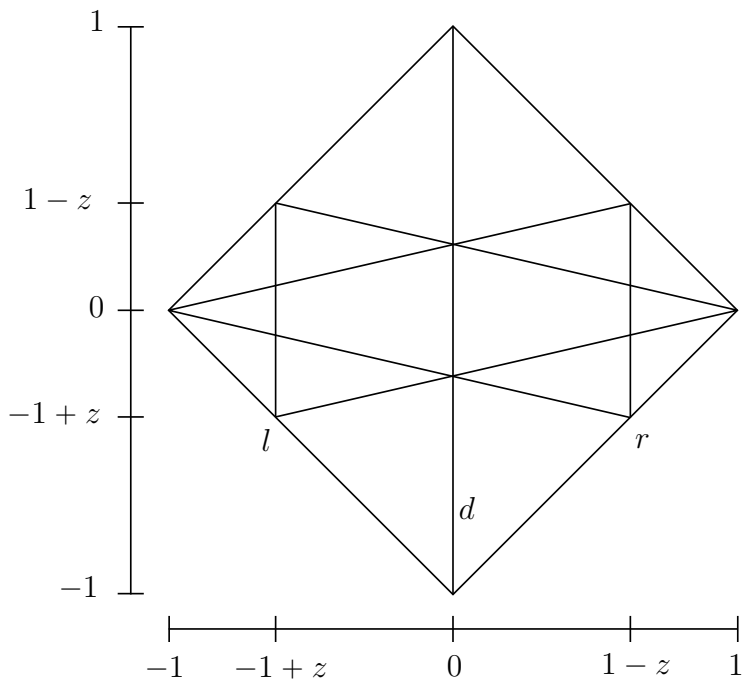
\caption{The original definition of the two triangles and the line segment depending on $z\in(0,1)$. The three apparent parallel vertical lines are denoted by $l$, $d$, and $r$, respectively--from left to right.}
\label{img:PackingProblem1}
\end{figure}

In order to show what remained to be proven in Theorems \ref{Thm:answer1} and \ref{Thm:answer2}(iii), we will prove the condition of Proposition \ref{Prop:quadrilateral12}(ii). For that, it suffices to prove the following proposition:

\bprop
Let $z\in(0,1)$. We consider the two triangles given by the vertices
\beqq (-1+z,z),\; (-1+z,-z),\; (1,0)\quad \text{ and }\quad (1-z,z),\; (1-z,-z),\; (-1,0)\eeqq
and the line segment given by the vertices $(0,-1)$ and $(0,1)$ (see Figure \ref{img:PackingProblem1}). Then, the square with vertices
\beqq (0,-1),\; (1,0),\; (0,1),\; (-1,0)\eeqq
--which is the convex hull of the three aforementioned sets--is the unique (up to translation) volume-minimizing convex hull of translates of the two triangles and the line segment.
\eprop

\bpf
We prove this statement by a case distinction based on the order of the three apparent parallel vertical lines. We have the following six cases:
\beqq dlr,\; lrd,\; ldr,\; rld,\; drl,\; rdl.\eeqq
Here, by $l$, $d$, and $r$, we denote the left, middle, and right vertical lines based on the original definition of the two triangles and the line segment, respectively. So, for instance, $dlr$ stands for the arrangement of translates of the three sets having the vertical line corresponding to the line segment on the left and the vertical lines corresponding to the two triangles in their original order following on the right.

Due to the symmetry of the convex hulls of translates of the three sets--due to the fact that the two triangles are rotated copies of each other--, the cases $dlr$ and $lrd$ as well as the cases $rld$ and $drl$ are equivalent. So, it is enough to focus on the four cases
\beqq dlr,\; ldr,\; drl,\; rdl.\eeqq


\begin{figure}[h!]
\centering
\def\svgwidth{380pt}
\begingroup%
  \makeatletter%
  \providecommand\color[2][]{%
    \errmessage{(Inkscape) Color is used for the text in Inkscape, but the package 'color.sty' is not loaded}%
    \renewcommand\color[2][]{}%
  }%
  \providecommand\transparent[1]{%
    \errmessage{(Inkscape) Transparency is used (non-zero) for the text in Inkscape, but the package 'transparent.sty' is not loaded}%
    \renewcommand\transparent[1]{}%
  }%
  \providecommand\rotatebox[2]{#2}%
  \newcommand*\fsize{\dimexpr\f@size pt\relax}%
  \newcommand*\lineheight[1]{\fontsize{\fsize}{#1\fsize}\selectfont}%
  \ifx\svgwidth\undefined%
    \setlength{\unitlength}{430.32073832bp}%
    \ifx\svgscale\undefined%
      \relax%
    \else%
      \setlength{\unitlength}{\unitlength * \real{\svgscale}}%
    \fi%
  \else%
    \setlength{\unitlength}{\svgwidth}%
  \fi%
  \global\let\svgwidth\undefined%
  \global\let\svgscale\undefined%
  \makeatother%
  \begin{picture}(1,0.35849035)%
    \lineheight{1}%
    \setlength\tabcolsep{0pt}%
    \put(0,0){\includegraphics[width=\unitlength,page=1]{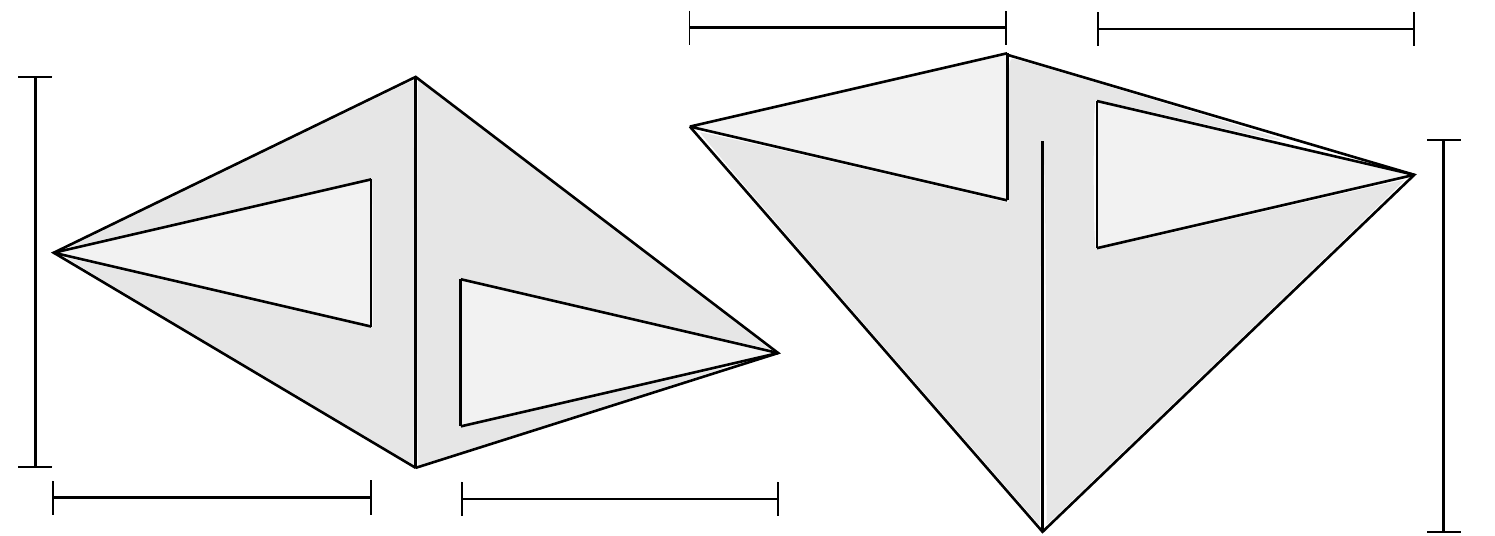}}%
    \put(0.11517485,0.00308923){\color[rgb]{0,0,0}\makebox(0,0)[lt]{\lineheight{1.25}\smash{\begin{tabular}[t]{l}$2-z$\end{tabular}}}}%
    \put(0.4185731,0.00170771){\color[rgb]{0,0,0}\makebox(0,0)[lt]{\lineheight{1.25}\smash{\begin{tabular}[t]{l}$2-z$\end{tabular}}}}%
    \put(-0.00096449,0.17226156){\color[rgb]{0,0,0}\makebox(0,0)[lt]{\lineheight{1.25}\smash{\begin{tabular}[t]{l}$2$\end{tabular}}}}%
    \put(0.55053457,0.34966245){\color[rgb]{0,0,0}\makebox(0,0)[lt]{\lineheight{1.25}\smash{\begin{tabular}[t]{l}$2-z$\end{tabular}}}}%
    \put(0.83873355,0.34779508){\color[rgb]{0,0,0}\makebox(0,0)[lt]{\lineheight{1.25}\smash{\begin{tabular}[t]{l}$2-z$\end{tabular}}}}%
    \put(0.97878694,0.11126054){\color[rgb]{0,0,0}\makebox(0,0)[lt]{\lineheight{1.25}\smash{\begin{tabular}[t]{l}$2$\end{tabular}}}}%
  \end{picture}%
\endgroup%

\caption{Case $rdl$.}
\label{img:PackingProblem2}
\end{figure}
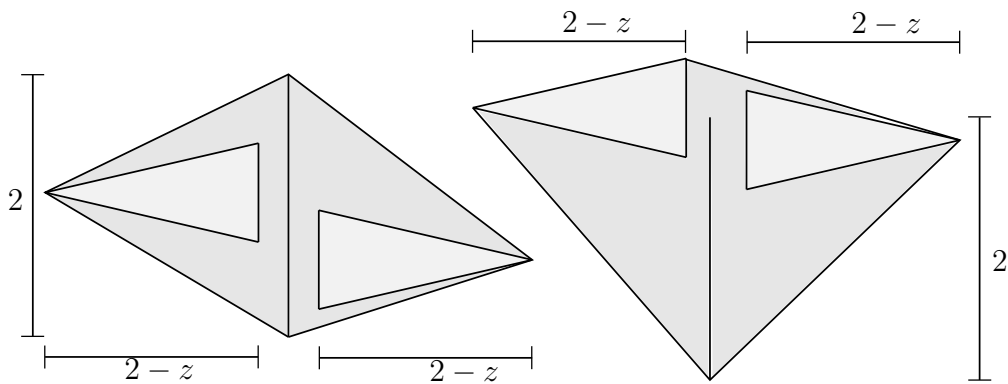

\underline{Case $rdl$}: We have a situation as shown in Figure \ref{img:PackingProblem2}. The volume of the convex hull of any arrangement of the two triangles and the line segment is greater or equal than
\beqq 2\frac{2(2-z)}{2}=4-2z\eeqq
which for $z\in(0,1)$ is greater than the volume of the square (which is $2$).

\begin{figure}[h!]
\centering
\def\svgwidth{340pt}
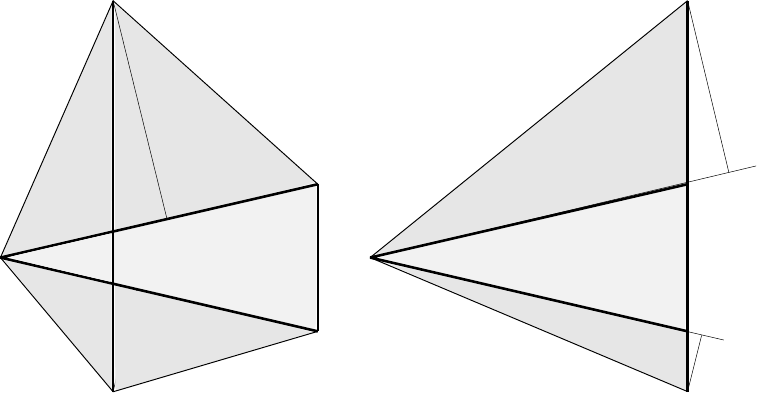
\caption{Case $drl$.}
\label{img:PackingProblem3}
\end{figure}

\underline{Case $drl$}: That the theorem's statement is true for the $drl$-case can be traced back to the fact that it is true for the $rdl$-case: As one can see in Figure \ref{img:PackingProblem3}, the area (additional to the convex hull of the two triangles) which has to be considered for any arrangement within the $drl$-case (marked dark grey) is greater or equal than the area (additional to the convex hull of the two triangles) which has to be considered for a certain arrangement within the $rdl$-case.

\underline{Case $dlr$}: We divide this case into two further cases:

\begin{figure}[h!]
\centering
\def\svgwidth{380pt}
\begingroup%
  \makeatletter%
  \providecommand\color[2][]{%
    \errmessage{(Inkscape) Color is used for the text in Inkscape, but the package 'color.sty' is not loaded}%
    \renewcommand\color[2][]{}%
  }%
  \providecommand\transparent[1]{%
    \errmessage{(Inkscape) Transparency is used (non-zero) for the text in Inkscape, but the package 'transparent.sty' is not loaded}%
    \renewcommand\transparent[1]{}%
  }%
  \providecommand\rotatebox[2]{#2}%
  \newcommand*\fsize{\dimexpr\f@size pt\relax}%
  \newcommand*\lineheight[1]{\fontsize{\fsize}{#1\fsize}\selectfont}%
  \ifx\svgwidth\undefined%
    \setlength{\unitlength}{287.24007727bp}%
    \ifx\svgscale\undefined%
      \relax%
    \else%
      \setlength{\unitlength}{\unitlength * \real{\svgscale}}%
    \fi%
  \else%
    \setlength{\unitlength}{\svgwidth}%
  \fi%
  \global\let\svgwidth\undefined%
  \global\let\svgscale\undefined%
  \makeatother%
  \begin{picture}(1,0.47155502)%
    \lineheight{1}%
    \setlength\tabcolsep{0pt}%
    \put(0,0){\includegraphics[width=\unitlength,page=1]{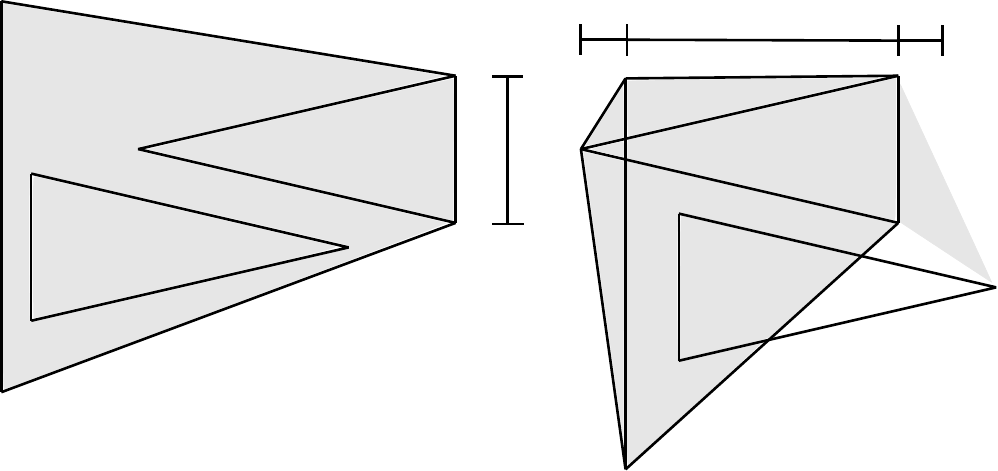}}%
    \put(0.51700894,0.3149798){\color[rgb]{0,0,0}\makebox(0,0)[lt]{\lineheight{1.25}\smash{\begin{tabular}[t]{l}$2z$\end{tabular}}}}%
    \put(0.59209934,0.45471241){\color[rgb]{0,0,0}\makebox(0,0)[lt]{\lineheight{1.25}\smash{\begin{tabular}[t]{l}$x$\end{tabular}}}}%
    \put(0.71374511,0.45207482){\color[rgb]{0,0,0}\makebox(0,0)[lt]{\lineheight{1.25}\smash{\begin{tabular}[t]{l}$2-z-x$\end{tabular}}}}%
    \put(0.91256339,0.45075605){\color[rgb]{0,0,0}\makebox(0,0)[lt]{\lineheight{1.25}\smash{\begin{tabular}[t]{l}$x$\end{tabular}}}}%
    \put(0,0){\includegraphics[width=\unitlength,page=2]{PackingProblem6.pdf}}%
  \end{picture}%
\endgroup%

\caption{Case $dlr$.}
\label{img:PackingProblem6}
\end{figure}

If the vertical line segment does not intersect the $r$-triangle (as on the left in Figure \ref{img:PackingProblem6}), then the trapezoid enclosed by the line segment and the $r$-triangle (which is contained in the convex hull of any such arrangement of the two triangles and the vertical line segment) has volume greater or equal
\beqq \frac{(2+2z)(2-z)}{2}=2+z(1-z)\eeqq
which for $z\in (0,1)$ is greater than the volume of the square.

If the vertical line segment intersects the $r$-triangle, then we denote the horizontal overhang on the left side by $x$ (see the picture on the right in Figure \ref{img:PackingProblem6}). Then, having the $dlr$-case forces the convex hull of every arrangement of the two triangles and the line segment to contain a triangle on the left (with basis-length $2$ and height $x$), a trapezoid in the middle (with parallel-lengths $2$ and $2z$ and height $2-z-x$), and a triangle on the right (with basis-length $2z$ and height greater or equal than $x$). This implies that every convex hull of arrangemets of this type has volume greater or equal
\beqq \frac{2x}{2}+\frac{(2+2z)(2-z-x)}{2}+\frac{2zx}{2}=2+z(1-z)\eeqq
which for $z\in(0,1)$ is greater than the volume of the square.

\underline{Case $ldr$}: Let $w$ be the distance between the two vertical lines corresponding to the $l$- and $r$-triangles. Then, one has $w\leq 2-2z$ or $2-2z\leq w$.

\begin{figure}[h!]
\centering
\def\svgwidth{380pt}
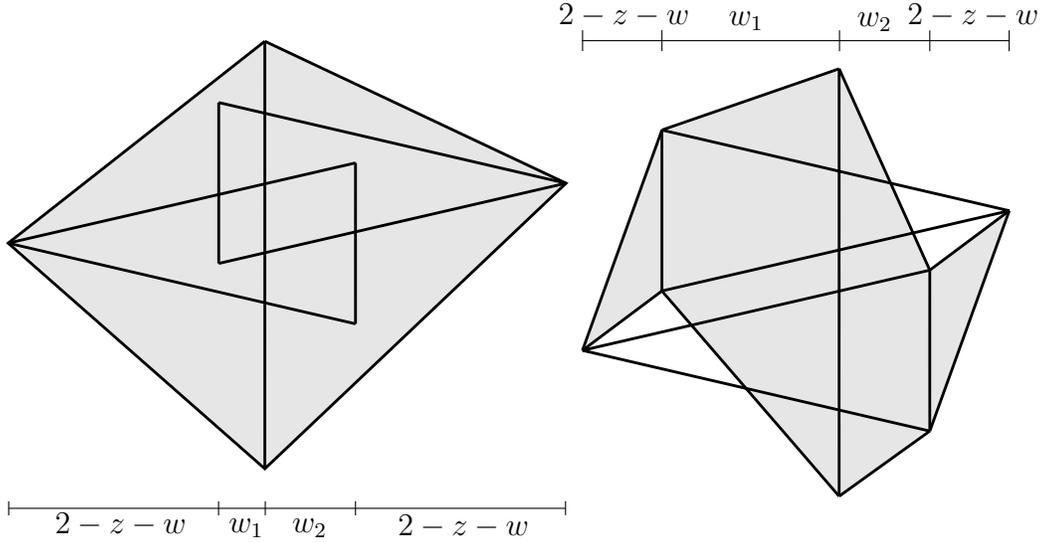
\caption{Case $ldr$ with $w\leq 2-2z$ on the left and $w\geq 2-2z$ on the right.}
\label{img:PackingProblem12}
\end{figure}

If $w\leq 2-2z$, then we are in a situation as shown on the left in Figure \ref{img:PackingProblem12}. We note that the convex hull of the end points of the line segment and the tips of the triangles pointing to the left and right, respectively, is contained in the convex hull of the line segment and the two triangles--for any arrangement. Let $w_1$ and $w_2$ be the distances between the three apparent parallel lines such that $w$ is their sum. Then, the volume of the convex hull of any arrangement of the line segment and the two triangles of this kind is greater or equal
\beqq \frac{2(2-z-w+w_1)}{2}+\frac{2(2-z-w+w_2)}{2}=4-2z-w \geq 2\eeqq
which is the volume of the square.

If $w\geq 2-2z$, then we are in a situation as shown on the right in Figure \ref{img:PackingProblem12}. Therein, the volume of the grey area--which in any arrangement is contained in the convex hull of the line segment and the two triangles--consists of two times the volume of a triangle with basis-length $2z$ and height $2-z-w$ and the volume of two trapezoids with heights $w_1$ and $w_2$ and each with parallel-lengths $2z$ and $2$. Added up, this implies that the volume of the convex hull of any arrangement of the line segment and the two triangles of this kind is greater or equal
\begin{align*}
2\frac{2z(2-z-w)}{2}+\frac{(2+2z)w_1}{2}+\frac{(2+2z)w_2}{2}&=2z(2-z-w)+(1+z)w\\
&=4z-2z^2 +w(1-z)\\
&\geq 4z-2z^2+(2-2z)(1-z)\\
&=2
\end{align*}
which is the volume of the square (in the estimate we used $z\in(0,1)$).

\begin{figure}[h!]
\centering
\def\svgwidth{400pt}
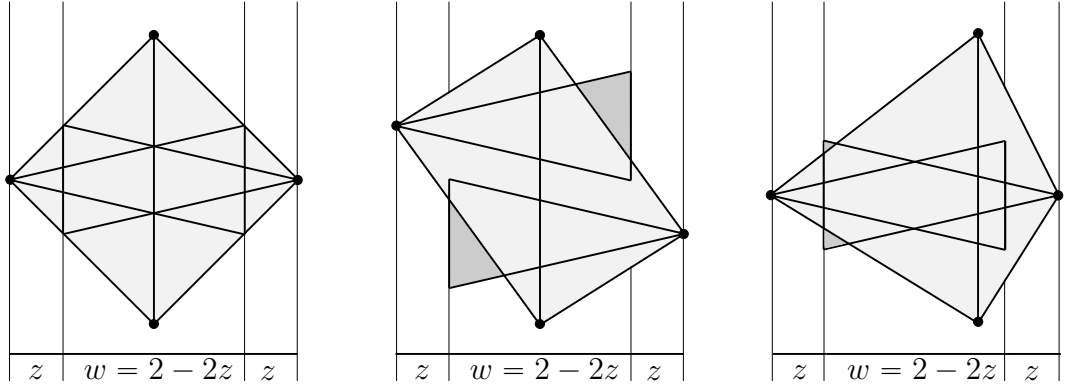
\caption{Case $ldr$ with $w=2-2z$.}
\label{img:PackingProblem13}
\end{figure}

It remains to prove that the square is the unique volume-minimizing convex hull (up to translations) of translates of the line segment and the two triangles. Based on the above investigations, we come to the conclusion that volume-minimizing convex hulls can only appear within the $ldr$-case when $w=2-2z$. Then, for any arrangement of the two triangles and the line segment under these circumstances the convex hull of the tips of the two triangles and the end points of the line segment is contained in the convex hull of the two triangles and the line segment and has volume $2$. Now, one notes that whenever the arrangement of the tips and end points does not coincide with the square-case, then the convex hull of the two tips and the end points is strictly contained in the cinvex hull of the two triangles and the line segement (see Figure \ref{img:PackingProblem13}).
\epf

Unfortunately, we are not able to prove the condition in Proposition \ref{Prop:quadrilateral12}(i). Nevertheless, we remark that Proposition \ref{Prop:quadrilateral12}(i) can be expressed as conjecture for a Euclidean covering/container problem which implies Viterbo's conjecture for all convex quadrilaterals which are not trapezoids:

\bconj
Let $a_1,a_2\in[0,1]$ be given such that $Q_{(\alpha,\beta,a_1,a_2,d_1,d_2,c_1,c_2)}$ is a convex quadrilateral in $\R^2$ which is not a trapezoid. Then, $\square$ is a volume-minimizing convex hull of translates of
\beqq \pm\lambda_1 J\Delta_{a_1,a_2,1}\; \text{ and }\;\pm\lambda_2 J\Delta_{a_1,a_2,2}.\eeqq
\econj

\begin{figure}[h!]
\centering
\def\svgwidth{410pt}
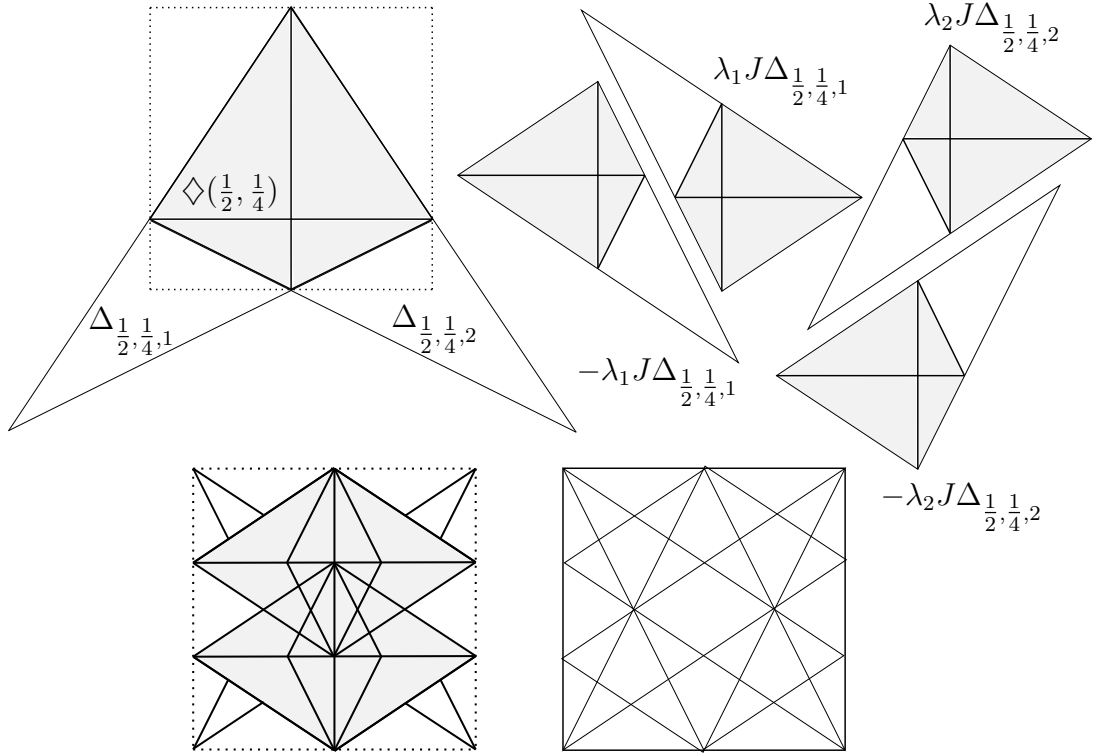
\caption[Illustration of the covering problem for $a_1=\frac{1}{2}$ and $a_2=\frac{1}{4}$.]{Starting from $\diamondsuit(\frac{1}{2},\frac{1}{4})$, we construct $\Delta_{\frac{1}{2},\frac{1}{4},1}$ and $\Delta_{\frac{1}{2},\frac{1}{4},2}$. Based on that, we find the corresponding $\lambda_1$ and $\lambda_2$ and get the four convex bodies $\pm\lambda_1 J \Delta_{\frac{1}{2},\frac{1}{4},1}$ and $\pm\lambda_2 J \Delta_{\frac{1}{2},\frac{1}{4},2}$ whose convex hull should be minimized in terms of its volume. If the square, which by construction is a convex hull of translates of the four aforementioned convex bodies, is volume-minimizing, then Viterbo's conjecture is true for all Lagrangian products $Q_{(\alpha,\beta,\frac{1}{2},\frac{1}{4},d_1,d_2,c_1,c_2)}\times T$, where $T$ is any convex body in $\R^2$.}
\label{img:greatdragon}
\end{figure}

In order to understand the difficulty of this covering/container problem, we will briefly present the following example: Let us consider convex quadrilaterals $Q_{(\alpha,\beta,a_1,a_2,d_1,d_2,c_1,c_2)}$ with $a_1=\frac{1}{2}$ and $a_2=\frac{1}{4}$. Then, Viterbo's conjecture is true for all Lagrangian products
\beqq Q_{(\alpha,\beta,a_1,a_2,d_1,d_2,c_1,c_2)} \times T,\eeqq
where $T$ can be any convex body in $\R^2$, if $\square$ is a volume-minimizing convex hull of translates of
\beqq \pm \lambda_1 J\Delta_{\frac{1}{2},\frac{1}{4},1}\; \text{ and } \; \pm \lambda_1 J\Delta_{\frac{1}{2},\frac{1}{4},1}.\eeqq
Figure \ref{img:greatdragon} illustrates this covering problem.

\section[Proof of Corollary \ref{Cor:sharpinequalities}]{Proof of Corollary \ref{Cor:sharpinequalities}}\label{Sec:equacaCorollaries}

From Theorem \ref{Thm:answer1}, we derive that
\beq \vol(K\times Q) \geq \frac{c_{EHZ}(K\times Q)^2}{2}\label{eq:Corshin1}\eeq
holds for all trapezoids $Q\subset\R^2$ and all convex bodies $K\subset\R^2$. Furthermore, from Theorem \ref{Thm:answer2}(iii), we conclude that this inequality is sharp. Assuming the volume of $K$ to be $1$ and referring to Theorem \ref{Thm:Chap7relationship}, we can rewrite \eqref{eq:Corshin1} by
\beq \max_{\vol(K)=1}\; \min_{q\in M_3(K,Q)}\ell_Q(q) \leq \sqrt{2}\label{eq:Corshin2}\eeq
which by \cite[Theorem 1.4]{Rudolf2022b} is equivalent to
\beq \max_{K\in A\left(Q,\sqrt{2}\right)}\vol(K) \geq 1.\label{eq:Corshin3}\eeq
The sharpness of \eqref{eq:Corshin2} and \eqref{eq:Corshin3} follows from the sharpness of \eqref{eq:Corshin1} together with the remarks in \cite[Theorem 1.4]{Rudolf2022b} concerning the equality cases.

\section{Zoll-property of equality cases}\label{Sec:Zollpropertys}

We begin with the proof of Theorem \ref{Thm:Zollproperty}:

\bpf[Proof of Theorem \ref{Thm:Zollproperty}]
First of all, we notice that it is enough to investigate the equality cases for the triangle- and parallelogram-configurations. This is due to the fact that the known equality cases of the trapezoid- and general-convex-quadrilateral-configurations are reflected copies (by interchanging $\R^2(x)$ and $\R^2(y)$) of certain equality cases within the parallelogram-configuration. In a similar way, we can argue that for the triangle-configuration it is enough to consider only the triangle-hexagon-cases since the triangle-parallelogram-cases are included within the equality cases of the parallelogram-configuration.

We begin by investigating regular Minkowski billiard trajectories whose dual billiard trajectories are also regular--in what follows we will call them just regular Minkowski billiard trajectories--by applying the Minkowski billiard reflection rule.

\begin{figure}[h!]
\centering
\def\svgwidth{390pt}
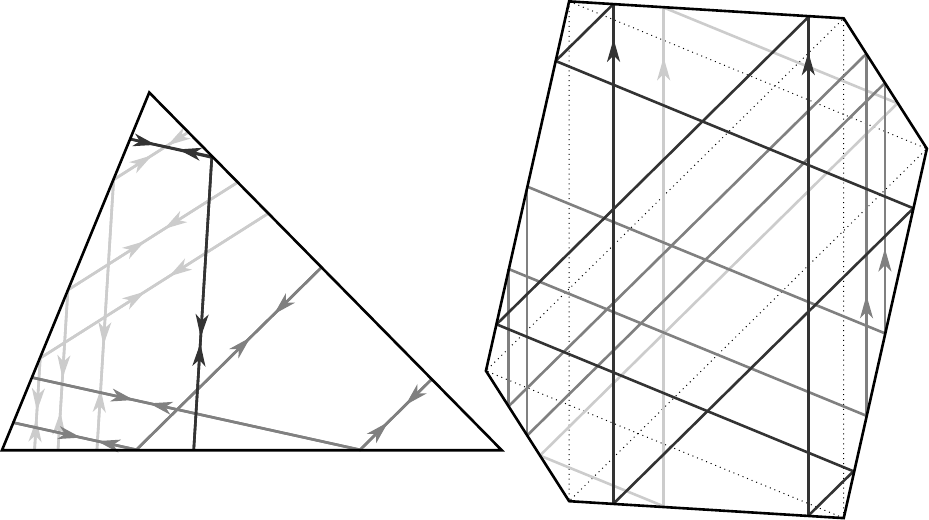
\caption{Illustration of some the regular closed $(\Delta,T)$-Minkowski billiard trajectories whose closed dual billiard trajectories are regular as well.}
\label{img:zollpropertytriangle}
\end{figure}

\begin{figure}[h!]
\centering
\def\svgwidth{390pt}
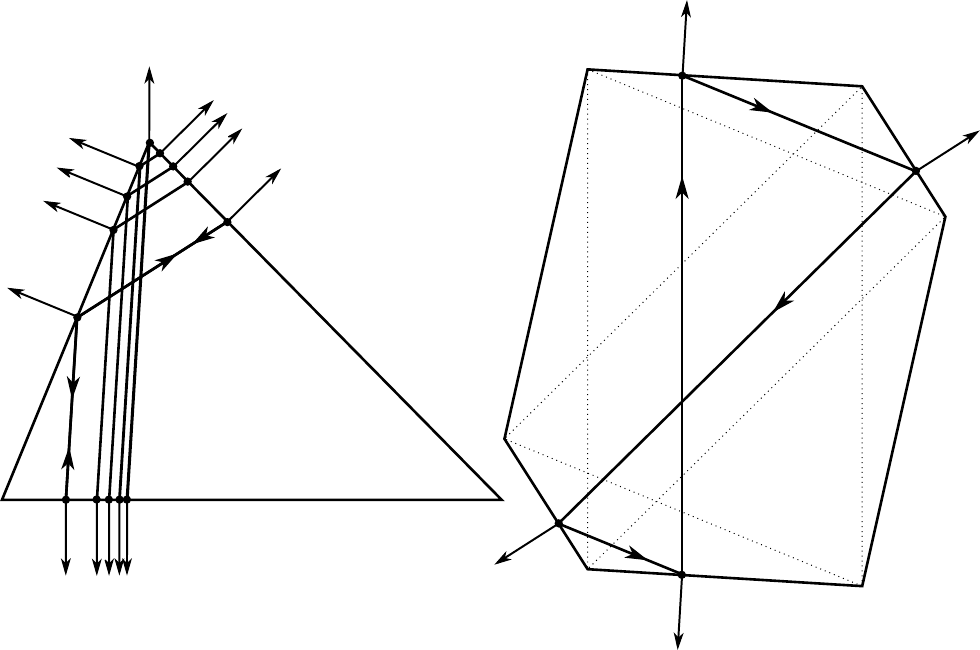
\caption{On the left, we see a sequence of regular closed $(\Delta,T)$-Minkowski billiard trajectories that approximates an $\ell_T$-minimizing closed $(\Delta,T)$-Minkowski billiard trajectory with two bouncing points. On the right, we see a corresponding regular closed dual billiard.}
\label{img:zollpropertytriangle1}
\end{figure}

Let us begin with the triangle-hexagon-configuration. One realizes that the respective regular $(\Delta,T)$-Minkowski billiard trajectories, where $\Delta$ is any triangle in $\R^2$ and $T$ a corresponding equality-case-hexagon, are of the form as shown in Figure \ref{img:zollpropertytriangle}.
The following aspects can be immediately derived from this pattern: Every (simple) regular $(\Delta,T)$-Minkowski billiard trajectory is forced to be closed and have exactly $4$ bouncing points. Furthermore, in $\Delta$ as well as in $T$, the union of all possible regular $(\Delta,T)$-Minkowski billiard and dual billiard trajectories is dense, respectively. It remains to justify that the regular $(\Delta,T)$-Minkowski billiard trajectories are actually (i) closed in the way Figure \ref{img:zollpropertytriangle} suggests and (ii) $\ell_T$-minimizing. The closedness of the $(\Delta,T)$-Minkowski billiard trajectories is clear by construction; the closedness of the dual billiard trajectories in $T$ follows from the subsequent Lemma \ref{Lem:closednesstrianglehexagon}. The $\ell_T$-minimality of the regular $(\Delta,T)$-Minkowski billiard trajectories follows from the continuity of the $\ell_T$-length-functional with respect to small perturbations of the trajectories and from the fact that there are straight lines connecting the vertices of $\Delta$ with their opposite sides which can be approximated by sequences of closed regular $(\Delta,T)$-Minkowski billiard trajectories as shown in Figure \ref{img:zollpropertytriangle1} and from which it is known by previous considerations that they are $\ell_T$-minimizing. Here, it is important to notice that when looking at one form of regular $(\Delta,T)$-Minkowski billiad trajectories, then they all have the same $\ell_T$-length since due to the Minkowski billiard reflection rule and the fact that their bouncing points are on the same sides of $\Delta$ they all can be associated to one closed dual billiard trajectory in $T$ (see for instance Figure \ref{img:zollpropertytriangle1}).

\begin{figure}[h!]
\centering
\def\svgwidth{310pt}
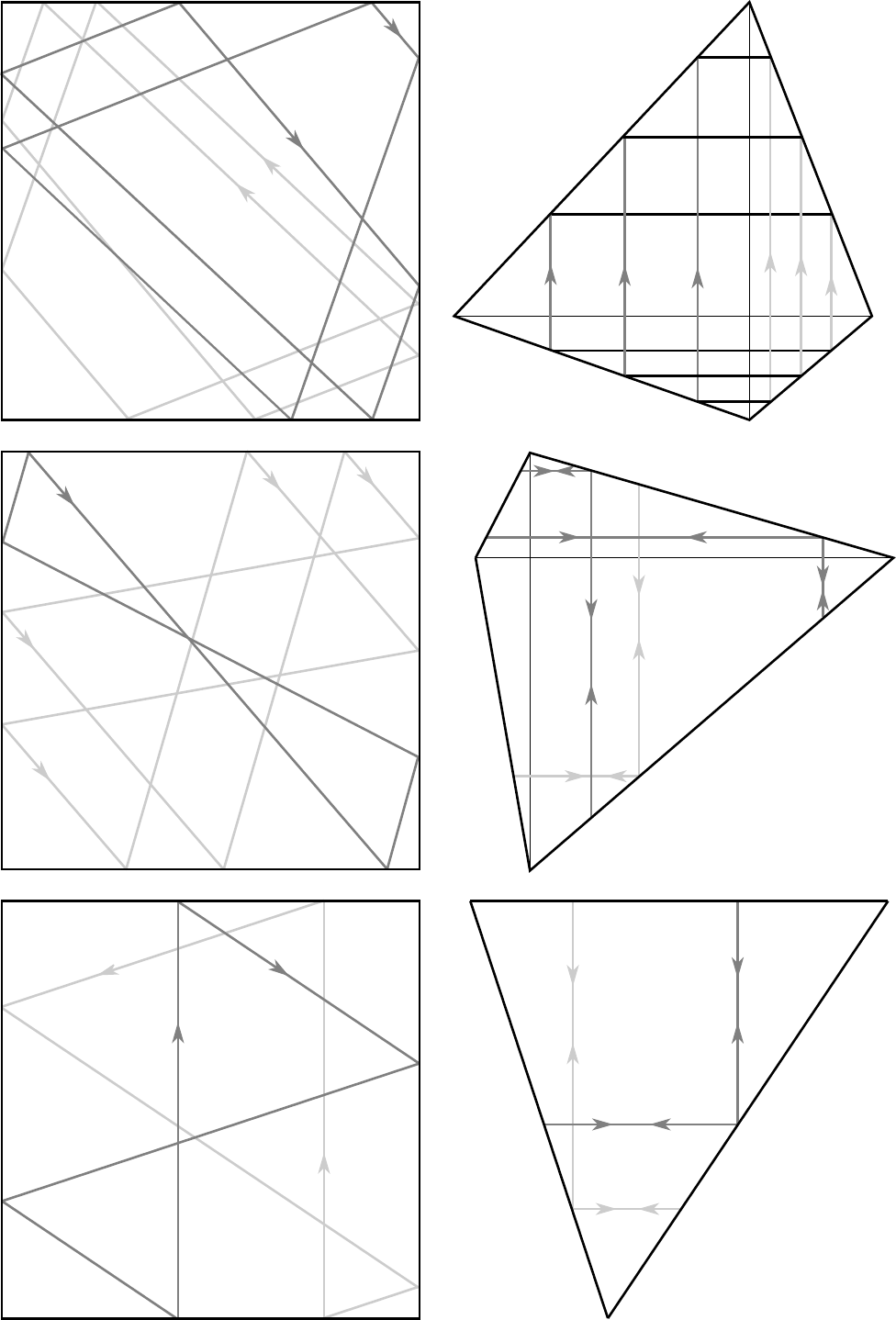
\caption{Illustration of three different patterns of regular $(\square,\diamondsuit)$-Minkowski billiard trajectories that fully represent the qualitative differences between the respective regular $(\square,\diamondsuit)$-Minkowski billiard trajectories and their regular dual trajectories in $\diamondsuit$.}
\label{img:zollpropertysquare}
\end{figure}

\begin{figure}[h!]
\centering
\def\svgwidth{390pt}
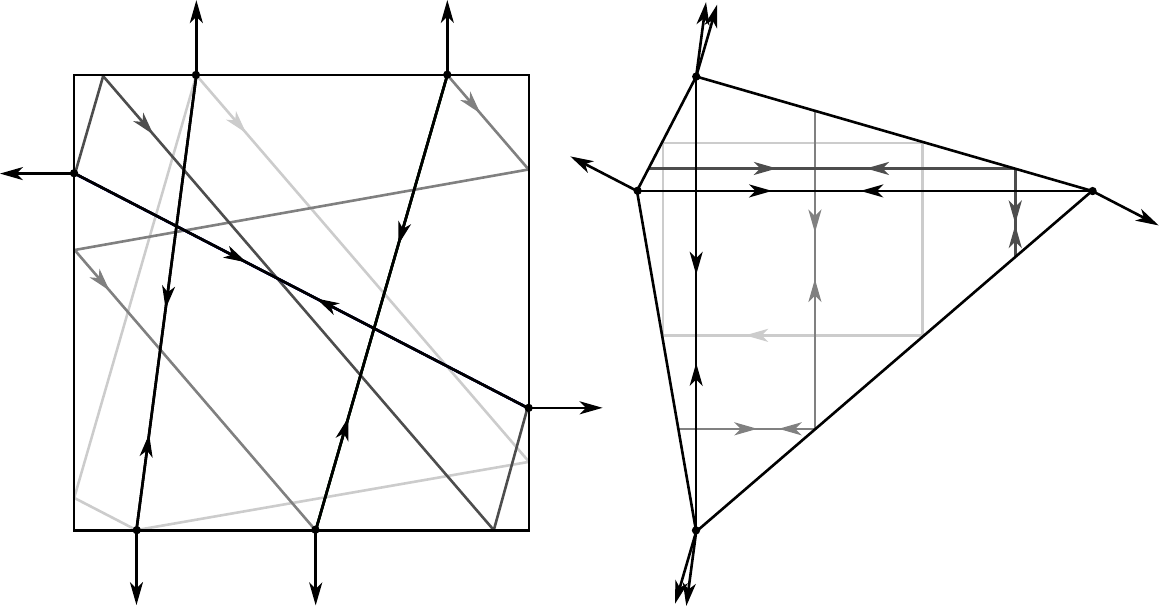
\caption{Illustration of the $\ell_\diamondsuit$-preserving modifications of three regular closed $(\square,\diamondsuit)$-Minkowski billiard trajectories (the red, blue, and green one). The three resulting closed $(\square,\diamondsuit)$-Minkowski billiard trajectories with two bouncing points are each $\ell_\diamondsuit$-minimizing.} 
\label{img:zollpropertysquare1}
\end{figure}

Now, we come to the parallelogram-configurations, while it actually suffices to investigate the square-diamond-configurations since all other can be traced back to them by applying suitable symplectomorphisms. As before, one realizes that the respective regular $(\square,\diamondsuit)$-Minkowski billiard trajectories, where $\square$ is a square in $\R^2$ and $\diamondsuit$ any corresponding equality-case-diamond, are of the form as shown in Figure \ref{img:zollpropertysquare}. In this figure, we indicate three different patterns of regular $(\square,\diamondsuit)$-Minkowski billiard trajectories that fully represent the qualitative differences between the respective regular $(\square,\diamondsuit)$-Minkowski billiard trajectories and their regular dual trajectories in $\diamondsuit$.
As before, the following aspects can be immediately derived from these patterns: Every (simple) regular $(\square,\diamondsuit)$-Minkowski billiard trajectory is forced to be closed and have exactly $4$ bouncing points. Furthermore, in $\square$ as well as in $\diamondsuit$, the union of all possible regular $(\square,\diamondsuit)$-Minkowski billiard trajectories and dual billiard trajctories is dense, respectively. As before, it remains to justify that the regular $(\square,\diamondsuit)$-Minkowski billiard trajectories are actually (i) closed in the way Figure \ref{img:zollpropertysquare} suggests and (ii) $\ell_\diamondsuit$-minimizing. The $\ell_\diamondsuit$-minimality follows from the fact that the regular closed $(\square,\diamondsuit)$-Minkowski billiard trajectories can be modified $\ell_\diamondsuit$-preservingly (here comes Lemma 5.4 in \cite{Rudolf2022b} into play; briefly summarized, one uses the property of the Minkowski functional that a polygonal line composed of directions belonging to the (one-dimensional) normal cones of two neighbouring facets of $\diamondsuit$ can be replaced by the line segment connecting start and end point of the polygonal line (which means that its direction is in the normal cone which is enclosed by the two aforementioned directions) while the $\ell_\diamondsuit$-length is preserved) by skipping some of the bouncing points in order to reach at closed $(\square,\diamondsuit)$-Minkowski billiard trajectories with two bouncing points from which we already know from previous considerations that they are $\ell_\diamondsuit$-minimizing. Figure \ref{img:zollpropertysquare1} shows these modifications for different instances.
The closedness of the regular $(\square,\diamondsuit)$-Minkowski billiard trajectories can be shown as follows: First, the closedness of the regular dual billiard trajectories in $\diamondsuit$ in the middle and lower picture of Figure \ref{img:zollpropertysquare} is clear by construction. The closedness of the regular $(\square,\diamondsuit)$-Minkowski billiard trajectories in the lower picture of Figure \ref{img:zollpropertysquare} follows from the subsequent Lemma \ref{Lem:closednesslower}. The closedness of the regular $(\square,\diamondsuit)$-Minkowski billiard trajectories in the upper picture of Figure \ref{img:zollpropertysquare} follows from the subsequent Lemma \ref{Lem:closednessupper}. And finally, the closedness of the regular $(\square,\diamondsuit)$-Minkowski billiard trajectories in the middle picture of Figure \ref{img:zollpropertysquare} follows from the subsequent Lemma \ref{Lem:closednessmiddlegreenblue}.
\epf

It remains to prove the lemmata which we mentioned within the proof of Theorem \ref{Thm:Zollproperty}. Although, certainly, there are nice geometric properties, to the author's knowledge, for most of the patterns, there are no obvious geometric proofs for the stated closedness of the respective regular Minkowski billiard trajectories and their dual billiard trajectories. So, for these patterns, we will prove the closedness algebraically.

\blem\label{Lem:closednesstrianglehexagon}
Let $\Delta$ be any triangle in $\R^2$ and $T$ a corresponding equality-case-hexagon. Then, the regular dual billiard trajectories in $T$ which correspond to the $(\Delta,T)$-Minkowski billiard trajectories as shown in Figure \ref{img:zollpropertytriangle} are closed.
\elem

\bpf
\begin{figure}[h!]
\centering
\def\svgwidth{300pt}
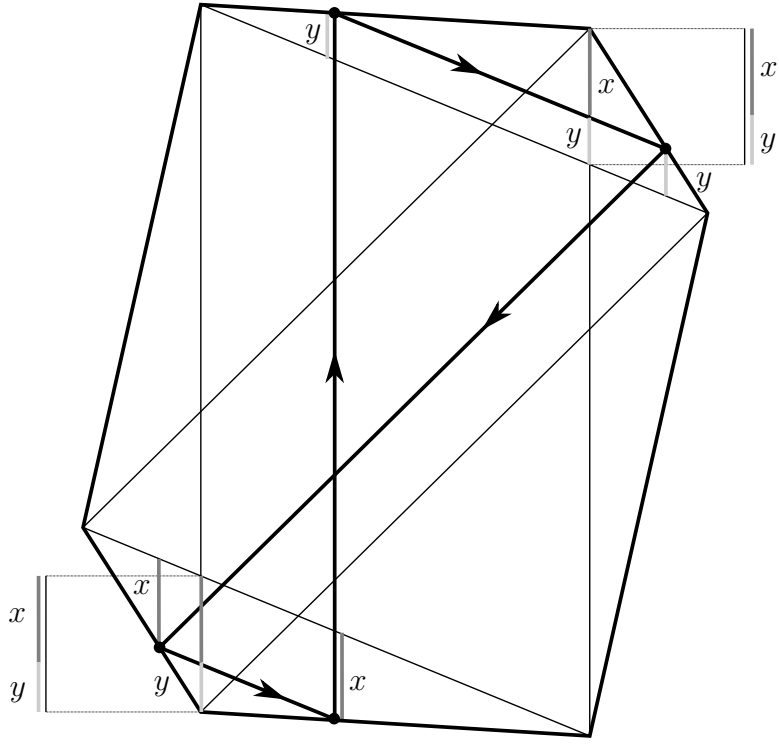
\caption{Illustration of a regular dual billiard trajectory which corresponds to a regular closed $(\Delta,T)$-Minkowski billiard trajectory, where $\Delta$ is a triangle and $T$ a corresponding equality-case-hexagon. The construction of $T$, its resulting symmetry, and the construction of the dual billiard trajectory in $T$ imply the cndition $x+y=\const$ which forces the dual billiard trajectory to be closed.} 
\label{img:closednesshexagon}
\end{figure}
The closedness of the regular dual billiard trajectories in $T$ follows from the central symmetry of $T$ together with the facts that, on the one hand, $T$ is the convex hull of $J\Delta$ and $-J\Delta+t$ for a certain $t\in \R^2$, and that, on the other hand, the directions of the trajectory segments of the dual billiard trajectories are parallel to the sides of the $J$- and $(-J)$-rotated copies of $\Delta$. More precisely and referring to Figure \ref{img:closednesshexagon}, due to the symmetry of $T$, there are certain vertical distances--which we denoted by $x$ and $y$--that satisfy the condition $x+y= \const$ This condition forces the dual billiard trajectory to be closed.
\epf

\blem\label{Lem:closednesslower}
The regular $(\square,\diamondsuit)$-Minkowski billiard trajectories in the lower picture of Figure \ref{img:zollpropertysquare}--these are all configurations for which $\diamondsuit$ is a $\diamondsuit(a_1,a_2)$ with $a_1\in\{0,1\}$ or $a_2\in\{0,1\}$--are closed. 
\elem

\bpf
\begin{figure}[h!]
\centering
\def\svgwidth{370pt}
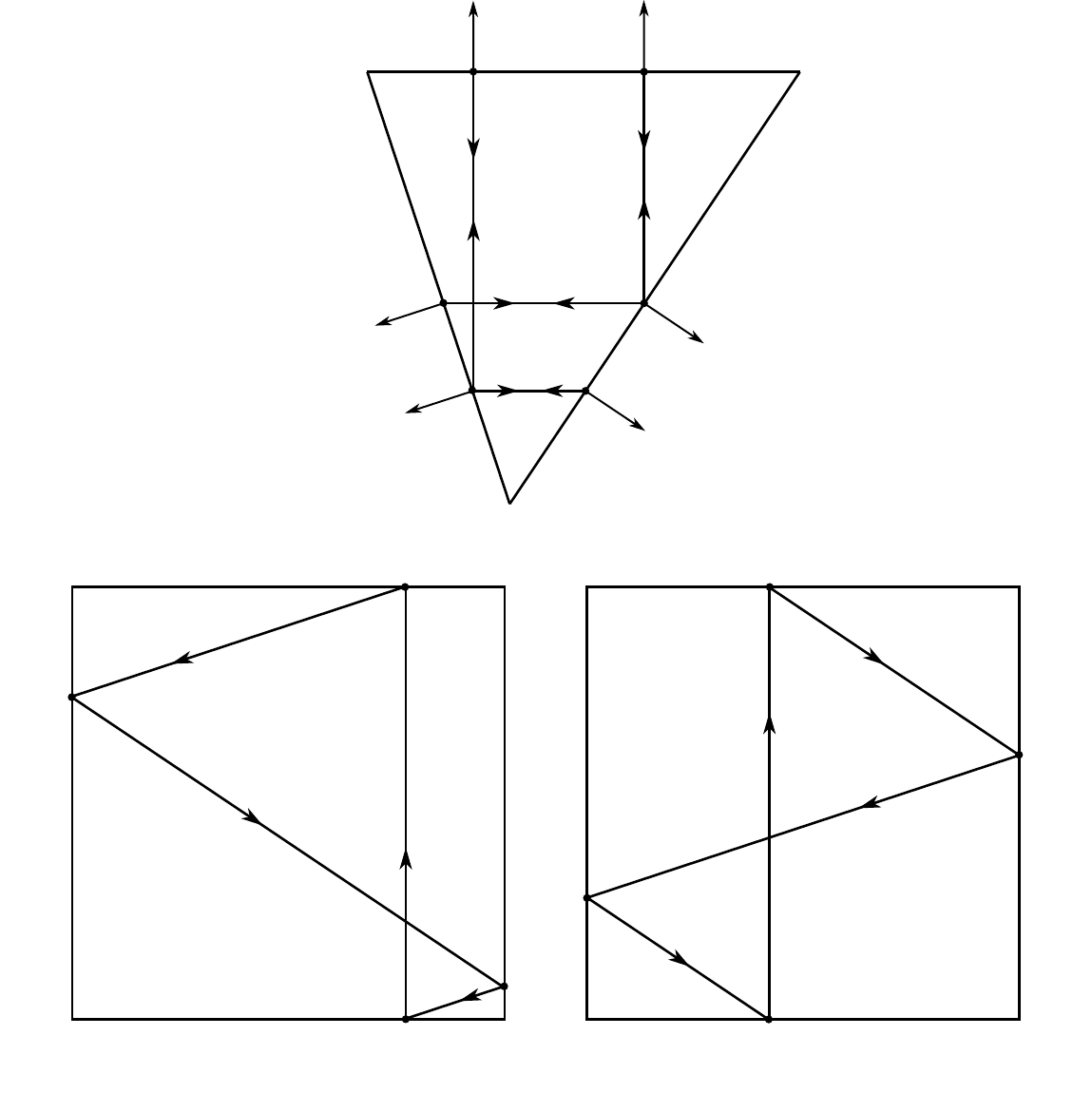
\caption{For both regular $(\square,\diamondsuit(a_1,a_2))$-Minkowski billiard trajectories $q$ and $q'$ (with regular dual billiard trajectories $p$ and $p'$) we have that for every $z_1\in (0,1)$, after four times application the $(\square,\diamondsuit(a_1,a_2))$-Minkowski billiard trajectory map maps $(z_1,0)$ to itself.}
\label{img:lastchange}
\end{figure}
We extend the respective relevant line segments of the regular $(\square,\diamondsuit(a_1,a_2))$-Minkowski billiard trajectories to straight lines and call them $g_1$, $g_2$, and $g_3$ as indicated in Figure \ref{img:lastchange}, respectively. Based on $a_1,a_2\in [0,1]$ (here: $a_2=1$ and $a_1\in(0,1)$) and $z_1\in (0,1)$, we can represent them algebraically, respectively.

For the lower right situation in Figure \ref{img:lastchange} we get
\beqq g_1(x)=(a_1-1)x+1-z_1(a_1-1),\quad g_2(x)=a_1x+z_1(1-a_1),\eeqq
\beqq g_3(x)=(a_1-1)x+z_1(1-a_1);\eeqq
for the lower left:
\beqq g_1(x)=a_1x+1-a_1z_1,\quad g_2(x)=(a_1-1)x+1-a_1z_1,\eeqq
\beqq g_3(x)=a_1x-a_1z_1.\eeqq

If we denote the projections along $g_i$, $i\in\{1,2,3\}$, onto a straight line $l(s)$, $s\in\R$, by $\pi_{g_i}^{l(s)}$, then, for the lower right picture in Figure \ref{img:lastchange}, we have:
\allowdisplaybreaks{\begin{align*}
\left(\pi_{g_3}^{(s,0)} \circ \pi_{g_2}^{(0,s)}\circ \pi_{g_1}^{(1,s)}\right)(z_1,1)=& \left(\pi_{g_3}^{(s,0)} \circ \pi_{g_2}^{(0,s)}\right)(1,a_1(1-z_1)+z_1)\\
=& \pi_{g_3}^{(s,0)} (0,z_1(1-a_1))\\
=& (z_1,0);
\end{align*}}%
and for the lower left:
\allowdisplaybreaks{\begin{align*}
\left(\pi_{g_3}^{(s,0)} \circ \pi_{g_2}^{(1,s)}\circ \pi_{g_1}^{(0,s)}\right)(z_1,1)=& \left(\pi_{g_3}^{(s,0)} \circ \pi_{g_2}^{(1,s)}\right)(0,1-a_1z_1)\\
=& \pi_{g_3}^{(s,0)} (1,a_1-a_1z_1)\\
=& (z_1,0);
\end{align*}}%

This immediately implies the closedness of the regular $(\square,\diamondsuit)$-Minkowski billiard trajectories in the lower picture of Figure \ref{img:zollpropertysquare}.
\epf

\blem\label{Lem:closednessupper}
The regular $(\square,\diamondsuit)$-Minkowski billiard trajectories in the upper picture of Figure \ref{img:zollpropertysquare} are closed.
\elem

\bpf
\begin{figure}[h!]
\centering
\def\svgwidth{390pt}
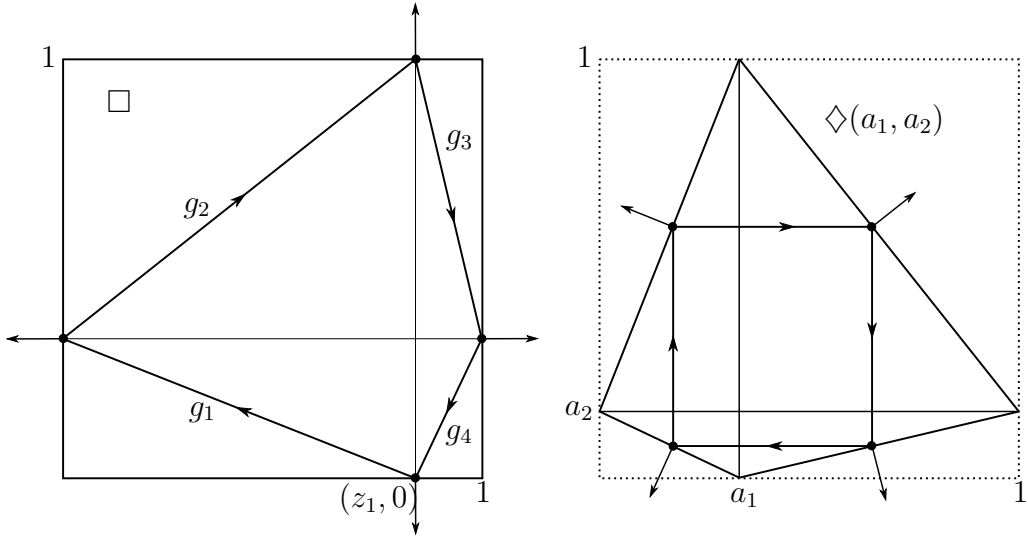
\caption{For every $z_1\in (0,1)$, after four times application the $(\square,\diamondsuit(a_1,a_2))$-Minkowski billiard trajectory map maps $(z_1,0)$ to itself.} 
\label{img:closednessupper}
\end{figure}
We extend the respective line segments of the regular $(\square,\diamondsuit(a_1,a_2))$-Minkowski billiard trajectories to straight lines and call them $g_1$, $g_2$, $g_3$, and $g_4$ as indicated in Figure \ref{img:closednessupper}. Based on $a_1,a_2\in (0,1)$ and $z_1\in (0,1)$, we can represent them algebraically:
\beqq g_1(x)=\frac{-a_1}{1-a_2}x+\frac{z_1a_1}{1-a_2},\quad g_2(x)=\frac{a_1-1}{a_2-1}x+\frac{z_1a_1}{1-a_2},\eeqq
\beqq g_3(x)=\frac{a_1-1}{a_2}x+\frac{1-z_1a_1}{a_2},\quad g_4(x)=\frac{a_1}{a_2}x -\frac{a_1z_1}{a_2}.\eeqq
If we denote the projections along $g_i$, $i\in\{1,2,3,4\}$, onto a straight line $l(s)$, $s\in\R$, by $\pi_{g_i}^{l(s)}$, then we have
\allowdisplaybreaks{\begin{align*}
\left(\pi_{g_4}^{(s,0)} \circ \pi_{g_3}^{(1,s)} \circ \pi_{g_2}^{(s,1)}\circ \pi_{g_1}^{(0,s)}\right)(z_1,0)=&\left(\pi_{g_4}^{(s,0)} \circ \pi_{g_3}^{(1,s)} \circ \pi_{g_2}^{(s,1)}\right) \left(0,\frac{z_1a_1}{1-a_2}\right)\\
=& \left(\pi_{g_4}^{(s,0)} \circ \pi_{g_3}^{(1,s)}\right) \left(\frac{a_2-1+z_1a_1}{a_1-1},1\right)\\
=& \pi_{g_4}^{(s,0)} \left(1,\frac{a_1-a_1z_1}{a_2}\right)\\
=& (z_1,0).
\end{align*}}%
This immediately implies the closedness of all regular $(\square,\diamondsuit)$-Minkowski billiard trajectories in the upper picture of Figure \ref{img:zollpropertysquare}.
\epf

\blem\label{Lem:closednessmiddlegreenblue}
The regular $(\square,\diamondsuit)$-Minkowski billiard trajectories in the middle picture of Figure \ref{img:zollpropertysquare} are closed.
\elem

\bpf
\begin{figure}[h!]
\centering
\def\svgwidth{390pt}
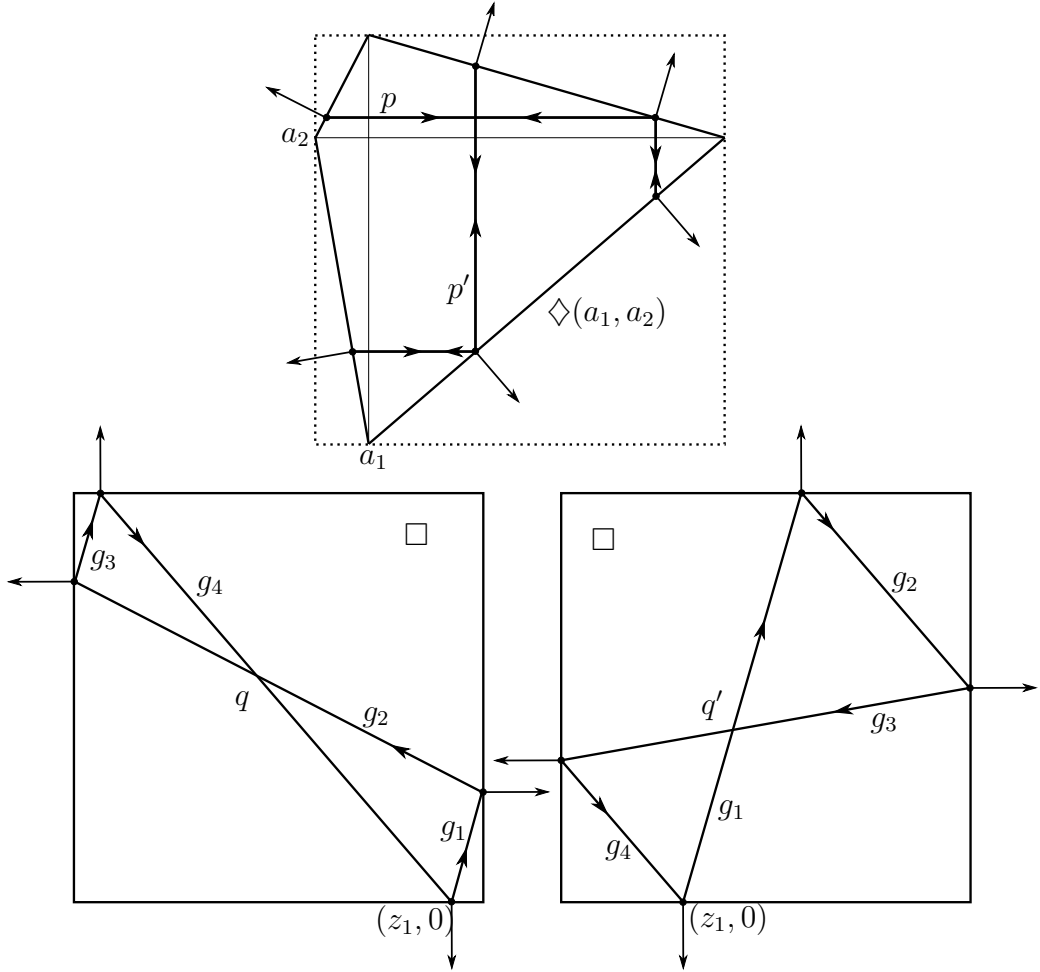
\caption{For both regular $(\square,\diamondsuit(a_1,a_2))$-Minkowski billiard trajectories $q$ and $q'$ (with regular dual billiard trajectories $p$ and $p'$) we have that for every $z_1\in (0,1)$, after four times application the $(\square,\diamondsuit(a_1,a_2))$-Minkowski billiard trajectory map maps $(z_1,0)$ to itself.} 
\label{img:closednessmiddlegreenblue}
\end{figure}
We extend the respective line segments of the regular $(\square,\diamondsuit(a_1,a_2))$-Minkowski billiard trajectories to straight lines and call them $g_1$, $g_2$, $g_3$, and $g_4$ as indicated in Figure \ref{img:closednessmiddlegreenblue}, respectively. Based on $a_1,a_2\in (0,1)$ and $z_1\in (0,1)$, we can represent them algebraically, respectively.

For the lower right situation in Figure \ref{img:closednessmiddlegreenblue} we get
\beqq g_1(x)=\frac{-a_1}{1-a_2}x+\frac{a_1z_1}{1-a_2},\quad g_2(x)=\frac{a_1-1}{a_2-1}a+\frac{a_1z_1}{1-a_2},\eeqq
\beqq g_3(x)=\frac{-a_1}{1-a_2}x +\frac{1+a_1z_1}{1-a_2},\quad g_4(x)=\frac{a_1}{a_2}x -\frac{a_1z_1}{a_2};\eeqq
for the lower left:
\beqq g_1(x)=\frac{a_1-1}{a_2-1}x+\frac{z_1-a_1z_1}{a_2-1},\quad g_2(x)=\frac{-a_1}{1-a_2}x+\frac{-1+z_1-a_1z_1}{a_2-1},\eeqq
\beqq g_3(x)=\frac{a_1-1}{a_2-1}x +\frac{-1+z_1-a_1z_1}{a_2-1},\quad g_4(x)=\frac{a_1-1}{a_2}x +\frac{z_1-a_1z_1}{a_2}.\eeqq

If we denote the projections along $g_i$, $i\in\{1,2,3,4\}$, onto a straight line $l(s)$, $s\in\R$, by $\pi_{g_i}^{l(s)}$, then, for the lower right picture in Figure \ref{img:closednessmiddlegreenblue}, we have:
\allowdisplaybreaks{\begin{align*}
\left(\pi_{g_4}^{(s,0)} \circ \pi_{g_3}^{(s,1)} \circ \pi_{g_2}^{(1,s)}\circ \pi_{g_1}^{(0,s)}\right)(z_1,0)=& \left(\pi_{g_4}^{(s,0)} \circ \pi_{g_3}^{(s,1)} \circ \pi_{g_2}^{(1,s)}\right)\left(0,\frac{a_1z_1}{1-a_2}\right)\\
=& \left(\pi_{g_4}^{(s,0)} \circ \pi_{g_3}^{(s,1)} \right)\left(1,\frac{1-a_1+a_1z_1}{1-a_2}\right)\\
=& \pi_{g_4}^{(s,0)} \left(\frac{a_2+a_1z_1}{a_1},1\right)\\
=& (z_1,0);
\end{align*}}%
and for the lower left:
\allowdisplaybreaks{\begin{align*}
 \left(\pi_{g_4}^{(s,0)} \circ \pi_{g_3}^{(s,1)} \circ \pi_{g_2}^{(0,s)}\circ \pi_{g_1}^{(1,s)}\right)(z_1,0)=& \left(\pi_{g_4}^{(s,0)} \circ \pi_{g_3}^{(s,1)} \circ \pi_{g_2}^{(0,s)}\right) \left(1,\frac{a_1-1+z_1-a_1z_1}{a_2-1}\right)\\
=& \left(\pi_{g_4}^{(s,0)} \circ \pi_{g_3}^{(s,1)} \right)\left(0,\frac{-1+z_1-a_1z_1}{a_2-1}\right)\\
=& \pi_{g_4}^{(s,0)} \left(\frac{a_2-z_1+a_1z_1}{a_1-1},1\right)\\
=&(z_1,0).
\end{align*}}%

This immediately implies the closedness of the regular $(\square,\diamondsuit)$-Minkowski billiard trajectories in the middle picture of Figure \ref{img:zollpropertysquare}.
\epf

\blem\label{Lem:closednessdualuppermiddle}
The regular dual billiard trajectories in the upper picture of Figure \ref{img:zollpropertysquare} are closed.
\elem

\begin{figure}[h!]
\centering
\def\svgwidth{235pt}
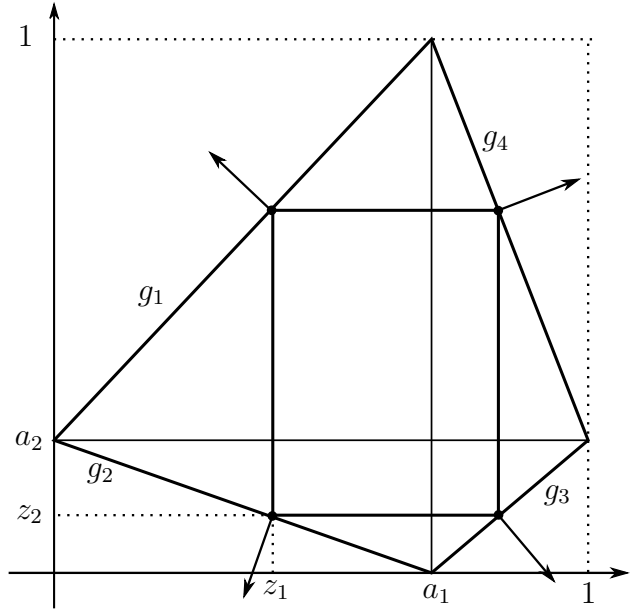
\caption{For every $z=(z_1,z_2)$ on $g_2 \cap \partial \diamondsuit$ with $z_1\in (0,a_1)$ which we take as starting point of a dual billiard trajectory in $\diamondsuit$, one has that with the $4$-th reflection the trajectory is back in $z$.} 
\label{img:zollpropertysquare3}
\end{figure}

\bpf
We consider any diamond $\diamondsuit(a_1,a_2)$ with $a_1,a_2\in (0,1)$. We understand the facets of $\diamondsuit(a_1,a_2)$ as straight lines $g_1$, $g_2$, $g_3$, and $g_4$. Furthermore, we denote the lower left vertex of the regular dual billiard trajectory (which is a rectangle) by $z=(z_1,z_2)$ (see Figure \ref{img:zollpropertysquare3}).
Now, we let $\pi_1^{g_1}$ and $\pi_1^{g_3}$ be the projections along the horizontal axis onto $g_1$ and $g_3$, respectively. Similarly, we let $\pi_2^{g_4}$ be the projections along the vertical axis onto $g_4$. Then, we can represent the straight lines $g_1$, $g_2$, $g_3$, and $g_4$ as graphs in $\R^2$:
\beqq g_1 (x)=a_2+\frac{1-a_2}{a_1}x, \quad g_2(x)=a_2-\frac{a_2}{a_1}x,\eeqq
\beqq g_3(x)=\frac{-a_1a_2}{1-a_1}+\frac{a_2}{1-a_1}x,\quad g_4(x)=\frac{1-a_1a_2}{1-a_1}-\frac{1-a_2}{1-a_1}x.\eeqq
Then, we just calculate:
{\allowdisplaybreaks\begin{align*}
\left( \pi_1^{g_1} \circ \pi_2^{g_4} \circ \pi_1^{g_3}\right) (z_1,z_2)=& \left( \pi_1^{g_1} \circ \pi_2^{g_4} \circ \pi_1^{g_3}\right) \left(z_1,a_2-\frac{a_2}{a_1}z_1\right)\\
=& \left( \pi_1^{g_1} \circ \pi_2^{g_4} \right) \left(1+\frac{a_1-1}{a_2}z_1,a_2-\frac{a_2}{a_1}z_1\right)\\
=&\pi_1^{g_1}\left(1+\frac{a_1-1}{a_2}z_1,\frac{a_2-a_1a_2}{1-a_1}+\frac{1-a_1-a_2+a_1a_2}{a_1(1-a_1)}z_1\right)\\
=&\left(z_1,\frac{a_2-a_1a_2}{1-a_1}+\frac{1-a_1-a_2+a_1a_2}{a_1(1-a_1)}z_1\right).
\end{align*}}%
Since after applying $\pi_1^{g_1} \circ \pi_2^{g_4} \circ \pi_1^{g_3}$ on $(z_1,z_2)$ the first coordinate remains unchanged, this proves the closedness.
\epf

In what follows, we will show that the statement of Theorem \ref{Thm:Zollproperty} is not true for non-regular Minkowski billiard trajectories. We give two simple examples, one within the triangle-, one within the parallelogram-configuration.

\begin{figure}[h!]
\centering
\def\svgwidth{380pt}
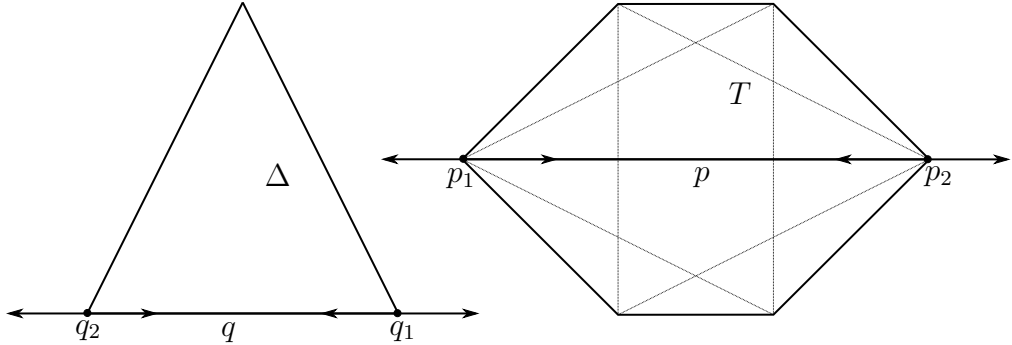
\caption{$\Delta\times T$ is an equality case of Viterbo's conjecture while $q=(q_1,q_2)$ is a closed $(\Delta,T)$-Minkowski billiard trajectory with closed dual billiard trajectory $p=(p_1,p_2)$ which is not $\ell_T$-minimizing.} 
\label{img:ce1triangle}
\end{figure}

Consider the triangle $\Delta$ given by the vertices
\beqq (-1,-1),\; (1,-1),\; (0,1)\eeqq
and the hexagon $T$ given by the vertices
\beqq \left(-\frac{1}{2},-1\right),\; \left(\frac{1}{2},-1\right),\; \left(\frac{3}{2},0\right),\; \left(\frac{1}{2},1\right),\; \left(-\frac{1}{2},1\right),\; \left(-\frac{3}{2},0\right)\eeqq
(see Figure \ref{img:ce1triangle}).
Then, $\Delta\times T$ is an equality case of Viterbo's conjecture, but if we consider the closed $(\Delta,T)$-Minkowski billiard trajectory $q=(q_1,q_2)$ given by the vertices $q_1=(1,-1)$ and $q_2=(-1,-1)$ and its corresponding dual billiard trajectory $p=(p_1,p_2)$ in $T$ given by the vertices $p_1=(-\frac{3}{2},0)$ and $p_2=(\frac{3}{2},0)$, then $q$ is not $\ell_T$-minimizing. One calculates
\beqq \ell_T(q)=\langle q_2-q_1,p_1-p_2\rangle = 6,\eeqq
where one notices that the minimal $\ell_T$-length of the closed $(\Delta,T)$-Minkowski billiard trajectories is $4$.

\begin{figure}[h!]
\centering
\def\svgwidth{380pt}
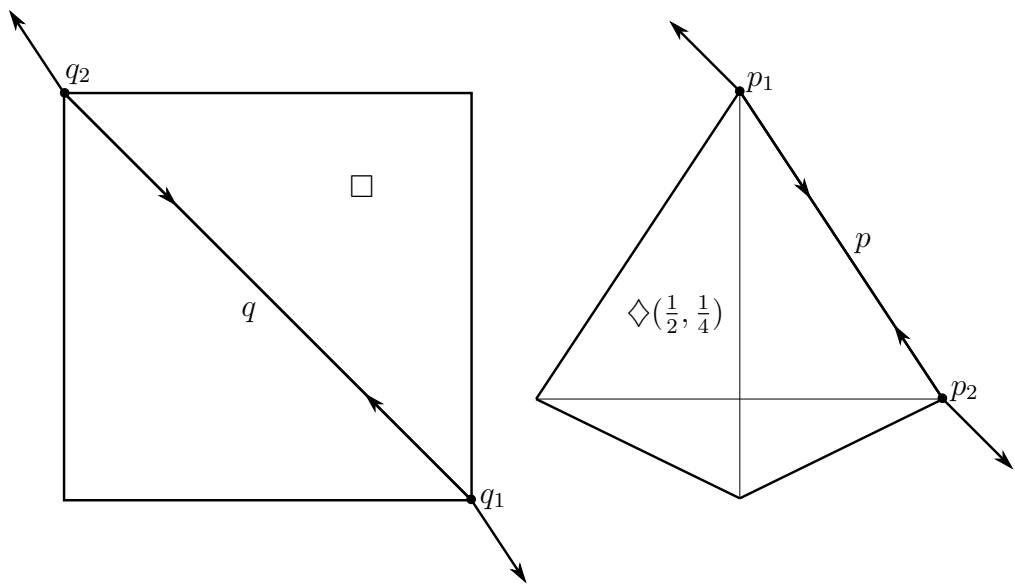
\caption{$\square\times \diamondsuit(a_1,a_2)$ is an equality case of Viterbo's conjecture while $q=(q_1,q_2)$ is a closed $(\square,\diamondsuit(a_1,a_2))$-Minkowski billiard trajectory with closed dual billiard trajectory $p=(p_1,p_2)$ which is not $\ell_{\diamondsuit(a_1,a_2)}$-minimizing.} 
\label{img:ce2square}
\end{figure}

Moreover, consider the square $\square$ given by the vertices
\beqq (1,-1),\; (1,1),\; (-1,1),\; (-1,-1)\eeqq
and the diamond $\diamondsuit(a_1,a_2)$ with $a_1=\frac{1}{2}$ and $a_2=\frac{1}{4}$ (see Figure \ref{img:ce2square}).
Then, $\square\times\diamondsuit$ is an equality cases of Viterbo's conjecture, but if we consider the closed $(\square,\diamondsuit(a_1,a_2))$-Minkowski billiard trajectory $q=(q_1,q_2)$ given by the vertices $q_1=(1,-1)$ and $q_2=(-1,1)$ and its corresponding dual billiard trajectory $p=(p_1,p_2)$ in $\diamondsuit(a_1,a_2)$ given by the vertices $p_1=(\frac{1}{2},1)$ and $p_2=(1,\frac{1}{4})$, then $q$ is not $\ell_{\diamondsuit(a_1,a_2)}$-minimizing. One calculates
\beqq \ell_{\diamondsuit(a_1,a_2)}(q)=\langle q_2-q_1,p_1-p_2\rangle = \frac{5}{2},\eeqq
where one notices that the minimal $\ell_{\diamondsuit(a_1,a_2)}$-length of the closed $(\square,\diamondsuit(a_1,a_2))$-Minkowski billiard trajectories is $2$.

\section*{Acknowledgement}
This research is supported by the SFB/TRR 191 'Symplectic Structures in Geometry, Algebra and Dynamics', funded by the \underline{German Research Foundation}, and was carried out under the supervision of Alberto Abbondandolo. The author is thankful to the supervisor's support and also would like to express his gratitude to Felix Schlenk for introducing him to the theory of symplectic embeddings. Moreover, the author thanks Alexey Balitskiy for his useful remarks on previous versions of this paper.


\medskip


\section*{Daniel Rudolf, Ruhr-Universit\"at Bochum, Fakult\"at f\"ur Mathematik, Universit\"atsstrasse 150, D-44801 Bochum, Germany.}
\center{E-mail address: daniel.rudolf@ruhr-uni-bochum.de}

\end{document}